\documentclass[11pt]{amsart}
\usepackage{amsmath,amsfonts,amssymb}
\usepackage[T2A]{fontenc}
\usepackage[cp1251]{inputenc}
\usepackage{amsfonts}
\usepackage[english]{babel}
\usepackage{amsmath,amsthm,amssymb}
\def\wh{\widehat}
\def\wt{\widetilde}
\def\R{\mathbb R}
\def\C{\mathbb C}
\def\Z{\mathbb Z}
\def\N{\mathbb N}
\def\A{\mathcal A}
\def\u{\mathbf u}
\def\e{\mathbf e}
\def\x{\mathbf x}
\def\v{\mathbf v}
\def\a{\mathbf a}
\def\b{\mathbf b}
\def\c{\mathbf c}
\def\k{\mathbf k}

\def\D{\mathbf D}
\def\FF{\mathbf F}
\def\H{\mathfrak H}
\def\NN{\mathfrak N}
\def\bxi{\boldsymbol \xi}
\def\bt{\boldsymbol \theta}
\def\bphi{\boldsymbol \phi}
\def\1{\bold 1}

\def\eps{\varepsilon}
\def\Dom{\mathrm{Dom}\,}
\def\Ker{\mathrm{Ker}\,}

\def\le{\leqslant}

\def\ge{\geqslant}

\begin{document}

\title[Homogenization of nonstationary Schr\"odinger type equations]{Homogenization of nonstationary Schr\"odinger type equations \\ with periodic coefficients}

\author{T.~A.~Suslina}


\address{St. Petersburg State University, Department of Physics, Ul'yanovskaya 3, Petrodvorets, St.~Petersburg, 198504, Russia}

\email{suslina@list.ru}


\begin{abstract}
In $L_2(\mathbb{R}^d;\C^n)$ we consider selfadjoint strongly elliptic second order differential operators ${\A}_\eps$
with periodic coefficients depending on ${\x}/\eps$. We study the behavior of the operator exponential
$\exp(-i \A_\eps \tau)$, $\tau \in \R$, for small $\eps$. Approximations for this exponential in the
$(H^s\to L_2)$-operator norm with a suitable $s$ are obtained. The results are applied to study the behavior of
the solution $\u_\eps$ of the Cauchy problem for the Schr\"odinger type equation $i \partial_\tau \u_\eps = \A_\eps \u_\eps$.
\end{abstract}

\thanks{The author would like to thank the Isaac Newton Institute for Mathematical Sciences, Cambridge, for support and hospitality
during the programme ``Periodic and Ergodic Spectral Problems'', where work on this paper was undertaken.}

\maketitle

\section*{Introduction}

The paper concerns  homogenization for periodic differential operators (DOs). A broad literature is devoted to homogenization problems;
first, we mention the books [BeLP], [BaPa], [ZhKO].

\subsection*{0.1. The class of operators}
We consider selfadjoint elliptic second order DOs in $L_2(\R^d;\C^n)$ admitting a factorization of the form
$$
\A = f(\x)^* b(\D)^* g(\x) b(\D) f(\x).
\eqno(0.1)
$$
Here $b(\D)=\sum_{l=1}^d b_l D_l$ is the $(m\times n)$-matrix first order DO with constant coefficients.
We assume that $m\ge n$ and that the symbol $b(\bxi)$ has maximal rank.
It is assumed that the matrix-valued functions  $g(\x)$ (of size $m\times m$) and $f(\x)$ (of size $n\times n$) are periodic with respect to
some lattice $\Gamma$ and such that
$$
g(\x)>0; \quad g,g^{-1} \in L_\infty; \quad f,f^{-1} \in L_\infty.
$$
It is convenient to start with a narrower class of operators
$$
\wh{\A} =  b(\D)^* g(\x) b(\D),
\eqno(0.2)
$$
corresponding to the case where $f=\1$.
Many operators of mathematical physics can be represented in the form (0.1) or (0.2);
the simplest example is the acoustics operator $\wh{\A} = - \hbox{\rm div}\, g(\x) \nabla = \D^* g(\x)\D$.
This and other examples are discussed in [BSu1] in detail.

Now we introduce the small parameter  $\eps >0$ and denote $\varphi^\eps(\x)= \varphi(\eps^{-1}\x)$ for any $\Gamma$-periodic function $\varphi(\x)$. 
Consider the operators
$$
\A_\eps = (f^\eps(\x))^* b(\D)^* g^\eps(\x) b(\D) f^\eps(\x),
\eqno(0.3)
$$
$$
\wh{\A}_\eps =  b(\D)^* g^\eps(\x) b(\D),
\eqno(0.4)
$$
whose coefficients oscillate rapidly as $\eps \to 0$.

\subsection*{0.2. Operator error estimates for elliptic and parabolic problems in $\R^d$}
 In a series of papers [BSu1--4] by  M.~Sh.~Birman and T.~A.~Suslina, an operator-theoretic approach to homogenization of elliptic equations in $\R^d$
 was suggested and developed. This approach is based on the scaling transformation, the Floquet-Bloch theory, and the analytic perturbation theory.

The homogenization problem for elliptic equations in $\R^d$ can be regarded as
 a problem of asymptotic description of the resolvent of $\A_\eps$ as $\eps \to 0$. For definiteness, let us talk about the simpler operators (0.4). 
 In [BSu1], it was shown that the resolvent  $(\wh{\A}_\eps +I)^{-1}$
converges to the resolvent \hbox{$(\wh{\A}^0+I)^{-1}$} in the  $L_2$-operator norm, as $\eps \to 0$.
Here $\wh{\A}^0= b(\D)^* g^0 b(\D)$ is the \textit{effective operator} with the constant \textit{effective matrix} $g^0$.
The formula for the effective matrix is well known in homogenization theory; in the case under consideration it is described below in
\S 8. In [BSu1], it was proved that
$$
\| (\wh{\A}_\eps +I)^{-1} - (\wh{\A}^0+I)^{-1} \|_{L_2(\R^d) \to L_2(\R^d)} \le C \eps.
\eqno(0.5)
$$
In [BSu2,3], a more accurate approximation of the resolvent of $\wh{\A}_\eps$ in the $L_2(\R^d;\C^n)$-operator norm with an error $O(\eps^2)$ was obtained,
and in [BSu4] an approximation of the same resolvent in the norm of operators acting from $L_2(\R^d;\C^n)$ to the Sobolev space $H^1(\R^d;\C^n)$
with an error $O(\eps)$ was found. In these approximations, some correction terms of first order (the \textit{correctors}) were taken into account.

Similarly, the homogenization problem for parabolic equations in $\R^d$ can be regarded as a problem of asymptotic description of the
operator exponential $\exp(- {\A}_\eps \tau)$ for $\tau >0$ and small $\eps$.
The operator-theoretic approach was applied to such problems in [Su1-3], [V], [VSu].
In [Su1,2], it was proved that
$$
\| \exp(-\wh{\A}_\eps \tau) - \exp(- \wh{\A}^0 \tau) \|_{L_2(\R^d) \to L_2(\R^d)} \le C \eps (\tau + \eps^2)^{-1/2}.
\eqno(0.6)
$$
In [V], a more accurate approximation of the operator $\exp(-\wh{\A}_\eps \tau)$
in the $L_2(\R^d;\C^n)$-operator norm with an error \hbox{$O(\eps^2)$} for  fixed $\tau$ was obtained, and in [Su3]
approximation of the same operator in the norm of operators acting from $L_2(\R^d;\C^n)$ to $H^1(\R^d;\C^n)$ with an error
{$O(\eps)$} for  fixed $\tau$ was proved. In these approximations, the first order correctors were taken into account.

Even more accurate approximations of the exponential and the resolvent of $\wh{\A}_\eps$
with the first and second correctors taken into account were found in [VSu].
In [BSu1--4], [Su1-3], [V], [VSu], similar (but more complicated) results were obtained also for more general operator (0.3); we will not dwell on this.

Estimates of the form (0.5), (0.6) are called \textit{operator error estimates} in homogenization theory.
A different approach to operator error estimates (the so called ``modified method of the first appproximation'')
was suggested by V.~V.~Zhikov. In [Zh], [ZhPas1], the acoustics operator and the operator of elasticity theory
(which have the form (0.4)) were studied;
approximations for the resolvents in the $(L_2 \to L_2)$-norm with an error $O(\eps)$
and in the $(L_2 \to H^1)$-norm with an error $O(\eps)$ were obtained.
In [ZhPas2],  estimate (0.6) was proved for the scalar elliptic operator $-{\rm div}\,g^\eps(\x) \nabla$.

\subsection*{0.3. Operator error estimates for nonstationary Schr\"odinger type and hyperbolic type equations}
So, in the case of elliptic and parabolic problems, the spectral approach to homogenization is developed in detail.
The situation with  homogenization of nonstationary  Schr\"odinger type and hyperbolic equations
is different. The paper [BSu5] is devoted to such problems. Again, we dwell on the results for the simpler operator (0.4).
In operator terms, the behavior of  the operator exponential $\exp(-i \tau \wh{\A}_\eps)$ and the operator cosine
$\cos (\tau \wh{\A}_\eps^{1/2})$ (where $\tau \in \R$) for small $\eps$ is studied.
For these operators it is impossible to obtain approximations in the $L_2(\R^d;\C^n)$-operator norm,
and we are forced to consider the norm of operators acting from the Sobolev space $H^s(\R^d;\C^n)$ (with appropriate $s$)
to $L_2(\R^d;\C^n)$. In [BSu5], the following estimates were proved:
$$
\| \exp(-i \tau \wh{\A}_\eps) - \exp(-i \tau \wh{\A}^0) \|_{H^3(\R^d) \to L_2(\R^d)}
\leqslant (C_1 + C_2 |\tau|) \eps,
\eqno(0.7)
$$
$$
\| \cos( \tau \wh{\A}_\eps^{1/2}) - \cos( \tau (\wh{\A}^0)^{1/2}) \|_{H^2(\R^d) \to L_2(\R^d)}
\leqslant (\wt{C}_1 + \wt{C}_2 |\tau|) \eps.
\eqno(0.8)
$$
By interpolation, we can also estimate the operator in (0.7) in the \hbox{$(H^s \to L_2)$}-norm by $O(\eps^{s/3})$ (where $0\le s \le 3$)
and the operator in (0.8) in the $(H^s \to L_2)$-norm by $O(\eps^{s/2})$ (where $0\le s \le 2$).
In [BSu5], approximations for the operator exponential and the operator cosine of the more general operator (0.3) were also obtained.
Note that for the operators  $\exp(-i \tau \wh{\A}_\eps)$ and
$\cos (\tau \wh{\A}_\eps^{1/2})$ there are no results concerning more accurate approximations with operator error estimates (and with some correctors taken into account).

The question about the sharpness of the resuts (0.7), (0.8) with respect to the type of operator norm (i.~e., the order of the Sobolev space)
remained open until now.

Let us explain the method of [BSu5]; we comment on the proof of estimate (0.7).
 Denote ${\mathcal H}_0 := -\Delta$.
Clearly, (0.7) is equivalent to
$$
\| \left( \exp(-i \tau \wh{\A}_\eps) - \exp(-i \tau \wh{\A}^0)\right) ({\mathcal H}_0 +I)^{-3/2} \|_{L_2(\R^d) \to L_2(\R^d)}
\leqslant (C_1 + C_2 |\tau|) \eps.
\eqno(0.9)
$$
In other words, in order to obtain an estimate in the $(L_2\to L_2)$-norm, we multiply the operator exponential by a
 ``smoothing factor'' \hbox{$({\mathcal H}_0 +I)^{-3/2}$}.
Next, the scaling transformation shows that (0.9) is equivalent to the estimate
$$
\begin{aligned}
\| \left( \exp(-i \eps^{-2} \tau \wh{\A}) - \exp(-i \eps^{-2} \tau \wh{\A}^0)\right) \eps^3 ({\mathcal H}_0 +\eps^2 I)^{-3/2} \|_{L_2(\R^d) \to L_2(\R^d)}
\\
\leqslant (C_1 + C_2 |\tau|) \eps.
\end{aligned}
\eqno(0.10)
$$
To prove (0.10), using the unitary Gelfand transformation, we expand $\wh{\A}$
in the direct integral of the operators $\wh{\A}(\k)$ acting in $L_2(\Omega;\C^n)$
(where $\Omega$ is the cell of the lattice $\Gamma$).
Here $\wh{\A}(\k)$ is given by the differential expression $b(\D+\k)^* g(\x)b(\D+\k)$ with periodic boundary conditions;
the spectrum of  $\wh{\A}(\k)$ is discrete. The family of operators $\wh{\A}(\k)$ is studied by means of the analytic perturbation theory
(with respect to the onedimensional parameter $t=|\k|$). For the operators $\wh{\A}(\k)$ the analog of estimate (0.10) is proved with
the constants independent of  $\k$. Then the inverse Gelfand transformation leads to (0.10).
A good deal of considerations in the study of the family $\wh{\A}(\k)$
is done in the abstract operator-theoretic setting.

\subsection*{0.4. Main results of the paper}
In the present paper, we study the behavior of the operator exponential $\exp(-i \tau {\A}_\eps)$ for small $\eps$, and next
we apply the results to study the behavior of the solution of the nonstationary Schr\"odinger type equation.
On the one hand, we confirm the sharpness of estimate (0.7) in the following sense.
We find a condition on the operator, under which the estimate
$$
\| \exp(-i \tau \wh{\A}_\eps) - \exp(-i \tau \wh{\A}^0) \|_{H^s(\R^d) \to L_2(\R^d)}
\leqslant C(\tau) \eps
$$
is false if $s<3$.
It is easy to formulate this condition in the spectral terms. We consider the operator family $\wh{\A}(\k)$ and put $\k = t \bt$, $t=|\k|$, $\bt \in {\mathbb S}^{d-1}$.
This family is analytic with respect to the parameter~$t$. For $t=0$ the point $\lambda_0=0$ is an eigenvalue of multiplicity $n$ of the ``unperturbed''
operator $\wh{\A}(0)$. Then for small $t$ there exist the real-analytic branches of the eigenvalues and the eigenvectors of  $\wh{\A}(t\bt)$.
For small $t$ the eigenvalues $\lambda_l(t,\bt)$, $l=1,\dots,n,$ admit the convergent power series expansions
$$
\lambda_l(t,\bt) = \gamma_l(\bt) t^2 + \mu_l(\bt) t^3 + \dots,\quad l=1,\dots,n,
\eqno(0.11)
$$
where $\gamma_l(\bt)>0$ and $\mu_l(\bt)\in \R$.
The condition is that $\mu_l(\bt_0) \ne 0$ for at least one $l$ and at least one point $\bt_0 \in {\mathbb S}^{d-1}$.
Examples of the operators satisfying this condition are provided; in particular, one example is
of the form $-{\rm div}\,g^\eps(\x) \nabla$, where $g(\x)$ is  Hermitian matrix with complex entries.

On the other hand, we distinguish conditions on the operator under which it is possible to improve the result and obtain the estimate
$$
\| \exp(-i \tau \wh{\A}_\eps) - \exp(-i \tau \wh{\A}^0) \|_{H^2(\R^d) \to L_2(\R^d)}
\leqslant (\check{C}_1+ \check{C}_2|\tau|) \eps.
\eqno(0.12)
$$
In the case where  $n=1$, for (0.12) it suffices that the coefficient
$\mu(\bt)=\mu_1(\bt)$ in (0.11) is identically zero. In particular, this is the case for the operator
$-{\rm div}\,g^\eps(\x) \nabla$, where $g(\x)$ is symmetric matrix with real entries.
In the matrix case (i.~e., for $n\ge 2$), besides the condition that all the coefficients $\mu_l(\bt)$
in (0.11) are equal to zero, we impose one more condition in terms of the coefficients
$\gamma_l(\bt)$, $l=1,\dots,n$. The simplest version of this condition is that the branches $\gamma_l(\bt)$ must not intersect:
for each pair $j\ne l$ either $\gamma_j(\bt)$ and $\gamma_l(\bt)$ are separated from each other or
they coincide identically for all $\bt\in {\mathbb S}^{d-1}$.

It turns out that for more general operator (0.3) it is convenient to study the operator exponential sandwiched between appropriate rapidly oscillating factors.
Namely, we study the operator $f^\eps e^{-i \tau \A_\eps } (f^\eps)^{-1}$ and obtain analogs of the results described above for this operator.

Next, we apply the results given in operator terms to study the behavior of the solution
$\u_\eps(\x,\tau)$, $\x\in \R^d$, $\tau \in \R$, of the following problem
$$
i \partial_\tau \u_\eps(\x,\tau) = (\wh{\A}_\eps \u_\eps)(\x,\tau) + \FF(\x,\tau),
\quad \u_\eps(\x,0)= \bphi(\x).
$$
A more general problem with the operator $\A_\eps$ is also studied.

We apply the general results to specific equations of mathematical physics.
In particular, we consider the  nonstationary Schr\"odinger equation
$i \partial_\tau u_\eps = -{\rm div} g^\eps(\x) \nabla u_\eps(\x,\tau)  + \eps^{-2} V^\eps(\x)  u_\eps(\x,\tau)$
with the singular potential $\eps^{-2} V^\eps$, and also the twodimensional Pauli equation with  the singular magnetic potential.

   Similar results have been obtained by the author jointly with M.~A.~Dorodnyi [DSu] for homogenization of the hyperbolic equations with
   rapidly oscillating coefficients.

\subsection*{0.5. Method}
The results are obtained by further development of the operator-theoretic approach.
We follow the plan described above in Subsection~0.3.
Considerations are based on the abstract operator-theoretic scheme. Let us dwell on this. In the abstract setting, we study the family of operators $A(t)=X(t)^* X(t)$ acting
in some Hilbert space $\H$. Here $X(t)=X_0 + tX_1$. (This family is modelling the operator family $\A(\k)= \A(t\bt)$,
but the parameter $\bt$ is absent in the abstract setting.) It is assumed that
the point $\lambda_0=0$ is an eigenvalue of $A(0)$ of finite multiplicity $n$.
Then for $|t|\le t^0$ the perturbed operator $A(t)$ has exactly $n$ eigenvalues (counted with multiplicities) on the interval $[0,\delta]$
(here $\delta$ and $t^0$ are controlled explicitly).
These eigenvalues and the corresponding eigenvectors are real-analytic functions of $t$.
The coefficients of the corresponding power series expansions for the eigenvalues and the eigenvectors are called \textit{threshold characteristics} of the operator $A(t)$.
We distinguish a finite rank operator $S$ (the so called  \textit{spectral germ} of the operator family $A(t)$) which acts in the space $\NN = \Ker A(0)$.
The spectral germ (see the definition in Subsection 1.2 below) contains the information about the threshold characteristics of principal order.
 Let $F(t)$ be the spectral projection of the operator $A(t)$ for the interval $[0,\delta]$.
 We rely on the threshold approximations for the projection $F(t)$ and for the operator $A(t)F(t)$ obtained in \hbox{[BSu1, Chapter~1]} and [BSu2].
  Note that in  [BSu5] the threshold approximations of principal order from [BSu1] were applied: $F(t)$ was approximated by the projection $P$ onto the subspace $\NN$,
  and the operator   $A(t)F(t)$ was approximated by $t^2SP$.
   It turns out that, in order to obtain more subtle results described above in Subsection 0.4, we need to use
  more accurate threshold approximations obtained in [BSu2]. Moreover,
  we need to divide the eigenvalues of $A(t)$ into clusters and find more detailed
  threshold approximations associated with this division (see \S 2).

 In terms of the spectral germ, it is possible to approximate the operator exponential $\exp(-i \eps^{-2} \tau A(t))$
 multiplied by an appropriate ``smoothing factor''.
  Application of the abstract results leads to the required estimates for differential operators.
  However, at this step additional difficulties arise. They concern the improvement of the results in the case where
  all the coefficients $\mu_l(\bt)$ are equal to zero. These difficulties are related to the fact that
  in the general (matrix) case we are not always able to make our constructions and estimates uniform in $\bt$, and we are forced to
  impose the additional assumptions of isolation of the branches $\gamma_l(\bt)$, $l=1,\dots,n$.

\subsection*{0.6. The scalar case}
For the case of the scalar elliptic operator $\wh{\A}_\eps = -{\rm div} g^\eps(\x) \nabla$
in $L_2(\R^d)$ main results of the paper can be obtained by a simpler method based on analyticity of the operator family
$\wh{\A}(\k)$ with respect to the multidimensional parameter $\k$ and variational arguments.
This is the subject of the  joint paper by I.~Kamotski and the author [KamSu].

\subsection*{0.7. Plan of the paper} The paper consists of three chapters. Chapter~1 (\S 1--5) contains the necessary operator-theoretic material.
In Chapter~2 (\S6--12), the periodic differential operators of the form (0.1), (0.2) are studied.
In \S 6 we describe the class of operators and introduce the Gelfand transformation. \S 7 is devoted to the direct integral expansion for periodic operators of the form (0.1);
the corresponding family of operators $\A(\k)$ acting in $L_2(\Omega;\C^n)$ is incorporated in the framework of the abstract scheme.
In \S 8 we describe the effective characteristics for the operator (0.2).
In \S 9, using the abstract results, we obtain approximation of the smoothed operator exponential
$\exp(- i \eps^{-2} \tau \wh{\A}(\k))$. The operator $\A(\k)$ is considered in \S 10.
In \S 11, with the help of the abstract results we find approximation of the smoothed sandwiched exponential
 $\exp(- i \eps^{-2} \tau {\A}(\k))$. Next, in \S 12 we return to the operators  (0.1), (0.2) acting in $L_2(\R^d;\C^n)$;
 applying the results of \S 9 and \S 11, we obtain approximations of the smoothed operator
 $\exp(- i \eps^{-2} \tau \wh{\A})$ and of the smoothed sandwiched operator
 $\exp(- i \eps^{-2} \tau {\A})$. Chapter~3 (\S 13--16) is devoted to homogenization problems.
In \S 13, by the scaling transformation, the results of Chapter~2 imply main results of the paper: approximations for the exponential  $\exp(- i \tau \wh{\A}_\eps)$
  and for the sandwiched exponential $\exp(- i \tau {\A}_\eps)$ in the $(H^s \to L_2)$-norm.
  In \S 14, the results are applied to study the behavior of the solution of the Cauchy problem for the nonstationary
  Schr\"odinger type equation. The last \S 15 and \S 16 are devoted to applications of the general results to the particular equations, namely, to
  the nonstationary Schr\"odinger equation and the two-dimensional nonstationary Pauli equation.

\subsection*{0.8. Notation}
Let $\H$ and $\H_*$ be complex separable Hilbert spaces.
The symbols $(\cdot,\cdot)_{\H}$ and $\|\cdot \|_{\H}$ stand for the inner product and the norm in $\H$, respectively;
the symbol $\|\cdot\|_{\H \to \H_*}$ stands for the norm of a linear continuous operator from $\H$ to $\H_*$.
Sometimes we omit the indices. By $I = I_\H$ we denote the identity operator in $\H$.
If $\NN$ is a subspace in $\H$, then $\NN^\perp := \H \ominus \NN$.
If $P$ is the orthogonal projection of $\H$ onto $\NN$, then $P^\perp$ is the orthogonal projection onto $\NN^\perp$.
If $A: \H \to \H_*$ is a linear operator, then $\Dom A$ stands for its domain and $\Ker A$ stands for its kernel.

The symbols $\langle\cdot, \cdot \rangle$ and $|\cdot|$ denote the inner product and the norm in $\C^n$;
$\1=\1_n$ is the unit $(n \times n)$-matrix. If $a$ is an $(n\times n)$-matrix, then the symbol $|a|$
denotes the norm of $a$ viewed as a linear operator in $\C^n$.
 Next, we use the notation $\x = (x_1,\dots, x_d) \in \R^d$, $i D_j = \partial_j = \partial / \partial x_j$,
$j=1,\dots,d$, $\D = -i \nabla = (D_1,\dots, D_d)$.

The $L_p$-classes of $\C^n$-valued functions in a domain ${\mathcal O} \subset \R^d$ are denoted by $L_p({\mathcal O}; \C^n)$, $1 \le p \le \infty$.
The Sobolev classes of $\C^n$-valued functions in a domain $\mathcal O$ of order $s$ and integrability index $p$
are denoted by  $W_p^s({\mathcal O}; \C^n)$. For $p=2$, we denote this space by $H^s({\mathcal O}; \C^n)$.
If $n=1$, we write simply $L_p({\mathcal O})$, $H^s({\mathcal O})$, but sometimes we use such abbreviated notation also for the spaces
of vector-valued or matrix-valued functions.

By $c, C, {\mathcal C}, {\mathfrak C}$ (possibly, with indices and marks) we denote various constants in estimates.

\section*{Chapter~1. Abstract operator-theoretic scheme}

\section*{§1.  Preliminaries}
Our approach to homogenization problems is based on the abstract operator-theoretic scheme.

\subsection*{1.1. The operators $X(t)$ and $A(t)$}
Let $\H$ and $\H_*$ be complex separable Hilbert spaces. Suppose that
$X_0: \H \to \H_*$ is a densely defined and closed operator, and $X_1:\H \to \H_*$ is a bounded operator.
Then the operator $X(t) := X_0 + tX_1$, $t \in \R$, is closed on the domain $\Dom X(t) = \Dom X_0$.
In the abstract setting, the family of selfadjoint (and nonnegative) operators
$$
A(t) := X(t)^* X(t),\quad t \in \R,
\eqno(1.1)
$$
in $\H$ is our main object. The operator (1.1) is generated by the closed quadratic form $\| X(t)u\|^2_{\H_*}$,
$u \in \Dom X_0$. Denote $A(0) = X_0^*X_0 =: A_0$ and
$\NN := \Ker A_0 = \Ker X_0$. We impose the following condition.

\smallskip\noindent
\textbf{Condition 1.1.} \textit{The point} $\lambda_0 = 0$ \textit{is an isolated point of the spectrum of}
$A_0$, \textit{and} $0 < n := \dim \NN < \infty$.

\smallskip
By $d^0$ we denote the \textit{distance from the point $\lambda_0 = 0$ to the rest of the spectrum of} $A_0$.
We put $\NN_*:= \Ker X_0^*$, $n_* := \dim \NN_*$, and \textit{assume that} $n \le n_* \le \infty$.
Let $P$ and $P_*$  be the orthogonal projections of $\H$ onto $\NN$ and of $\H_*$ onto $\NN_*$, respectively.

The operator family $A(t)$ has been studied in [BSu1, Chapter~1; BSu2; BSu4, Chapter~1] in detail.

Denote by $F(t;[a,b])$ the spectral projection of $A(t)$ for the interval $[a,b]$, and put
${\mathfrak F}(t;[a,b]):=F(t; [a,b]) \H$. \textit{We fix a number $\delta>0$ such that} $8 \delta < d^0$.
We often write $F(t)$ in place of $F(t;[0,\delta])$ and ${\mathfrak F}(t)$ in place of ${\mathfrak F}(t;[0,\delta])$.
Next, we choose a number $t^0 >0$ such that
$$
t^0 \le \delta^{1/2} \|X_1 \|^{-1}.
\eqno(1.2)
$$
According to [BSu1, Chapter~1, (1.3)],
$$
F(t;[0,\delta]) = F(t;[0,3\delta]),\quad  {\rm rank}\,F(t;[0,\delta])=n,\quad |t|\le t^0.
$$

\subsection*{1.2. The operators $Z$, $R$, and $S$}
Now we introduce some operators appearing in the analytic perturbation theory considerations; see [BSu1, Chapter~1, \S 1; BSu2, \S 1].

Let $\omega \in \NN$, and let $\psi = \psi(\omega) \in \Dom X_0 \cap \NN^\perp$ be a (weak) solution of the equation
$$
X_0^*(X_0 \psi + X_1 \omega) = 0.
$$
We define a bounded operator $Z : \H \to \H$ by the following relation
$$
Zu = \psi(Pu), \ \ u\in \H.
$$
Note that $Z$  takes $\NN$ to $\NN^\perp$ and $\NN^\perp$ to $\{0\}$.

Let $R : \NN \to \NN_*$ be the operator defined by
$$
R \omega = X_0 \psi(\omega) +X_1 \omega = (X_0 Z + X_1)\omega, \ \ \omega \in \NN.
$$
Another representation for $R$ is given by $R = P_* X_1\vert_{\NN}$.

According to [BSu1, Chapter~1, Subsection~1.3], the operator  $S:=R^*R: \NN \to \NN$
is called the \textit{spectral germ of the operator family $A(t)$ at} $t=0$.
The germ $S$ can be represented as
$$
S= P X_1^* P_* X_1\vert_{\NN}.
\eqno(1.3)
$$
The spectral germ is said to be \textit{nondegenerate} if ${\rm Ker}\, S = \{ 0\}$.

We have
$$
\|Z\| \le (8\delta)^{-1/2} \|X_1\|,\quad
\|R\| \le \|X_1\|,\quad
\|S\| \le \|X_1\|^2.
\eqno(1.4)
$$

\subsection*{1.3. The analytic branches of eigenvalues and eigenvectors  of $A(t)$} According to the general analytic perturbation theory (see [Ka]), for $|t|\le t^0$ there exist real-analytic functions $\lambda_l(t)$
(the branches of the eigenvalues) and real-analytic $\H$-valued functions $\varphi_l(t)$ (the branches of the eigenvectors) such that
$$
A(t)\varphi_l(t) = \lambda_l(t)\varphi_l(t),\quad l=1,\dots,n,\quad |t|\le t^0,
\eqno(1.5)
$$
and the $\varphi_l(t)$, $l=1,\dots,n$, form an \textit{orthonormal basis} in ${\mathfrak F}(t)$.
Moreover, for $|t|\le t_*$, where $0< t_* \le t^0$ \textit{is sufficiently small}, we have the following convergent power series expansions:
$$
\lambda_l(t) = \gamma_l t^2 + \mu_l t^3 +  \dots,\quad \gamma_l \ge 0,\quad \mu_l \in \R,\quad l=1,\dots,n,
\eqno(1.6)
$$
$$
\varphi_l(t) = \omega_l + t \psi_l^{(1)}  + \dots,\quad l=1,\dots,n.
\eqno(1.7)
$$
The elements $\omega_l= \varphi_l(0)$, $l=1,\dots,n,$ form an orthonormal basis in $\NN$.

Substituting expansions (1.6), (1.7) in (1.5) and comparing the coefficients of  powers $t$ and $t^2$, we arrive at the relations
$$
\widetilde{\omega}_l = \psi_l^{(1)} - Z \omega_l \in \NN,\quad l=1,\dots,n,
\eqno(1.8)
$$
$$
S \omega_l = \gamma_l \omega_l,\quad l=1,\dots,n.
\eqno(1.9)
$$
(Cf. [BSu1, Chapter~1, \S 1]; BSu2, \S 1].)
Thus, the \textit{numbers $\gamma_l$ and the elements $\omega_l$ defined by}  (1.5)--(1.7) \textit{are eigenvalues and eigenvectors of the germ} $S$.
We have
$$
P= \sum_{l=1}^n (\cdot, \omega_l) \omega_l,
\eqno(1.10)
$$
$$
S P= \sum_{l=1}^n \gamma_l (\cdot, \omega_l) \omega_l.
\eqno(1.11)
$$

\smallskip\noindent\textbf{Remark 1.2.}
Relations (1.9) give another (spectral) definition of the germ.
These relations show that the germ of $A(t)$ at $t=0$ does not depend on the (possibly, nonunique) choice of factorization
(1.1). At the same time,  invariant representation (1.3) shows that the germ does not depend on the (possibly, nonunique)
choice of the (analytic) basis $\{\varphi_l(t)\}$ in ${\mathfrak F}(t)$.
If all the eigenvalues $\gamma_l$ of $S$ are simple, then the ``embrios'' $\omega_l$ in (1.7) are defined by $S$ uniquely
(up to  phase factors).
If there are multiple eigenvalues among $\gamma_l$, then the knowledge of $S$ is not sufficient for
determining the elements $\omega_l$.

   \smallskip
   Note that
   $$
   (\widetilde{\omega}_l , \omega_j) + (\omega_l, \widetilde{\omega}_j)=0,\quad l,j =1,\dots,n,
   \eqno(1.12)
   $$
   which follows from the relations $(\varphi_l(t),\varphi_j(t))=\delta_{lj}$
by substituting  (1.7), comparing the coefficients of power $t$, and taking (1.8) into account.

\smallskip\noindent\textbf{1.4. Threshold approximations.}
The spectral projection $F(t)$ and the operator $A(t)F(t)$ are real-analytic operator-valued functions for $|t|\le t^0$.
We have
$$
F(t) = \sum_{l=1}^n (\cdot, \varphi_l(t)) \varphi_l(t), \quad |t|\le t^0,
$$
$$
A(t)F(t) = \sum_{l=1}^n \lambda_l(t) (\cdot, \varphi_l(t)) \varphi_l(t), \quad |t|\le t^0.
$$
Together with (1.6), (1.7), (1.10), and (1.11) this yields the power series expansions
$F(t) = P + tF_1+\dots$ and $A(t)F(t)= t^2 SP+ t^3 K +\dots$, convergent for $|t|\le t_*$.
However, for our purposes not expansions, but
approximations (with one or several first terms)
 with error estimates on the whole interval $|t|\le t^0$ are needed.

The following statement was obtained in [BSu1, Chapter~1, Theorems~4.1~and~4.3]). In what follows, we agree
to denote by $\beta_j$ various  absolute constants (which can be controlled explicitly) assuming that $\beta_j \ge 1$.

\smallskip\noindent\textbf{Theorem 1.3.}
\textit{Suppose that the assumptions of Subsection} 1.1 \textit{are satisfied. Then we have}
$$
\|F(t) - P\| \le C_1|t|,\quad |t|\le t^0,
\eqno(1.13)
$$
$$
\|A(t)F(t) -t^2 S  P\| \le C_2|t|^3,\quad |t|\le t^0.
\eqno(1.14)
$$
\textit{Here $t^0$ is subject to} (1.2), \textit{and the constants $C_1$, $C_2$ are given by
$$
C_1= \beta_1 \delta^{-1/2} \|X_1\|,\quad C_2 = \beta_2 \delta^{-1/2} \|X_1\|^3.
\eqno(1.15)
$$
}

We also need a more precise approximation for the operator $A(t)F(t)$ obtained in  [BSu2,~Theorem~4.1].

\smallskip\noindent\textbf{Theorem 1.4.}
\textit{Suppose that the assumptions of Subsection} 1.1 \textit{are satisfied. Then for $|t|\le t^0$ we have
$$
A(t) F(t) = t^2 SP + t^3 K + \Psi(t),
\eqno(1.16)
$$
and
$$
\| \Psi(t)\| \le C_3 t^4,\quad |t|\le t^0.
\eqno(1.17)
$$
The operator  $K$ is represented as
$$
K = K_0 + N,\quad N=N_0 + N_*,
$$
where $K_0$ takes $\NN$ to $\NN^\perp$ and $\NN^\perp$ to $\NN$, while $N_0$ and $ N_*$ take $\NN$ to itself and $\NN^\perp$ to
$\{ 0\}$. In terms of the power series coefficients, these operators are given by
$$
K_0 = \sum_{l=1}^n \gamma_l \left(  (\cdot, Z \omega_l) \omega_l + (\cdot, \omega_l) Z \omega_l \right),
$$
$$
N_0 = \sum_{l=1}^n \mu_l  (\cdot,  \omega_l) \omega_l,
\eqno(1.18)
$$
$$
N_* = \sum_{l=1}^n \gamma_l  \left((\cdot,  \widetilde{\omega}_l) \omega_l + (\cdot,\omega_l) \widetilde{\omega}_l \right).
\eqno(1.19)
$$
In the invariant terms, we have $K_0 = ZSP + SP Z^*$ and $N = Z^* X_1^* RP + (RP)^* X_1Z$.
The constant in} (1.17) \textit{is given by $C_3 = \beta_3 \delta^{-1} \|X_1\|^4$.
We have}
$$
\|K_0\|\le (2\delta)^{-1/2} \|X_1\|^3,\quad
\|N\|\le (2\delta)^{-1/2} \|X_1\|^3.
\eqno(1.20)
$$
\smallskip

Note that
$$
P^\perp K P = ZSP,\quad PKP=N.
\eqno(1.21)
$$

\smallskip\noindent
\textbf{Remark 1.5.}
$1^\circ$. If $Z=0$, then $K_0=0$, $N=0$, and $K=0$.

\noindent $2^\circ$. In the basis $\{\omega_l\}$, the operators $N$, $N_0$, and $N_*$ (restricted to $\NN$)
are represented by $(n\times n)$-matrices. The operator $N_0$ is diagonal:
$$
 (N_0 \omega_j, \omega_k) = \mu_j \delta_{jk}, \quad j,k=1,\dots,n.
$$
The matrix entries of $N_*$ are given by
$$
(N_* \omega_j, \omega_k) = \gamma_k (\omega_j,\widetilde{\omega}_k) +
\gamma_j (\widetilde{\omega}_j, \omega_k)= (\gamma_j - \gamma_k) (\widetilde{\omega}_j, \omega_k),
  \quad j,k=1,\dots,n.
$$
Here we have taken (1.12) into account.  It is seen that the diagonal elements of $N_*$ are equal to zero:
$(N_* \omega_j, \omega_j)=0$, $j=1,\dots,n$. Moreover,
$$
(N_* \omega_j, \omega_k)=0 \quad \text{if}\ \gamma_j=\gamma_k.
$$
Thus, in the basis $\{\omega_l\}$, the diagonal part of $N$ coincides with $N_0$, and the off-diagonal part coincides with  $N_*$.
Moreover, in the case where some eigenvalues of $S$ are multiple,
the off-diagonal elements of $N$ corresponding to different eigenvectors with the same eigenvalue are equal to zero.

\noindent $3^\circ$.
If $n=1$, then $N_*=0$, i.~e., $N=N_0$.

\smallskip\noindent\textbf{1.5. The nondegeneracy condition.}
Below we impose the following additional condition on the operator $A(t)$.

\smallskip\noindent\textbf{Condition 1.6.} \textit{There exists a constant $c_*>0$ such that}
$$
A(t) \ge c_* t^2 I,\quad |t|\le t^0.
\eqno(1.22)
$$

\smallskip
From (1.22) it follows that $\lambda_l(t) \ge c_* t^2$, $l=1,\dots,n$, for $|t|\le t^0$.
By (1.6), this implies
$$
\gamma_l \ge c_* >0,\quad l=1,\dots,n,
\eqno(1.23)
$$
i.~e., the germ $S$ is nondegenerate.

\section*{§2. The clusters of eigenvalues of $A(t)$}

This section concerns the case where $n \ge 2$.

\subsection*{2.1. Renumbering of eigenvalues} Suppose that Condition~1.6 is satisfied.
Then the spectrum of the operator $SP$ consists of the eigenvalue $\lambda_0=0$
(with the eigenspace $\NN^\perp$) and the eigenvalues $\gamma_1,\dots, \gamma_n$
 satisfying (1.23). Now it is convenient to change the notation tracing
the multiplicities of the eigenvalues.
Let $p$ be the number of different eigenvalues of the germ.
We enumerate these eigenvalues in the increasing order
and denote them by $\gamma_j^\circ$, $j=1,\dots,p$.
Their multiplicities are denoted by $k_1,\dots, k_p$ (obviously, $k_1+\dots+k_p =n$).
Then, in the previous notation,
$$
\gamma_1=\dots=\gamma_{k_1} < \gamma_{k_1+1}=\dots = \gamma_{k_1+k_2} < \dots
< \gamma_{n- k_p+1}=\dots = \gamma_n.
$$
We have
$\gamma_1^\circ=\gamma_1=\dots=\gamma_{k_1}$,
$\gamma_2^\circ=\gamma_{k_1+1}=\dots=\gamma_{k_1+k_2}$, etc.
Let $\NN_j = {\rm Ker}\, (S - \gamma_j^\circ I_{\NN})$, $j=1,\dots,p$. Then
$$
\NN = \sum_{j=1}^p \oplus \NN_j.
$$
Let $P_j$ be the orthogonal projection of $\H$ onto $\NN_j$. Then
$$
P = \sum_{j=1}^p P_ j, \quad P_j P_l =0 \quad \text{for}\ j\ne l.
\eqno(2.1)
$$

We also change the notation for the eigenvectors of the germ
(which are the ``embrios'' in (1.7)) dividing them in $p$ parts so that
$\omega_1^{(j)},\dots, \omega^{(j)}_{k_j}$ correspond to the eigenvalue
 $\gamma_j^\circ$ and form an orthonormal basis in $\NN_j$.
 (In the previous notation, these are $\omega_{k_1+\dots+k_{j-1}+1}, \dots, \omega_{k_1+\dots+k_j}$.)

We also change the notation for the analytic branches of the eigenvalues and the eigenvectors of $A(t)$.
The eigenvalue and the eigenvector whose expansions (1.6) and (1.7) start with the terms $\gamma_j^\circ t^2$ and $\omega^{(j)}_q$
are denoted by $\lambda^{(j)}_q(t)$ and $\varphi^{(j)}_q(t)$, respectively.
 For $|t|\le t_*$ we have
 $$
 \lambda^{(j)}_q(t)= \gamma_j^\circ t^2 + \mu^{(j)}_q t^3 +  \dots, \quad q=1,\dots,k_j,
 $$
 $$
 \varphi_q^{(j)}(t) = \omega^{(j)}_q + t \psi^{(j)}_q + \dots,\quad q=1,\dots, k_j.
 $$

 \smallskip\noindent\textbf{2.2. Refinement of threshold approximations.}
 Suppose that Condition 1.6 is satisfied. For $|t|\le t^0$ we consider
 the following bounded selfadjoint operator in $\H$:
$$
{\mathfrak A}(t) =  \begin{cases}
 t^{-2} A(t)F(t), & t\ne 0,\cr
SP, & t=0.
\end{cases}
$$
We apply the spectral perturbation theory arguments,
treating ${\mathfrak A}(t)$ as a perturbation of the operator $SP$. By (1.14),
$$
\| {\mathfrak A}(t) - SP \| \le C_2 |t|,\quad |t| \le t^0.
\eqno(2.2)
$$

The results about $A(t)$ (see \S 1) show that for $|t|\le t^0$
the point $\lambda_0=0$ is an eigenvalue of the perturbed operator ${\mathfrak A}(t)$
(with the eigenspace ${\mathfrak F}(t)^\perp$), and ${\mathfrak A}(t)$ has
positive eigenvalues of total multiplicity $n$.
We divide these positive eigenvalues in $p$ clusters
which for small $|t|$ are located near the eigenvalues $\gamma_1^\circ,\dots, \gamma_p^\circ$
of the unperturbed operator. Clearly, the $j$-th  cluster
consists of the eigenvalues $\nu^{(j)}_q(t) = t^{-2} \lambda^{(j)}_q(t)$ of ${\mathfrak A}(t)$, $q=1,\dots,k_j$,
since $\nu^{(j)}_q(t)$ are continuous (and even analytic) in $t$
and $\nu^{(j)}_q(0)=\gamma_j^\circ$.
The corresponding orthonormal eigenvectors are $\varphi^{(j)}_q(t)$, $q=1,\dots,k_j$.

For sufficiently small $|t|$ these clusters are separated from each other.
However, it will be more convenient for our purposes, for each pair of indices $j\ne l$, to divide
the clusters in two parts separated  from each other and such that one part contains
the $j$-th cluster and another part contains the  $l$-th cluster.
For each pair of indices $(j,l)$,  $1\le j,  l \le p$, $j \ne l$, we denote
  $$
 c^\circ_{jl} := \min\{ c_*,  n^{-1}| \gamma_{l}^\circ - \gamma^\circ_{j}| \}.
 \eqno(2.3)
 $$
Clearly, there exists a number $i_0 = i_0(j,l)$, where $j \le i_0 \le l-1$ if $j<l$
and $l \le i_0 \le j-1$ if $l<j$, such that $\gamma^\circ_{i_0+1} - \gamma_{i_0}^\circ \ge c^\circ_{jl}$.
It means that on the interval between $\gamma_{j}^\circ$ and $\gamma_{l}^\circ$
there is a gap in the spectrum of $S$ of length at least $c^\circ_{jl}$.
If such $i_0$ is not unique, we agree to take the minimal possible $i_0$ (for definiteness).

  We choose a number $t^{00}_{jl}\le t^0$ such that  (see (1.15))
$$
t^{00}_{jl} \le  (4 C_2)^{-1} c^\circ_{jl} = (4\beta_2)^{-1} \delta^{1/2} \| X_1 \|^{-3} c^\circ_{jl}.
\eqno(2.4)
$$

By (2.2), for $|t| \le t^{00}_{jl}$ we have $\| {\mathfrak A}(t) - SP\| \le c^\circ_{jl}/4$.
Hence, the segments $[\gamma_1^\circ - c_{jl}^\circ/4, \gamma^\circ_{i_0} + c_{jl}^\circ/4]$ and
\hbox{$[\gamma_{i_0+1}^\circ - c_{jl}^\circ/4, \gamma^\circ_{p} + c_{jl}^\circ/4]$}
are disjoint, and  the distance between them is at least $c_{jl}^\circ/2$.
Consequently, for $|t| \le t^{00}_{jl}$ the perturbed operator ${\mathfrak A}(t)$
has exactly $k_1+\dots + k_{i_0}$ eigenvalues (counted with multiplicities)
in the segment \hbox{$[\gamma_1^\circ - c_{jl}^\circ/4, \gamma^\circ_{i_0} + c_{jl}^\circ/4]$}.
These are $\nu^{(1)}_1(t),\dots, \nu^{(1)}_{k_1}(t); \dots; \nu^{(i_0)}_1(t),\dots, \nu^{(i_0)}_{k_{i_0}}(t)$.
We denote the corresponding eigenspace by ${\mathfrak F}_{jl}^{(1)}(t)$;
for $t\ne 0$ it coincides with the eigenspace \hbox{${\mathfrak F}(t;[0, (\gamma^\circ_{i_0} + c^\circ_{jl}/4)t^2])$}
 of $A(t)$. The elements  $\varphi^{(1)}_{1}(t), \dots, \varphi^{(1)}_{k_1}(t);\dots; \varphi^{(i_0)}_{1}(t), \dots, \varphi^{(i_0)}_{k_{i_0}}(t)$ form an orthonormal basis in ${\mathfrak F}_{jl}^{(1)}(t)$.
Similarly, for $|t| \le t^{00}_{jl}$ the perturbed operator ${\mathfrak A}(t)$
has exactly $k_{i_0+1}+\dots + k_{p}$ eigenvalues (counted with multiplicities)
in the segment  $[\gamma_{i_0+1}^\circ - c_{jl}^\circ/4, \gamma^\circ_{p} + c_{jl}^\circ/4]$.
These are $\nu^{(i_0+1)}_1(t),\dots, \nu^{(i_0+1)}_{k_{i_0+1}}(t); \dots; \nu^{(p)}_1(t),\dots, \nu^{(p)}_{k_{p}}(t)$.
The corresponding eigenspace is denoted by ${\mathfrak F}_{jl}^{(2)}(t)$; for $t\ne 0$ it coincides with the eigenspace
${\mathfrak F}(t;[(\gamma^\circ_{i_0+1} - c^\circ_{jl}/4 )t^2, (\gamma^\circ_{p} + c^\circ_{jl}/4 )t^2 ])$
of $A(t)$. The elements  $\varphi^{(i_0+1)}_{1}(t), \dots, \varphi^{(i_0+1)}_{k_{i_0+1}}(t);\dots; \varphi^{(p)}_{1}(t), \dots, \varphi^{(p)}_{k_p}(t)$ form an orthonormal basis in ${\mathfrak F}_{jl}^{(2)}(t)$.
Let $F^{(r)}_{jl}(t)$ be the orthogonal projections onto ${\mathfrak F}_{jl}^{(r)}(t)$,  $r=1,2$.
Then the spectral projection $F(t)$ of the operator $A(t)$ for the interval $[0,\delta]$
can be represented as
$$
F(t) = F_{jl}^{(1)}(t) + F_{jl}^{(2)}(t),\quad |t| \le t^{00}_{jl}.
\eqno(2.5)
$$

\smallskip\noindent\textbf{Proposition 2.1.}
\textit{For $|t|\le t^{00}_{jl}$ we have}
$$
\| F_{jl}^{(1)}(t) -( P_1 + \dots + P_{i_0}) \| \le C_{4,jl} |t|,
\eqno(2.6)
$$
$$
\| F_{jl}^{(2)}(t) - ( P_{i_0+1} + \dots + P_p) \| \le C_{4,jl} |t|.
\eqno(2.7)
$$
\textit{The number $t^{00}_{jl}$ is subject to} (2.3), (2.4), \textit{and the constant~$C_{4,jl}$ is given by}
$$
C_{4,jl} = \beta_4 \delta^{-1/2} \|X_1\|^5 (c^\circ_{jl})^{-2}.
\eqno(2.8)
$$

\smallskip\noindent\textbf{Proof.}
Consider the contour $\Gamma_1 \subset \C$ that envelops the interval
\hbox{$[\gamma_1^\circ - c^\circ_{jl}/4, \gamma_{i_0}^\circ + c^\circ_{jl}/4 ]$}
equidistantly at the distance $c^\circ_{jl}/4$.
For $|t|\le t^{00}_{jl}$ this contour encloses the first $i_0$ clusters of the eigenvalues of ${\mathfrak A}(t)$
and is separated from the other clusters. We have
$$
F_{jl}^{(1)}(t) = - \frac{1}{2\pi i} \intop_{\Gamma_1} ({\mathfrak A}(t) - z I)^{-1}\,dz,
\eqno(2.9)
$$
$$
P_1+\dots + P_{i_0} = - \frac{1}{2\pi i} \intop_{\Gamma_1} (SP - z I)^{-1}\,dz,
\eqno(2.10)
$$
where we  integrate  in the positive direction.
For $|t|\le t^{00}_{jl}$ and $z\in \Gamma_1$ we have
$$
\| (SP - z I)^{-1} \|  \le 2 (c^\circ_{jl})^{-1},
\eqno(2.11)
$$
$$
\| ({\mathfrak A}(t) - z I)^{-1} \|\le  4 (c^\circ_{jl})^{-1}.
\eqno(2.12)
$$

Next, by the resolvent identity,
$$
({\mathfrak A}(t) - z I)^{-1} - (SP - zI)^{-1} = ({\mathfrak A}(t) - z I)^{-1} \left( SP - {\mathfrak A}(t) \right) (SP - zI)^{-1}.
$$
Combining this with (2.2), (2.11), and (2.12), we obtain
$$
\| ({\mathfrak A}(t) - z I)^{-1} - (SP - zI)^{-1} \| \le 8 C_2 (c^\circ_{jl})^{-2} |t|,\quad z \in \Gamma_1,\quad |t|\le t^{00}_{jl}.
\eqno(2.13)
$$
It is easy to estimate the length of the contour $\Gamma_1$ by $2\|S\|$.
Now, relations (2.9), (2.10), (2.13), (1.4), and (1.15) imply (2.6) with
$C_{4,jl} = 8 \pi^{-1} C_2 \|X_1\|^2 (c^\circ_{jl})^{-2} = \beta_4 \delta^{-1/2} \|X_1\|^5 (c^\circ_{jl})^{-2}$.

Estimate (2.7) is proved similarly by integration over the contour $\Gamma_2$ that envelops the interval
$[\gamma_{i_0+1}^\circ - c^\circ_{jl}/4, \gamma_p^\circ + c^\circ_{jl}/4]$
equidistantly at the distance $c^\circ_{jl}/4$.
$\bullet$

\section*{§3. Threshold approximations for the operator exponential}

\noindent\textbf{3.1. Approximation of the operator $e^{- i \tau A(t)} P$.}
In this subsection, we approximate the operator $e^{- i \tau A(t)}  P$  by $e^{-i \tau t^2 SP}P$ for $\tau \in \R$ and $|t|\le t^0$.
Such approximation was found in [BSu5, \S 2].
Now we will repeat the proof of this result tracing more carefully how different terms are estimated.
This will be needed in what follows to confirm the sharpness of the result.
On the other hand,  we will distinguish an important case where the result of [BSu5] can be refined.

Consider the operator
$$
E(t,\tau):= e^{-i \tau A(t)} P - e^{-i \tau t^2 SP}P.
\eqno(3.1)
$$
We have
$$
E(t,\tau) = E_1(t,\tau) + E_2(t,\tau),
\eqno(3.2)
$$
$$
E_1(t,\tau) =  e^{-i \tau A(t)} F(t)^\perp P - F(t)^\perp e^{-i \tau t^2 SP}P,
\eqno(3.3)
$$
$$
E_2(t,\tau) = e^{-i \tau A(t)} F(t)P  - F(t) e^{-i \tau t^2 SP}P.
\eqno(3.4)
$$
Since  $F(t)^\perp P = (P - F(t))P$, (1.13)  implies the following estimate for the operator (3.3):
$$
\| E_1(t,\tau) \| \le 2 C_1 |t|,\quad |t| \le t^0.
\eqno(3.5)
$$

The operator (3.4) can be written as
$$
E_2(t,\tau) = e^{-i \tau A(t)} \Sigma(t,\tau),
\eqno(3.6)
$$
$$
\Sigma(t,\tau):= F(t)P - e^{i \tau A(t)} F(t) e^{-i \tau t^2 SP} P.
$$
Obviously,  $\Sigma(t,0)=0$ and
$$
\Sigma'(t,\tau):= \frac{d\Sigma}{d\tau}(t,\tau) = - i e^{i \tau A(t)} F(t) \left(A(t)F(t) - t^2 S P \right) e^{-i \tau t^2 SP} P.
\eqno(3.7)
$$
Since $\Sigma(t,\tau) = \int_0^\tau \Sigma'(t,\wt{\tau}) \,d\wt{\tau}$, by (3.7) and (1.14), we have
$$
\| \Sigma(t,\tau)\| \le C_2  |\tau| |t|^3,\quad |t| \le t^0.
\eqno(3.8)
$$

Relations (3.1), (3.2), (3.5), (3.6), and (3.8) imply the following result which is close to Theorem~2.1 from~ [BSu5].

\smallskip\noindent\textbf{Theorem 3.1.} \textit{Under the assumptions of Subsection} 1.1, \textit{for $\tau \in \R$ and $|t|\le t^0$ we have}
$$
\| e^{-i \tau A(t)} P - e^{-i \tau t^2 SP}P \| \le 2  C_1 |t| + C_2 |\tau| |t|^3.
\eqno(3.9)
$$
\textit{The number $t^0$ is subject to} (1.2), \textit{and the constants $C_1$ and $C_2$ are defined by}~(1.15).

\smallskip
Now we proceed to more subtle considerations that will allow us to improve the result under the additional assumptions.
Using representation (1.16), from (3.7) we obtain
$$
\Sigma(t,\tau) = - i \intop_0^{\tau} e^{i \widetilde{\tau} A(t)} F(t) \left(t^3 K + \Psi(t) \right) e^{-i \widetilde{\tau} t^2 SP} P
\,d \widetilde{\tau}.
\eqno(3.10)
$$
By (1.21), the operator (3.10) can be represented as
$$
\Sigma(t,\tau) = \widetilde{\Sigma}(t,\tau) + \widehat{\Sigma}(t,\tau),
\eqno(3.11)
$$
$$
\widetilde{\Sigma}(t,\tau) = - i \intop_0^{\tau} e^{i \widetilde{\tau} A(t)} F(t) \left(t^3 ZSP + \Psi(t) \right) e^{-i \widetilde{\tau} t^2 SP} P
\,d \widetilde{\tau},
\eqno(3.12)
$$
$$
\widehat{\Sigma}(t,\tau) = - i t^3 \intop_0^{\tau} e^{i \widetilde{\tau} A(t)} F(t)  N  e^{-i \widetilde{\tau} t^2 SP} P
\,d \widetilde{\tau}.
$$
Since $P Z=0$, then $F(t)ZSP  = (F(t) - P)ZSP$. Hence,  relations (1.4), (1.13), and (1.17) imply the following estimate for the term (3.12):
$$
\| \widetilde{\Sigma}(t,\tau)\| \le C_5 |\tau| t^4, \quad |t|\le t^0,
\eqno(3.13)
$$
where $C_5 = C_1 \|X_1\|^3 (8\delta)^{-1/2} + C_3 = \beta_5 \delta^{-1} \|X_1\|^4$.
(We have used expressions for $C_1$ and $C_3$.)

Now (3.6) and (3.11) imply that
$$
E_2(t,\tau) = \widetilde{E}_2(t,\tau)+ \widehat{E}_2(t,\tau),
\eqno(3.14)
$$
where $\widetilde{E}_2(t,\tau) = e^{-i \tau A(t)}  \widetilde{\Sigma}(t,\tau)$ and $\widehat{E}_2(t,\tau) = e^{-i \tau A(t)}  \widehat{\Sigma}(t,\tau)$.
By (3.13),
$$
\|\widetilde{E}_2(t,\tau)\| \le C_5 |\tau| t^4, \quad |t|\le t^0, \quad C_5 = \beta_5 \delta^{-1} \|X_1\|^4.
\eqno(3.15)
$$

Finally, relations (3.1), (3.2), (3.5), (3.14), (3.15), together with (1.20) imply the following result.

\smallskip\noindent\textbf{Theorem 3.2.} \textit{Under the assumptions of Subsection} 1.1, \textit{for $\tau \in \R$  and $|t|\le t^0$ we have
$$
e^{-i \tau A(t)} P - e^{-i \tau t^2 SP}P = E_1(t,\tau) + \widetilde{E}_2(t,\tau) + \widehat{E}_2(t,\tau),
$$
where the first two terms satisfy estimates} (3.5) \textit{and} (3.15), \textit{respectively}.
\textit{The third term admits the following representation}
$$
\widehat{E}_2(t,\tau) = - i t^3 e^{- i \tau A(t)}  \intop_0^{\tau} e^{i \widetilde{\tau} A(t)} F(t)  N  e^{-i \widetilde{\tau} t^2 SP} P \, d  \widetilde{\tau},
\eqno(3.16)
$$
\textit{where the operator $N$ is defined in Theorem} 1.4. \textit{We have}
$$
\| \widehat{E}_2(t,\tau) \| \le C_6 |\tau| |t|^3,\quad |t| \le t^0, \quad C_6 = (2 \delta)^{-1/2} \| X_1 \|^3.
\eqno(3.17)
$$

\smallskip\noindent\textbf{Corollary 3.3.} \textit{Suppose that the assumptions of Subsection} 1.1 \textit{are satisfied.
Suppose that $N=0$. Then for $\tau \in \R$ and $|t|\le t^0$ we have}
$$
\| e^{-i \tau A(t)} P  - e^{-i \tau t^2 SP}P \| \le  2 C_1 |t| + C_5 |\tau| t^4.
\eqno(3.18)
$$

\smallskip\noindent\textbf{3.2. Estimate of the term containing $N_*$.}
Assume that Condition 1.6 is satisfied. We will use the notation and the results of \S 2.
Recall that $N=N_0 + N_*$. By Remark 1.5, we have
$$
 P_j N_* P_j=0,\quad j=1,\dots,p; \quad P_l N_0 P_j = 0 \quad \text{for}\ l\ne j.
\eqno(3.19)
$$
Thus, relations (2.1) and (3.19) imply the following invariant representations for the operators $N_0$ and $N_*$:
$$
N_0 = \sum_{j=1}^p P_j N P_j, \quad N_* = \sum_{1\le l,j \le p: \, j\ne l} P_j N P_l.
\eqno(3.20)
$$

The term (3.16) can be written as
$$
\widehat{E}_2(t,\tau) = {E}_0(t,\tau) + {E}_*(t,\tau),
\eqno(3.21)
$$
where
$$
{E}_0(t,\tau) =
- i t^3 e^{- i \tau A(t)}  \intop_0^{\tau} e^{i \widetilde{\tau} A(t)} F(t)  N_0  e^{-i \widetilde{\tau} t^2 SP} P \, d  \widetilde{\tau},
\eqno(3.22)
$$
$$
{E}_*(t,\tau) =
- i t^3 e^{- i \tau A(t)}  \intop_0^{\tau} e^{i \widetilde{\tau} A(t)} F(t)  N_*  e^{-i \widetilde{\tau} t^2 SP} P \, d  \widetilde{\tau}.
\eqno(3.23)
$$

In this subsection, we obtain the analog of estimate (3.18) under the weaker assumption that $N_0=0$.
For this, we have to estimate the operator (3.23). However, we are able to do this only for a
smaller interval of $t$. By (3.20), the term (3.23) can be represented as
$$
{E}_*(t,\tau) = -i e^{- i \tau A(t)} \sum_{1\le j,l \le p:\ j \ne l} J_{jl}(t,\tau),
\eqno(3.24)
$$
$$
J_{jl}(t,\tau)=
 t^3    \intop_0^{\tau} e^{i \widetilde{\tau} A(t)} F(t) P_j  N P_l e^{-i \widetilde{\tau} t^2 SP} P \, d  \widetilde{\tau}.
\eqno(3.25)
$$

We have to estimate only those terms in (3.24) for which $P_j N P_l \ne 0$.
So, let $j\ne l$, and let $P_j N P_l \ne 0$. Suppose that $c^\circ_{jl}$ is defined by (2.3),
and $t^{00}_{jl}$ is subject to (2.4).
By (2.5), the operator (3.25) can be represented as
$$
J_{jl}(t,\tau)= J_{jl}^{(1)}(t,\tau) + J_{jl}^{(2)}(t,\tau),
\eqno(3.26)
$$
$$
J_{jl}^{(r)}(t,\tau) =  t^3   \intop_0^{\tau} e^{i \widetilde{\tau} A(t)} F_{jl}^{(r)}(t)  P_j  N  P_l e^{-i \widetilde{\tau} t^2 SP} P \, d  \widetilde{\tau},\quad r=1,2.
 \eqno(3.27)
$$

For definiteness, assume that $j<l$. Then $j < i_0+1$ and, by (2.1) and (2.7),
$$
\| F^{(2)}_{jl}(t) P_j\| = \| \left( F^{(2)}_{jl}(t) - (P_{i_0+1} + \dots + P_p) \right)  P_j\|  \le C_{4,jl} |t|.
\eqno(3.28)
$$
Combining (3.27), (3.28), and (1.20), (2.8), we obtain
$$
\| J^{(2)}_{jl}(t, \tau)\| \le C_{4,jl} t^4 |\tau| \|N\| \le C_{7,jl} |\tau| t^4,\quad |t| \le t^{00}_{jl},
\eqno(3.29)
$$
where
$
C_{7,jl} = C_{4,jl} (2\delta)^{-1/2} \|X_1\|^3 = \beta_7  \delta^{-1} \|X_1\|^8 (c^\circ_{jl})^{-2}.
$

It remains to consider the term $J_{jl}^{(1)}(t,\tau)$. Obviously, $P_l e^{-i \widetilde{\tau} t^2 SP} P= e^{-i \widetilde{\tau} \gamma_l^\circ t^2} P_l$.
Recall that (see Subsection 2.2) the  (nonzero) eigenvalues of the operator $A(t)F_{jl}^{(1)}(t)$ are
$t^2 \nu^{(1)}_1(t),\dots, t^2 \nu^{(1)}_{k_1}(t); \dots; t^2 \nu_1^{(i_0)}(t),\dots, t^2 \nu_{k_{i_0}}^{(i_0)}(t)$, and
$\nu^{(r)}_q(t) \in [ \gamma_1^\circ - c_{jl}^\circ /4,  \gamma_{i_0}^\circ + c_{jl}^\circ /4]$.
The corresponding orthonormal eigenvectors are
$\varphi_1^{(1)}(t),\dots, \varphi^{(1)}_{k_1}(t);\dots; \varphi_1^{(i_0)}(t),\dots, \varphi^{(i_0)}_{k_{i_0}}(t)$.
Then
$$
e^{i \widetilde{\tau} A(t)} F_{jl}^{(1)}(t) = \sum_{r=1}^{i_0} \sum_{q=1}^{k_r} e^{i \widetilde{\tau} t^2 \nu_q^{(r)}(t)} (\cdot, \varphi^{(r)}_q(t)) \varphi^{(r)}_q(t).
$$
As a result, the operator  $J_{jl}^{(1)}(t,\tau)$ can be written as
$$
J_{jl}^{(1)}(t,\tau) =  t^3  \sum_{r=1}^{i_0} \sum_{q=1}^{k_r} \left(\int_0^{\tau} e^{i \widetilde{\tau} t^2 (\nu_q^{(r)}(t)-\gamma_l^\circ)}  \, d  \widetilde{\tau}\right) ( P_j  N P_l \cdot, \varphi^{(r)}_q(t)) \varphi^{(r)}_q(t).
 \eqno(3.30)
$$
Calculating the integral in (3.30) and taking into account that $|\nu^{(r)}_q(t) - \gamma_l^\circ| \ge 3c^\circ_{jl}/4$
for $|t|\le t^{00}_{jl}$, we obtain
$$
\begin{aligned}
&\left|\int_0^{\tau} e^{i \widetilde{\tau} t^2 (\nu_q^{(r)}(t)-\gamma_l^\circ)}  \, d  \widetilde{\tau}\right|
\cr
&=
t^{-2} |\nu_q^{(r)}(t)-\gamma_l^\circ|^{-1} \left| e^{i {\tau} t^2 (\nu_q^{(r)}(t)-\gamma_l^\circ)} -1\right| \le 8 (3 c_{jl}^\circ)^{-1}t^{-2}.
\end{aligned}
\eqno(3.31)
$$
Now relations (3.30) and (3.31) together with (1.20) imply that
$$
\| J_{jl}^{(1)}(t,\tau)\| \le 8  (3c_{jl}^\circ)^{-1} |t| \|P_j N P_l\| \le C_{8,jl}  |t|,\quad |t|\le t^{00}_{jl},
\eqno(3.32)
$$
where
$
C_{8,jl} = 8 (2\delta)^{-1/2} \|X_1\|^3 (3c^\circ_{jl})^{-1}= \beta_8 \delta^{-1/2} \|X_1\|^3 (c_{jl}^\circ)^{-1}.
$
The case where $j>l$ can be treated similarly.

Let ${\mathcal Z} = \{ (j,l): 1\le j,l \le p,\ j\ne l,\ P_j N P_l \ne 0\}.$
Let $c^\circ_{jl}$ be defined by (2.3). We put
$$
c^\circ := \min_{(j,l)\in {\mathcal Z}} c^\circ_{jl},
\eqno(3.33)
$$
and choose a number $t^{00} \le t^0$ such that
$$
t^{00} \le (4 \beta_2)^{-1} \delta^{1/2} \|X_1\|^{-3} c^\circ.
\eqno(3.34)
$$
We may assume that  $t^{00} \le t^{00}_{jl}$ for all $(j,l)\in {\mathcal Z}$ (see (2.4)).

Now, relations (3.24), (3.26), (3.29), and (3.32), together with  expressions for the constants $C_{7,jl}$ and $C_{8,jl}$ yield
$$
\| {E}_*(t,\tau) \| \le  C_{7} |\tau| t^4 + C_{8} |t|, \quad |t|\le t^{00}.
\eqno(3.35)
$$
Here
$$
C_7 = \beta_7 n^2  \delta^{-1} \|X_1\|^8 (c^\circ)^{-2},
\quad
C_8 = \beta_8 n^2  \delta^{-1/2} \|X_1\|^3 (c^\circ)^{-1}.
\eqno(3.36)
$$

Finally, combining Theorem 3.2 and relations (3.21), (3.22), (3.35), and denoting
$\check{E}(t,\tau):= E_1(t,\tau) + \widetilde{E}_2(t,\tau) + {E}_*(t,\tau)$, we arrive at the following result.

\smallskip\noindent\textbf{Theorem 3.4.} \textit{Suppose that the assumptions of Subsection} 1.1 \textit{and Condition} 1.6
\textit{are satisfied. Suppose that the number $t^{00}\le t^0$ is subject to} (3.33), (3.34). \textit{Then for $\tau \in \R$ and $|t|\le t^{00}$ we have
$$
e^{-i \tau A(t)} P - e^{-i \tau t^2 SP}P = E_0(t,\tau) + \check{E}(t,\tau),
$$
where the second term satisfies}
$$
\| \check{E}(t,\tau)\| \le C_{9} |t| + C_{10} |\tau| t^4,\quad |t|\le t^{00}.
$$
\textit{The constants $C_{9}$ and $C_{10}$ are given by}
$C_{9}= 2 C_1 + C_8$, $C_{10}= C_5 + C_7$, \textit{where $C_1$, $C_5$ are defined by} (1.15), (3.15),
\textit{and $C_7$, $C_8$ are defined by} (3.36).
\textit{The operator $E_0(t,\tau)$ is given by} (3.22)
\textit{and satisfies the estimate
$$
\| {E}_0(t,\tau) \| \le C_6 |\tau| |t|^3,
$$
where $C_6$ is as in} (3.17).

\smallskip\noindent\textbf{Corollary 3.5.} \textit{Suppose that the assumptions of Subsection} 1.1 \textit{and Condition} 1.6
\textit{are satisfied. If $N_0=0$, then for $\tau \in \R$ and $|t|\le t^{00}$ we have}
$$
\| e^{-i \tau A(t)} P - e^{-i \tau t^2 SP}P  \| \le C_{9} |t| + C_{10} |\tau| t^4.
$$

\smallskip\noindent\textbf{Remark 3.6.}
Let $\mu_l$, $l=1,\dots,n,$  be the coefficients at $t^3$ in the expansions (1.6).
By Remark 1.5, the condition $N_0=0$ is equivalent to the relations $\mu_l=0$ for all $l=1,\dots,n$.

\section*{§4. Approximation of the operator  $e^{- i \eps^{-2} \tau A(t)}$}

\noindent\textbf{4.1. Approximation of the operator  $e^{- i \eps^{-2} \tau A(t)}$ in the general case.}
Let $\eps>0$. We study the behavior of the operator $e^{- i \eps^{-2} \tau A(t)}$ for small $\eps$.
We multiply this operator by the ``smoothing factor'' $\eps^s (t^2 + \eps^2)^{-s/2}P$, where $s>0$.
(The term is explained by the fact that in applications to differential operators  this factor turns into the smoothing operator.)
Our goal is to find approximation for the smoothed operator exponential with an error $O(\eps)$ for minimal possible $s$.

Let $|t|\le t^0$. We apply Theorem 3.1. By (3.9) (with $\tau$ replaced by $\eps^{-2}\tau$),
$$
\begin{aligned}
&\| e^{- i \eps^{-2} \tau A(t)} P  - e^{-i \eps^{-2} \tau t^2 SP} P  \| \eps^3 (t^2 + \eps^2)^{-3/2}
\cr
&\le  (2 C_1 |t| + C_2 \eps^{-2} |\tau| |t|^3 ) \eps^3 (t^2 + \eps^2)^{-3/2} \le
({C_1} + C_2 |\tau|) \eps.
\end{aligned}
$$
Here we take $s=3$. We arrive at the following result which has been proved before in [BSu5, Theorem 2.6].

\smallskip\noindent\textbf{Theorem 4.1.} \textit{Suppose that the assumptions of Subsection} 1.1
\textit{are satisfied. Then for $\eps>0$, $\tau \in \R$, and $|t|\le t^0$ we have}
$$
\| e^{- i \eps^{-2} \tau A(t)} P - e^{-i \eps^{-2} \tau t^2 SP} P  \| \eps^3 (t^2 + \eps^2)^{-3/2}
\le  (C_1  + C_2  |\tau| ) \eps.
\eqno(4.1)
$$
\textit{The number $t^0$ is subject to} (1.2), \textit{and the constants $C_1$, $C_2$ are given by}  (1.15).

\noindent\textbf{4.2. Refinement of approximation for  $e^{- i \eps^{-2} \tau A(t)}$ under the additional assumptions.}
Corollary~3.3 allows us to improve the result of Theorem~4.1 in the case where~$N=0$.

\smallskip\noindent\textbf{Theorem 4.2.} \textit{Suppose that the assumptions of Theorem}~4.1
\textit{are satisfied. Suppose that the operator $N$ defined in Theorem}~1.4 \textit{is equal to zero}: $N=0$.
\textit{Then for $\eps>0$, $\tau \in \R$, and $|t|\le t^0$ we have}
$$
\| e^{- i \eps^{-2} \tau A(t)} P  - e^{-i \eps^{-2} \tau t^2 SP} P  \| \eps^2 (t^2 + \eps^2)^{-1}
\le  ({C}'_1  + C_5  |\tau| ) \eps.
\eqno(4.2)
$$
\textit{Here ${C}'_1 = \max\{2, C_1\}$ and $C_5 = \beta_5 \delta^{-1} \|X_1\|^4$.}

\smallskip\noindent\textbf{Proof.}
Note that for $|t| \ge \sqrt{\eps}$ we have $\eps^2 (t^2 + \eps^2)^{-1} \le \eps$, whence the left-hand side of (4.2) does not exceed $2\eps$.

Thus, we may assume that $|t|< \sqrt{\eps}$.
Using (3.18) with $\tau$ replaced by $\eps^{-2}\tau$, for $|t|< \sqrt{\eps}$ we obtain
$$
\begin{aligned}
&\|  e^{- i \eps^{-2} \tau A(t)} P  - e^{-i \eps^{-2} \tau t^2 SP} P  \| \eps^2 (t^2 + \eps^2)^{-1}
\cr
&\le  (2 C_1 |t| + C_5 \eps^{-2} |\tau| t^4 ) \eps^2 (t^2 + \eps^2)^{-1} \le
{C_1} \eps + C_5 |\tau| t^2 \le ({C_1}+ C_5 |\tau|) \eps.
\end{aligned}
$$
The required statement follows. $\bullet$

\smallskip
Similarly, Corollary 3.5 yields the following result.

\smallskip\noindent\textbf{Theorem 4.3.} \textit{Suppose that the assumptions of Subsection}~1.1 \textit{and Condition}~1.6
\textit{are satisfied. Suppose that the operator $N_0$ defined in Theorem}~1.4 \textit{is equal to zero}: $N_0 =0$.
\textit{Then for $\eps>0$, $\tau \in \R$, and $|t|\le t^{00}$ we have}
$$
\|  e^{- i \eps^{-2} \tau A(t)} P  - e^{-i \eps^{-2} \tau t^2 SP} P  \| \eps^2 (t^2 + \eps^2)^{-1}
\le  ({C}'_{9}  + C_{10}  |\tau| ) \eps.
$$
\textit{Here the number $t^{00}\le t^0$ is subject to} (3.34), \textit{the constant ${C}'_{9}$ is given by}
$C_9' = \max\{ 2, \frac{1}{2} C_{9} \}$, \textit{and $C_9$, $C_{10}$ are defined in Theorem} 3.4.

\smallskip\noindent\textbf{4.3. Sharpness of the result in the general case.}
Now we show that the result of Theorem~4.1 is sharp in the general case.
Namely, if $N_0 \ne 0$, the exponent $s$ in the smoothing factor can not be taken smaller than 3.

\smallskip\noindent\textbf{Theorem 4.4.} \textit{Let $N_0 \ne 0$. Let $0\ne \tau \in \R$. Then for any $1\le s<3$
it is impossible that the estimate}
$$
\|  e^{- i \eps^{-2} \tau A(t)} P  - e^{-i \eps^{-2} \tau t^2 SP} P  \| \eps^s (t^2 + \eps^2)^{-s/2}
\le  C(\tau) \eps
\eqno(4.3)
$$
\textit{holds for all sufficiently small $|t|$ and $\eps >0$.}

\smallskip\noindent\textbf{Proof.} We start with a preliminary remark.
Since $F(t)^\perp P = (P - F(t)) P$, from (1.13) it follows that
$$
\| e^{- i \eps^{-2} \tau A(t)} F(t)^\perp P \| \eps (t^2 + \eps^2)^{-1/2} \le C_1 |t|\eps (t^2 + \eps^2)^{-1/2} \le C_1 \eps, \quad |t| \le t^0.
\eqno(4.4)
$$

Let us fix $0 \ne \tau \in \R$. We prove by contradiction. Suppose that for some $1\le s<3$
there exists a constant $C(\tau)>0$ such that (4.3) is valid for all sufficiently small $|t|$ and $\eps$.
By (4.4), this assumption is equivalent to
the existence of a constant $\widetilde{C}(\tau)$ such that
$$
\| \left( e^{- i \eps^{-2} \tau A(t)} F(t)  - e^{-i \eps^{-2} \tau t^2 SP} P \right)P \| \eps^s (t^2 + \eps^2)^{-s/2}
\le  \widetilde{C}(\tau) \eps
\eqno(4.5)
$$
for all sufficiently small $|t|$ and $\eps$.

Consider the interval  $|t|\le t_*$ of convergence of the power series expansions (1.6), (1.7). (Now we use the initial enumeration.)
We have
$$
 e^{- i \eps^{-2} \tau A(t)} F(t)   = \sum_{l=1}^n e^{-i\eps^{-2} \tau \lambda_l(t)} (\cdot, \varphi_l(t)) \varphi_l(t).
 \eqno(4.6)
$$
From the convergence of the power series expansions (1.7) it follows that
$$
\| \varphi_l(t) - \omega_l \| \le c_1 |t|,\quad |t|\le t_*,\quad l=1,\dots,n.
\eqno(4.7)
$$
Relations (4.6) and (4.7) show that
$$
\bigl\| e^{- i \eps^{-2} \tau A(t)} F(t) - \sum_{l=1}^n e^{-i\eps^{-2} \tau \lambda_l(t)} (\cdot, \omega_l) \omega_l \bigr\|
\le c_2 |t|,\quad |t|\le t_*.
\eqno(4.8)
$$

Comparing (4.8) and (4.5), we see that there exists a constant $\widehat{C}(\tau)$ such that
$$
\bigl\|  \sum_{l=1}^n \left( e^{- i \eps^{-2} \tau \lambda_l(t)}   - e^{-i \eps^{-2} \tau \gamma_l t^2} \right)  (\cdot, \omega_l) \omega_l \bigr\| \eps^s (t^2 + \eps^2)^{-s/2}
\le  \widehat{C}(\tau) \eps
\eqno(4.9)
$$
for all sufficiently small $|t|$ and $\eps$.

The condition $N_0\ne 0$ means that $\mu_j \ne 0$ at least for one $j$.
Applying the operator under the norm sign in (4.9) to $\omega_j$, we obtain
$$
\left| e^{- i \eps^{-2} \tau \lambda_j(t)}   - e^{-i \eps^{-2} \tau \gamma_j t^2} \right|  \eps^s (t^2 + \eps^2)^{-s/2} \le \widehat{C}(\tau) \eps
\eqno(4.10)
$$
for all sufficiently small $|t|$ and $\eps$.
The left-hand side of  (4.10) can be written as $2 |\sin \frac{1}{2}\eps^{-2} \tau (\lambda_j(t) - \gamma_j t^2)| \eps^s (t^2 + \eps^2)^{-s/2}$.
 Using that the expansion (1.6) for $\lambda_j(t)$  is convergent and $\mu_j \ne 0$,
 we may assume that
 $$
 \frac{1}{2} |\mu_j| |t|^3 \le |\lambda_j (t) - \gamma_j t^2| \le \frac{3}{2} |\mu_j| |t|^3,\quad |t| \le t_*,
  $$
   possibly diminishing $t_*$. Hence,
   $$
   \frac{1}{4} \eps^{-2}|\mu_j \tau| |t|^3 \le \left| \frac{1}{2}\eps^{-2} \tau (\lambda_j (t) - \gamma_j t^2) \right| \le \frac{3}{4} \eps^{-2}|\mu_j \tau| |t|^3,\quad |t| \le t_*.
   $$
Now, for a fixed $\tau \ne 0$, assuming that $\eps$ is sufficiently small
(namely, such that $\pi^{1/3} |\mu_j \tau|^{-1/3} \eps^{2/3} \le t_*$), we put
$t = t(\eps) = \pi^{1/3} |\mu_j \tau|^{-1/3} \eps^{2/3}= c \,\eps^{2/3}$. For such $t$ we have
$2 |\sin \frac{1}{2}\eps^{-2} \tau (\lambda_j(t) - \gamma_j t^2)| \ge \sqrt{2}$,
whence (4.10) implies that
$\sqrt{2} \eps^s (c^2 \eps^{4/3} + \eps^2)^{-s/2} \le \widehat{C} \eps$.
This means that the function \hbox{$\eps^{s/3 -1} (c^2 + \eps^{2/3})^{-s/2}$}
is uniformly bounded for small $\eps$. But this is not true provided that $s<3$.
This contradiction completes the proof. $\bullet$

\section*{\S 5. Approximation of the sandwiched operator exponential}

\smallskip\noindent\textbf{5.1. The operator family  $A(t)= M^* \widehat{A}(t)M$.}
Let $\widehat{\H}$ be yet another separable Hilbert space. Let $\widehat{X}(t) = \wh{X}_0 + t \wh{X}_1: \wh{\H} \to \H_*$
be the family of operators of the same form as  $X(t)$, and suppose that $\wh{X}(t)$ satisfies the assumptions of Subsection~1.1.
Let $M:\H \to \widehat{\H}$ be an isomorphism. Suppose that
$M {\rm Dom}\, X_0 = {\rm Dom}\, \widehat{X}_0$, $X(t) = \widehat{X}(t) M$, and then also $X_0 = \widehat{X}_0  M$,
$X_1 = \widehat{X}_1  M$. In $\widehat{\H}$, we consider the family of operators $\wh{A}(t)= \wh{X}(t)^* \wh{X}(t)$.
Then
$$
A(t)= M^* \widehat{A}(t)M.
\eqno(5.1)
$$
In what follows, all the objects corresponding to the family $\wh{A}(t)$ are marked by ``hat''. Note that  $\wh{\NN} = M \NN$ and
$\wh{\NN}_* =\NN_*$.

In $\wh{\H}$ we consider the positive definite operator $Q := (M M^*)^{-1}$.
Let $Q_{\wh{\NN}} = \wh{P} Q \vert_{\wh{\NN}}$ be the block of  $Q$ in the subspace~$\wh{\NN}$.
Obviously, $Q_{\wh{\NN}}$ is an isomorphism in $\wh{\NN}$.

According to [Su2, Proposition 1.2], the orthogonal projection $P$ of $\H$ onto $\NN$ and the orthogonal projection
 $\wh{P}$ of $\wh{\H}$ onto $\wh{\NN}$ satisfy the following relation
 $$
 P = M^{-1} (Q_{\wh{\NN}})^{-1} \wh{P} (M^*)^{-1}.
 \eqno(5.2)
 $$
 Let $\wh{S}: \wh{\NN} \to \wh{\NN}$ be the spectral germ of $\wh{A}(t)$ at $t=0$, and let $S$ be the germ of  $A(t)$.
According to [BSu1, Chapter 1, Subsection 1.5], we have
$$
S = P M^* \wh{S} M \vert_{\NN}.
\eqno(5.3)
$$

\smallskip\noindent\textbf{5.2. The operators  $\widehat{Z}_Q$ and $\wh{N}_Q$.}
For the operator family $\wh{A}(t)$ we introduce the operator $\wh{Z}_Q$ acting in $\wh{\H}$
and taking an element $\wh{u} \in \wh{\H}$ to the solution $\wh{\psi}_Q$ of the problem
$$
\wh{X}_0^* (\wh{X}_0 \wh{\psi}_Q + \wh{X}_1 \wh{\omega})=0,\quad Q \wh{\psi}_Q \perp \wh{\NN},
$$
where $\wh{\omega} = \wh{P}\wh{u}$.
As shown in  [BSu2, \S 6], the operator $Z$ for $A(t)$ and the operator $\wh{Z}_Q$ introduced above
satisfy $\wh{Z}_Q = M Z M^{-1} \wh{P}$. Next, we put
$$
\wh{N}_Q := \wh{Z}_Q^* \wh{X}^*_1 \wh{R} \wh{P} + (\wh{R} \wh{P})^* \wh{X}_1 \wh{Z}_Q.
\eqno(5.4)
$$
According to [BSu2, \S 6], the operator $N$ for $A(t)$ and the operator (5.4)  satisfy
$$
\wh{N}_Q = \wh{P} (M^*)^{-1} N M^{-1} \wh{P}.
\eqno(5.5)
$$
Recall that $N=N_0 + N_*$ and introduce the operators
$$
\wh{N}_{0,Q} = \wh{P} (M^*)^{-1} N_0 M^{-1} \wh{P}, \quad
\wh{N}_{*,Q} = \wh{P} (M^*)^{-1} N_* M^{-1} \wh{P}.
\eqno(5.6)
$$
Then $\wh{N}_Q = \wh{N}_{0,Q} + \wh{N}_{*,Q}$.

\smallskip\noindent\textbf{Lemma 5.1.} \textit{Suppose that the assumptions of Subsection} 5.1
\textit{are satisfied. Suppose that the operators $N$ and $N_0$ are defined in Theorem}~1.4, \textit{and the operators $\wh{N}_Q$ and $\wh{N}_{0,Q}$ are defined in Subsection}~5.2.
\textit{Then the relation \hbox{$N=0$} is equivalent to the relation $\wh{N}_Q=0$.
The relation $N_0=0$ is equivalent to the relation $\wh{N}_{0,Q}=0$.}

\smallskip\noindent\textbf{Proof.} By (5.5), the relation $N=0$ implies that  $\wh{N}_Q=0$.

Conversely, let $\wh{N}_Q=0$. Then, by (5.5),
$\wh{P} (M^*)^{-1} N M^{-1} \wh{\omega}=0$ for any $\wh{\omega}\in \wh{\NN}$.
Since $M^{-1}$ is an isomorphism of $\wh{\NN}$ onto $\NN$, then
$\wh{P} (M^*)^{-1} N {\omega}=0$ for any ${\omega}\in {\NN}$. Multiplying the last relation by $\wh{\eta} \in \wh{\NN}$,
we obtain \hbox{$(N \omega, M^{-1}\wh{\eta})=0$} for any $\wh{\eta} \in \wh{\NN}$.
Using again that $M^{-1}\wh{\NN} = \NN$, we conclude that the block of $N$ in the subspace $\NN$ is equal to zero.
It remains to recall that $N$ takes $\NN$ to $\NN$, and $\NN^\perp$ to $\{0\}$. Hence, $N=0$.

The second statement can be checked in a similar way. $\bullet$

\smallskip\noindent\textbf{5.3. Relations between the operators and the coefficients of the power series expansions.}
Now we describe relations between the coefficients of the power series expansions (1.6), (1.7)
and the operators $\wh{S}$ and $Q_{\wh{\NN}}$. (See [BSu3, Subsections 1.6, 1.7].)  We denote $\zeta_l := M \omega_l \in \wh{\NN}$, $l=1, \dots,n$.
Then relations (1.9), (5.2), and (5.3) show that
$$
\wh{S} \zeta_l = \gamma_l Q_{\wh{\NN}}\zeta_l,\quad l=1,\dots,n.
\eqno(5.7)
$$
The set $\zeta_1,\dots, \zeta_n$ forms a basis in $\wh{\NN}$ that is orthonormal with the weight $Q_{\wh{\NN}}$:
   $$
   ( Q_{\wh{\NN}}\zeta_l, \zeta_j) = \delta_{lj},\quad l,j=1, \dots,n.
   \eqno(5.8)
   $$

The operators $\wh{N}_{0,Q}$ and $\wh{N}_{*,Q}$ can be described in terms of the coefficients
of the expansions (1.6) and (1.7); cf. (1.18), (1.19). We put $\wt{\zeta}_l := M \wt{\omega}_l \in \wh{\NN}$, $l=1,\dots,n$. Then
$$
\wh{N}_{0,Q} = \sum_{k=1}^n \mu_k (\cdot, Q_{\wh{\NN}} \zeta_k) Q_{\wh{\NN}} \zeta_k,
\eqno(5.9)
$$
$$
\wh{N}_{*,Q} = \sum_{k=1}^n \gamma_k \left( (\cdot, Q_{\wh{\NN}} \wt{\zeta}_k) Q_{\wh{\NN}} \zeta_k
+   (\cdot, Q_{\wh{\NN}} \zeta_k) Q_{\wh{\NN}} \wt{\zeta}_k \right).
\eqno(5.10)
$$

\smallskip\noindent\textbf{Remark 5.2.} By (5.8) and (5.9), we have
$$
(\wh{N}_{0,Q} \zeta_j, \zeta_l ) = \mu_l \delta_{jl},\quad j,l=1,\dots, n.
$$
From (5.8) and (5.10) it follows that
$$
(\wh{N}_{*,Q} \zeta_j, \zeta_l ) = \gamma_l (\zeta_j, Q_{\wh{\NN}} \wt{\zeta}_l)  + \gamma_j ( Q_{\wh{\NN}} \wt{\zeta}_j, \zeta_l),\quad j,l=1,\dots,n.
$$
Relations (1.12) imply that
$$
( Q_{\wh{\NN}} \wt{\zeta}_j , \zeta_l)  + ( \zeta_j, Q_{\wh{\NN}} \wt{\zeta}_l)=0 ,\quad j,l=1,\dots,n.
$$
Hence,
$$
(\wh{N}_{*,Q} \zeta_j, \zeta_l ) =0 \quad \text{if}\ \gamma_j = \gamma_l.
$$

Now we return to the notation of \S 2. Recall that the different eigenvalues of the germ $S$ are denoted by
$\gamma_j^\circ$, $j=1,\dots,p$, and the corresponding eigenspaces by $\NN_j$.
The vectors $\omega^{(j)}_i$, $i=1,\dots,k_j$, form an orthonormal basis in  $\NN_j$.
Then the same numbers $\gamma_j^\circ$, $j=1,\dots,p$, are different eigenvalues of the problem (5.7), and $M\NN_j$
are the corresponding eigenspaces. The vectors $\zeta^{(j)}_i = M \omega_i^{(j)}$, $i=1,\dots,k_j$, form a basis in $M\NN_j$
(orthonormal with the weight $Q_{\wh{\NN}}$). By ${\mathcal P}_j$ we denote the ``skew'' projection onto $M\NN_j$
that  is orthogonal with respect to the inner product $(Q_{\wh{\NN}} \cdot,\cdot)$, i.~e.,
$$
{\mathcal P}_j = \sum_{i=1}^{k_j} (\cdot, Q_{\wh{\NN}} \zeta^{(j)}_i) \zeta^{(j)}_i,\quad j=1,\dots,p.
$$
It is easily seen that ${\mathcal P}_j = M P_j M^{-1} \wh{P}$.

Using (3.20), (5.5), and (5.6), it is easy to check that
$$
\wh{N}_{0,Q} = \sum_{j=1}^p {\mathcal P}^*_j \wh{N}_Q {\mathcal P}_j,
\quad
\wh{N}_{*,Q} = \sum_{1 \le l,j \le p:\,l\ne j} {\mathcal P}^*_l \wh{N}_Q {\mathcal P}_j.
\eqno(5.11)
$$
Relations (5.11) are similar to (3.20); they  give the invariant representations for the operators $\wh{N}_{0,Q}$ and $\wh{N}_{*,Q}$.

\smallskip\noindent\textbf{5.4. Approximation of the sandwiched exponential.}
In this subsection, we find an approximation for the operator exponential $e^{-i \tau A(t)}$ of the family (5.1) in terms of the germ $\wh{S}$ of $\wh{A}(t)$
and the isomorphism $M$. It is convenient to border the exponential by appropriate factors.

We put $M_0 = (Q_{\wh{\NN}})^{-1/2}$. According to [BSu5, Proposition 3.1], we have
$$
M e^{-i \tau t^2 SP} P M^* = M_0 e^{-i \tau t^2 M_0 \wh{S} M_0} M_0 \wh{P}.
\eqno(5.12)
$$

\smallskip\noindent\textbf{Lemma 5.3.} \textit{Under the assumptions of Subsection}~5.1, \textit{we have}
$$
\begin{aligned}
&\|M e^{-i \tau A(t)} M^{-1} \wh{P} - M_0 e^{-i \tau t^2 M_0 \wh{S} M_0} M_0^{-1} \wh{P}\|
\cr
&\le \|M\|^2 \|M^{-1}\|^2 \| e^{-i \tau A(t)} P - e^{-i \tau t^2 S P} P \|,
\end{aligned}
\eqno(5.13)
$$
$$
\begin{aligned}
&\| e^{-i \tau A(t)} P - e^{-i \tau t^2 S P} P \|
\cr
&\le
\|M\|^2 \|M^{-1}\|^2 \|M e^{-i \tau A(t)} M^{-1} \wh{P} - M_0 e^{-i \tau t^2 M_0 \wh{S} M_0} M_0^{-1} \wh{P}\|.
\end{aligned}
\eqno(5.14)
$$

\smallskip\noindent\textbf{Proof.}
Denote the left-hand side of (5.13) by $J(t,\tau)$.
Since $M_0 = Q_{\wh{\NN}}^{-1/2}$, then
$$
J(t,\tau) = \|\left(M e^{-i \tau A(t)} M^{-1} Q_{\wh{\NN}}^{-1}\wh{P} - M_0 e^{-i \tau t^2 M_0 \wh{S} M_0} M_0 \wh{P}\right) Q_{\wh{\NN}} \wh{P}\|.
$$
Next, using the identity $M^{-1} Q_{\wh{\NN}}^{-1}\wh{P} = P M^*$ (see (5.2)) and (5.12), we obtain
$$
J(t,\tau) = \|  \left(M e^{-i \tau A(t)} P M^*  -  M e^{-i \tau t^2  {S} P} P M^* \right)  Q_{\wh{\NN}} \wh{P}\|.
$$
Hence,
$$
J(t,\tau) \le  \|M\|^2   \|Q_{\wh{\NN}} \wh{P}\| \|   e^{-i \tau A(t)} P  -  e^{-i \tau t^2  {S} P} P  \|.
$$
Since $\|Q_{\wh{\NN}}\| \le \|Q\| = \|M^{-1}\|^2$, we arrive at (5.13).

Estimate (5.14) can be checked similarly in the ``inverse way''.
Obviously,
$$
\| e^{-i \tau A(t)} P - e^{-i \tau t^2 S P} P \| \le \|M^{-1}\|^2 \|M e^{-i \tau A(t)} P M^*   - M e^{-i \tau t^2  {S} P} P M^* \|.
$$
By the identity $P M^* = M^{-1} Q_{\wh{\NN}}^{-1}\wh{P}$ and (5.12), the right-hand side can be written as
$\|M^{-1}\|^2 \|M e^{-i \tau A(t)} M^{-1}  Q_{\wh{\NN}}^{-1} \wh{P}  -
M_0 e^{-i \tau t^2  M_0 \wh{S} M_0} M_0^{-1}  Q_{\wh{\NN}}^{-1} \wh{P}\|.$
Together with the inequality $\|Q_{\wh{\NN}}^{-1} \wh{P}\| \le \| Q^{-1} \| = \|M\|^2$, this implies (5.14).
$\bullet$

  \smallskip
  Now, Theorem 3.1 and inequality (5.13) directly imply the following result (which has been obtained before in [BSu5, Subsection~3.2]).

\smallskip\noindent\textbf{Theorem 5.4.} \textit{Under the assumptions of Subsection}~5.1, \textit{for $\tau \in \R$ and $|t|\le t^0$ we have}
$$
\begin{aligned}
&\| M e^{-i \tau A(t)} M^{-1} \wh{P} - M_0 e^{-i \tau t^2 M_0 \wh{S} M_0} M_0^{-1} \wh{P} \|
\cr
&\le \| M \|^2 \| M^{-1}\|^2 (2 C_1 |t| + C_2 |\tau| |t|^3).
\end{aligned}
\eqno(5.15)
$$
\textit{The number $t^0$ is subject to} (1.2), \textit{and the constants $C_1$, $C_2$ are given by} (1.15).

\smallskip
Similarly,  combining Corollary 3.3, Lemma 5.1, and Lemma 5.3, we arrive at the following result.

\smallskip\noindent\textbf{Theorem 5.5.} \textit{Suppose that the assumptions of Subsection} 5.1 \textit{are satisfied.
Suppose that the operator $\wh{N}_{Q}$ defined in Subsection}~5.2 \textit{is equal to zero}: $\wh{N}_{Q}=0$.
\textit{Then for $\tau \in \R$ and $|t|\le t^0$ we have}
$$
\begin{aligned}
&\| M e^{-i \tau A(t)} M^{-1} \wh{P} - M_0 e^{-i \tau t^2 M_0 \wh{S} M_0} M_0^{-1} \wh{P} \|
\cr
&\le \| M \|^2 \| M^{-1}\|^2 (2 C_1 |t| + C_5 |\tau| t^4).
\end{aligned}
$$
\textit{The number $t^0$ is subject to} (1.2), \textit{and the constants $C_1$, $C_5$ are defined by} (1.15) \textit{and} (3.15).

\smallskip
Finally,  from Corollary 3.5, Lemma 5.1, and Lemma 5.3 we deduce the following statement.

\smallskip\noindent\textbf{Theorem 5.6.} \textit{Suppose that the assumptions of Subsection}~5.1 \textit{and Condition}~1.6 \textit{are satisfied.
Suppose that the operator $\wh{N}_{0,Q}$ defined in Subsection} 5.2 \textit{is equal to zero}: $\wh{N}_{0,Q}=0$.
\textit{Then for $\tau \in \R$ and $|t|\le t^{00}$ we have }
$$
\| M e^{-i \tau A(t)} M^{-1} \wh{P} - M_0 e^{-i \tau t^2 M_0 \wh{S} M_0} M_0^{-1} \wh{P} \|
\le \| M \|^2 \| M^{-1}\|^2 (C_{9} |t| + C_{10} |\tau| t^4).
$$
\textit{The number $t^{00}$ is subject to} (3.34), \textit{and the constants $C_{9}$, $C_{10}$ are as in Theorem} 3.4.

\subsection*{5.5. Approximation of the smoothed sandwiched exponential}
Writing down (5.15) with $\tau$ replaced by $\eps^{-2} \tau$ and multiplying it by the ``smoothing factor'', we arrive at the following result,
which has been proved before in  [BSu5, Theorem 3.2].

\smallskip\noindent\textbf{Theorem 5.7.} \textit{Under the assumptions of Subsection} 5.1, \textit{for $\tau \in \R$, $\eps >0$, and $|t|\le t^0$ we have}
$$
\begin{aligned}
\| M e^{-i \eps^{-2}\tau A(t)} M^{-1} \wh{P} - M_0 e^{-i \eps^{-2} \tau t^2 M_0 \wh{S} M_0} M_0^{-1} \wh{P} \|
\eps^3 (t^2 + \eps^2)^{-3/2}
\cr
\le \| M \|^2 \| M^{-1}\|^2 (C_1  + C_2 |\tau| ) \eps.
\end{aligned}
\eqno(5.16)
$$
\textit{The number $t^0$ is subject to} (1.2), \textit{and the constants $C_1$, $C_2$ are defined by} (1.15).

\smallskip
Similarly to the proof of Theorem 4.2, from Theorem 5.5 we deduce the following statement.

\smallskip\noindent\textbf{Theorem 5.8.} \textit{Suppose that the assumptions of Subsection}~5.1 \textit{are satisfied. Suppose that the operator $\wh{N}_{Q}$ defined in Subsection}~5.2
\textit{is equal to zero}: $\wh{N}_{Q}=0$. \textit{Then for $\tau \in \R$, $\eps >0$, and $|t|\le t^0$ we have}
$$
\begin{aligned}
\| M e^{-i \eps^{-2} \tau A(t)} M^{-1} \wh{P} - M_0 e^{-i \eps^{-2} \tau t^2 M_0 \wh{S} M_0} M_0^{-1} \wh{P} \|
\eps^2 (t^2 + \eps^2)^{-1}
\cr
\le \| M \|^2 \| M^{-1}\|^2 ({C}'_1  + C_5 |\tau|) \eps.
\end{aligned}
$$
\textit{The number $t^0$ is subject to} (1.2), \textit{and the constants ${C}'_1$ and $C_5$ are as in Theorem} 4.2.

\smallskip
Finally, Theorem 5.6 implies the following result.

\smallskip\noindent\textbf{Theorem 5.9.} \textit{Suppose that the assumptions of Subsection}~5.1 \textit{and Condition}~1.6 \textit{are satisfied. Suppose that the operator
$\wh{N}_{0,Q}$ defined in Subsection}~5.2 \textit{is equal to zero}: $\wh{N}_{0,Q}=0$.
\textit{Then for $\tau \in \R$, $\eps >0$, and $|t|\le t^{00}$ we have}
$$
\begin{aligned}
\| M e^{-i \eps^{-2}\tau A(t)} M^{-1} \wh{P} - M_0 e^{-i \eps^{-2}\tau t^2 M_0 \wh{S} M_0} M_0^{-1} \wh{P} \|  \eps^2 (t^2 + \eps^2)^{-1}
\cr
\le \| M \|^2 \| M^{-1}\|^2 ({C}'_{9} + C_{10} |\tau| ) \eps.
\end{aligned}
$$
\textit{The number $t^{00}$ is subject to} (3.33), (3.34), \textit{and the constants ${C}'_{9}$, $C_{10}$ are as in Theorem}~4.3.

\subsection*{5.6. The sharpness of the result}
Now we confirm that the result of Theorem~5.7 is sharp in the general case.
Namely, if $\wh{N}_{0,Q} \ne 0$, then the exponent $s$ in the smoothing factor can not be taken smaller than 3.

\smallskip\noindent\textbf{Theorem 5.10.} \textit{Suppose that the assumptions of Subsection} 5.1 \textit{are satisfied. Let $\wh{N}_{0,Q} \ne 0$. Let $0\ne \tau \in \R$.
Then for any $1\le s<3$ it is impossible that the estimate}
$$
\| M e^{-i \eps^{-2}\tau A(t)} M^{-1} \wh{P} - M_0 e^{-i \eps^{-2}\tau t^2 M_0 \wh{S} M_0} M_0^{-1} \wh{P} \|  \eps^s (t^2 + \eps^2)^{-s/2}
\le  C(\tau) \eps
\eqno(5.17)
$$
\textit{holds for all sufficiently small $|t|$ and $\eps >0$.}

\smallskip\noindent\textbf{Proof.} By Lemma~5.1, under our assumptions we have $N \ne 0$.
We prove by contradiction. Let us fix $\tau \ne 0$. Suppose that for some $1\le s<3$ there exists a constant $C(\tau)>0$ such that (5.17) holds for all
sufficiently small $|t|$ and $\eps$. By (5.14), this means that the inequality of the form (4.3)  also holds (with some other constant).
But this contradicts the statement of Theorem~4.4. $\bullet$

\section*{Chapter 2. Periodic differential operators in $L_2(\R^d;\C^n)$}

\section*{§6. The class of operators.\\ Lattices and the Gelfand transformation}

\subsection*{6.1. Factorized second order operators}
Let $b(\D) = \sum_{l=1}^d b_l D_l$ be a matrix first order differential operator; here $b_l$ are constant  $(m \times n)$-matrices
(in general, with complex entries).
\textit{Assume that} $m \ge n$. Consider the symbol $b(\bxi)= \sum_{l=1}^d b_l \xi_l$, $\bxi \in \R^d$, and \textit{assume that}
  $$
  {\rm rank}\, b(\bxi) =n,\quad 0 \ne \bxi \in \R^d.
  \eqno(6.1)
  $$
Condition (6.1) is equivalent to the inequalities
$$
\alpha_0 \1_n \le b(\bt)^* b(\bt) \le \alpha_1 \1_n,\quad |\bt|=1,\quad 0< \alpha_0 \le \alpha_1 < \infty,
\eqno(6.2)
$$
with some positive constants $\alpha_0$ and $\alpha_1$.

Suppose that an $(m\times m)$-matrix-valued function $h(\x)$ and an $(n\times n)$-matrix-valued function $f(\x)$ (in general, with complex entries) are such that
$$
h, h^{-1} \in L_\infty(\R^d);\quad f,f^{-1} \in L_\infty(\R^d).
\eqno(6.3)
$$
Consider the DO
$$
\begin{aligned}
{\mathcal X}&:= h b(\D) f: L_2(\R^d;\C^n) \to L_2(\R^d;\C^m),
\cr
{\rm Dom}\, {\mathcal X}& := \{ \u \in L_2(\R^d;\C^n):\ f\u \in H^1(\R^d;\C^n)\}.
\end{aligned}
$$
The operator $\mathcal X$ is closed. In $L_2(\R^d;\C^n)$, consider the selfadjoint operator ${\mathcal A}:= {\mathcal X}^*{\mathcal X}$
generated by the closed quadratic form
$a[\u,\u] := \| {\mathcal X} \u \|^2_{L_2(\R^d)}$, $\u \in {\rm Dom}\, {\mathcal X}$. Formally, we have
$$
{\mathcal A} = f(\x)^* b(\D)^*  g(\x) b(\D) f(\x),
\eqno(6.4)
$$
where $g(\x):= h(\x)^* h(\x)$. Note that the Hermitian matrix-valued function $g(\x)$ is bounded and uniformly positive definite.
Using the Fourier transformation and (6.2), (6.3), it is easy to check that
$$
\begin{aligned}
&c' \intop_{\R^d} |\D (f \u)|^2\, d\x \le a[\u,\u] \le c'' \intop_{\R^d} |\D (f \u)|^2\, d\x,
\cr
&\u \in {\rm Dom}\,{\mathcal X},\quad c'= \alpha_0 \|g^{-1}\|_{L_\infty}^{-1}, \quad c''= \alpha_1 \|g\|_{L_\infty}.
\end{aligned}
\eqno(6.5)
$$

\smallskip\noindent\textbf{6.2. Lattices in $\R^d$.} In what follows, the matrix-valued functions $h(\x)$ and $f(\x)$ are assumed to be \textit{periodic with respect to some lattice} $\Gamma \subset \R^d$.
Let $\a_1,\dots, \a_d \in \R^d$ be the basis in $\R^d$ that generates the lattice $\Gamma$:
$$
\Gamma = \{ \a \in \R^d: \ \a = \sum_{j=1}^d \nu^j \a_j,\ \nu^j \in \Z \},
$$
and let $\Omega$ be the (elementary) cell of this lattice:
$$
\Omega = \{ \x \in \R^d: \ \x = \sum_{j=1}^d \rho^j \a_j,\ 0< \rho^j <1 \}.
$$
The basis  $\b^1,\dots, \b^d$ in $\R^d$ dual to $\a_1,\dots, \a_d$ is defined by the relations $\langle \b^i, \a_j \rangle = 2 \pi \delta^i_j$.
This basis generates the \textit{lattice  $\wt{\Gamma}$ dual to} $\Gamma$:
$$
\wt{\Gamma} = \{ \b \in \R^d: \ \b = \sum_{i=1}^d \kappa_i \b^i,\ \kappa_i \in \Z \}.
$$
Let $\wt{\Omega}$ be the central Brillouin zone of the lattice $\wt{\Gamma}$:
$$
\wt{\Omega} = \{  \k \in \R^d:\ |\k| < |\k-\b|,\ 0 \ne \b \in \wt{\Gamma}\}.
$$
Note that $\wt{\Omega}$ is a fundamental domain of the lattice $\wt{\Gamma}$.
Denote  $|\Omega|= {\rm meas}\,\Omega$, $|\wt{\Omega}|= {\rm meas}\,\wt{\Omega}$. We have
$|\Omega|\, |\wt{\Omega}| = (2\pi)^{d}$. Let $r_0$ be the maximal radius of the ball containing in ${\rm clos}\, \wt{\Omega}$; then
$$
2 r_0 = \min_{0 \ne \b \in \wt{\Gamma}} |\b|.
$$

With the lattice $\Gamma$, we associate the discrete Fourier transformation $ \{ \hat{\v}_{\b}\} \mapsto \v$:
$$
\v(\x) = |\Omega|^{-1/2} \sum_{\b \in \wt{\Gamma}} \hat{\v}_{\b} e^{i \langle \b,\x\rangle},\quad \x \in \Omega,
\eqno(6.6)
$$
which is a unitary mapping of $l_2(\wt{\Gamma};\C^n)$ onto $L_2(\Omega;\C^n)$:
$$
\intop_{\Omega} |\v(\x)|^2\, d\x = \sum_{\b \in \wt{\Gamma}} |\hat{\v}_{\b}|^2.
$$

By $\wt{W}_p^s(\Omega)$ we denote the subspace in $W_p^s(\Omega)$ consisting of the functions in $W_p^s(\Omega)$
whose $\Gamma$-periodic extension to $\R^d$ belongs to $W_{p, \text{loc}}^s(\R^d)$.
If $p=2$, we use the notation $\wt{H}^s(\Omega) = \wt{W}_2^s(\Omega)$.

\smallskip\noindent\textbf{6.3. The Gelfand transformation.}
Initially, the Gelfand transformation is defined on the functions of the Schwartz class by the formula
$$
\begin{aligned}
\wt{\v}(\k,\x) = ({\mathcal U}\v)(\k,\x) = |\wt{\Omega}|^{-1/2} \sum_{\a \in \Gamma} e^{-i \langle \k, \x + \a \rangle} \v(\x+\a),
\cr
\v \in {\mathcal S}(\R^d;\C^n),\quad \x \in \Omega,\quad \k \in \wt{\Omega}.
\end{aligned}
$$
Since
$$
\intop_{\wt{\Omega}}\intop_{\Omega} | \wt{\v} (\k, \x)|^2 \,d\x \, d\k = \intop_{\R^d} |\v(\x)|^2 \, d\x,
$$
the transformation  $\mathcal U$ extends by continuity up to a \textit{unitary mapping}
$$
{\mathcal U}: L_2(\R^d;\C^n) \to \intop_{\wt{\Omega}} \oplus L_2(\Omega;\C^n)\, d\k =: {\mathcal K}.
$$
The relation $\v \in H^1(\R^d;\C^n)$ is equivalent to the fact that $\wt{\v}(\k,\cdot) \in \wt{H}^1(\Omega;\C^n)$ for almost~every $\k \in \wt{\Omega}$ and
$$
\intop_{\wt{\Omega}} \intop_{\Omega} \left( | (\D + \k)\wt{\v}(\k,\x) |^2 + |\wt{\v}(\k,\x)|^2 \right) \, d\x \,d\k < \infty.
$$
Under the Gelfand transformation $\mathcal U$, the operator of multiplication by a bounded periodic function in
$L_2(\R^d;\C^n)$  turns into
multiplication by the same function on the fibers of the direct integral $\mathcal K$. The operator $b(\D)$
applied to $\v \in H^1(\R^d;\C^n)$
turns into  the operator $b(\D + \k)$ applied to  $\wt{\v}(\k,\cdot) \in \wt{H}^1(\Omega;\C^n)$.

\section*{§7. The direct integral expansion for the operator $\mathcal A$}

\noindent\textbf{7.1. The forms $a(\k)$ and the operators ${\mathcal A}(\k)$.}
 Putting $\H = L_2(\Omega;\C^n)$, $\H_* = L_2(\Omega;\C^m)$, we consider the closed operator
${\mathcal X}(\k): \H \to \H_*$ depending on the parameter $\k \in \R^d$:
$$
{\mathcal X}(\k) = h b(\D + \k) f,\quad {\rm Dom}\, {\mathcal X}(\k) = \{ \u \in \H:\ f \u \in \wt{H}^1(\Omega;\C^n)\} =: {\mathcal D}.
$$
The selfadjoint operator ${\mathcal A}(\k)= {\mathcal X}(\k)^* {\mathcal X}(\k)$ in $\H$ is generated by the closed quadratic form
$a(\k)[\u,\u] := \| {\mathcal X}(\k) \u\|^2_{\H_*}$, $\u \in {\mathcal D}$. Using the Fourier series expansion and conditions (6.2), (6.3), it is easy to check
that
$$
\begin{aligned}
c' \intop_{\Omega} | (\D + \k)\v |^2\, d\x \le a(\k)[\u,\u] \le c'' \intop_{\Omega} | (\D + \k)\v |^2\, d\x,
\cr
\v = f \u \in \wt{H}^1(\Omega;\C^n),
\end{aligned}
\eqno(7.1)
$$
with the same constants $c'$ and $c''$ as in (6.5). From (7.1) and the compactness of the embedding  of $\wt{H}^1(\Omega;\C^n)$ in $\H$
it follows that the resolvent of the operator ${\mathcal A}(\k)$ is compact and depends on $\k$
continuously (in the operator norm).

Let
$$
\NN := {\rm Ker}\, {\mathcal A}(0) = {\rm Ker}\, {\mathcal X}(0).
\eqno(7.2)
$$
Relations (7.1) with $\k=0$ show that
$$
\NN = \{ \u \in L_2(\Omega;\C^n):\ f\u = \c,\ \c \in \C^n \},\quad {\rm dim}\, \NN =n.
\eqno(7.3)
$$

\noindent\textbf{7.2. The band functions.} The consecutive eigenvalues
$E_j(\k)$, $j\in \N$, of the operator ${\mathcal A}(\k)$ (counted with multiplicities) are called \textit{band functions}:
$$
E_1(\k) \le E_2(\k) \le \dots \le E_j(\k) \le \dots,\quad \k \in \R^d.
$$
The band functions $E_j(\k)$ are continuous and $\wt{\Gamma}$-periodic.

We put
$$
c_* = \alpha_0 \|g^{-1}\|_{L_\infty}^{-1} \|f^{-1}\|^{-2}_{L_\infty}.
\eqno(7.4)
$$
As shown in [BSu1, Chapter 2, Subsection~2.2] (by simple variational arguments), the band functions satisfy the following estimates:
$$
E_j(\k) \ge c_* |\k|^2,\quad \k \in {\rm clos}\, \wt{\Omega}, \quad j=1,\dots,n,
\eqno(7.5)
$$
$$
E_{n+1}(\k) \ge c_* r_0^2,\quad \k \in {\rm clos}\, \wt{\Omega},
\eqno(7.6)
$$
$$
E_{n+1}(0) \ge 4 c_* r_0^2.
$$

\smallskip\noindent\textbf{7.3. The direct integral expansion for the operator $\mathcal A$.}
With the help of the Gelfand transformation, the operator $\mathcal A$ is represented as
$$
{\mathcal U} {\mathcal A}{\mathcal U}^{-1} = \intop_{\wt{\Omega}} \oplus {\mathcal A}(\k) \,d\k.
\eqno(7.7)
$$
This means the following.
If $\u \in {\rm Dom}\, a$, then $\wt{\u}(\k,\cdot) \in {\mathcal D}$ for almost~every $\k \in \wt{\Omega}$, and
$$
a[\u,\u] = \intop_{\wt{\Omega}}  a(\k)[\wt{\u}(\k,\cdot),\wt{\u}(\k,\cdot) ]\,d\k.
\eqno(7.8)
$$
Conversely, if $\wt{\u} \in {\mathcal K}$ satisfies
$\wt{\u}(\k,\cdot) \in {\mathcal D}$ for a.~e. $\k \in \wt{\Omega}$ and the integral in (7.8) is finite, then
$\u \in {\rm Dom} \,a$ and (7.8) is valid.

From (7.7) it follows that the spectrum of $\mathcal A$ is the union of segments  (spectral bands) ${\rm Ran}\, E_j$, $j\in \N$.
By (7.2) and (7.3),
$$
\min_{\k} E_j(\k) = E_j(0) =0,\quad j=1,\dots,n,
$$
i.~e., the first $n$ spectral bands of $\mathcal A$ overlap and have the common bottom $\lambda_0=0$,
while the $(n+1)$-th band is separated from zero (see (7.6)).

\smallskip\noindent\textbf{7.4. Incorporation of the operators ${\mathcal A}(\k)$ into the pattern of~\S 1.}
For $\k \in \R^d$ we put $\k = t \bt$, $t = |\k|$, $\bt \in {\mathbb S}^{d-1}$, and view $t$ as the main parameter.
Then all constructions and estimates will depend on the additional parameter $\bt$, which will  often be reflected in the notation.
We have to make our estimates uniform in $\bt$.

We will apply the scheme of  \S 1 putting $\H = L_2(\Omega;\C^n)$ and $\H_* = L_2(\Omega;\C^m)$. The role of the operator $X(t)$
is played by $X(t,\bt)= {\mathcal X}(t\bt)$. Then $X(t,\bt)= X_0 + t X_1(\bt)$, ${\rm Dom}\, X(t,\bt)= {\rm Dom}\, X_0 = {\mathcal D}$,
where $X_0 = h b(\D) f$ and $X_1(\bt)= h b(\bt) f$. The role of the operator $A(t)$ is played by
$A(t,\bt)= {\mathcal A}(t\bt)$. We have $A(t,\bt)= X(t,\bt)^* X(t,\bt)$.
The kernel $\NN = {\rm Ker}\,X_0$ is described by (7.3). We have ${\rm dim}\, \NN =n$.
Together with (7.6) this shows that Condition~1.1 is satisfied.
The distance $d^0$ from the point $\lambda_0=0$ to the rest of the spectrum of ${\mathcal A}(0)$ is equal to $E_{n+1}(0)$ and satisfies the estimate
$$
d^0 \ge 4 c_* r_0^2.
\eqno(7.9)
$$
Here $c_*$ is defined by (7.4). The condition $m \ge n$ ensures that $n \le n_*$ (see [BSu1, Chapter~2, \S 3]).

   In Subsection 1.1, it was required to choose the number $\delta \in (0, d^0/8)$. Taking (7.9) into account, we fix $\delta$ as follows:
   $$
   \delta = \frac{c_* r_0^2}{4}
   = \frac{1}{4} \alpha_0 r_0^2 \|g^{-1}\|^{-1}_{L_\infty} \| f^{-1}\|^{-2}_{L_\infty}.
   \eqno(7.10)
   $$
Next, by (6.2), the operator $X_1(\bt)= h b(\bt) f$ satisfies
$$
\|X_1(\bt)\| \le \alpha_1^{1/2} \|h\|_{L_\infty} \|f\|_{L_\infty}.
\eqno(7.11)
$$
This allows us to take $t^0$ (see (1.2)) equal to the following number independent of $\bt$:
$$
\begin{aligned}
t^0 &= \delta^{1/2} \alpha_1^{-1/2} \|h\|^{-1}_{L_\infty} \|f\|^{-1}_{L_\infty}
\cr
&= \frac{r_0}{2} \alpha_0^{1/2}\alpha_1^{-1/2} \left( \|h\|_{L_\infty}\|h^{-1}\|_{L_\infty} \|f\|_{L_\infty}\|f^{-1}\|_{L_\infty}\right)^{-1}.
\end{aligned}
\eqno(7.12)
$$
Note that $t^0 \le r_0/2$. Thus, the ball $|\k|\le t^0$ lies inside $\wt{\Omega}$.
It is important that $c_*$, $\delta$, and $t^0$ (see
 (7.4), (7.10), (7.12)) are independent of $\bt$.

The variational estimates (7.5) for the eigenvalues of ${\mathcal A}(\k)$ imply that
$$
{\mathcal A}(\k) = A(t,\bt) \ge c_* t^2 I,\quad \k = t\bt \in \wt{\Omega}.
$$
Thus, Condition 1.6 is now satisfied with the constant $c_*$ defined by (7.4).
The germ ${S}(\bt)$ of the operator $A(t,\bt)$ is nondegenerate uniformly in $\bt$:  we have ${S}(\bt) \ge {c}_* I_{{\NN}}$ (cf. (1.23)).

\section*{§8. The effective characteristics \\ of the operator $\wh{\mathcal A}=b(\D)^* g(\x)b(\D)$}

\smallskip\noindent\textbf{8.1. The operator $A(t,\bt)$ in the case where $f = \1_n$.}
In the case where $f=\1_n$, the operator $A(t,\bt)$ plays a special role. In this case, all the objects will be marked by ``hat''. For instance, for the operator
$$
\wh{\mathcal A}= b(\D)^* g(\x)b(\D)
\eqno(8.1)
$$
the family $\wh{\mathcal A}(\k)$ is denoted by $\wh{A}(t,\bt)$. The kernel~(7.3) takes the form
$$
\wh{\NN} = \{\u \in L_2(\Omega;\C^n): \ \u(\x)= \c,\ \c \in \C^n\},
\eqno(8.2)
$$
i.~e., $\wh{\NN}$ consists of constant vector-valued functions.
The orthogonal projection $\wh{P}$ of the space $L_2(\Omega;\C^n)$ onto the subspace (8.2)
is the operator of averaging over the cell:
$$
\wh{P} \u = |\Omega|^{-1} \intop_\Omega \u(\x)\, d\x.
\eqno(8.3)
$$

If $f=\1_n$, the constants (7.4), (7.10), and (7.12) take the form
$$
\wh{c}_* = \alpha_0 \|g^{-1}\|^{-1}_{L_\infty},
$$
$$
\wh{\delta} = \frac{r_0^2}{4} \alpha_0  \|g^{-1}\|^{-1}_{L_\infty},
\eqno(8.4)
$$
$$
\wh{t}^{\,0} = \frac{r_0}{2} \alpha_0^{1/2}  \alpha_1^{-1/2} \|g\|^{-1/2}_{L_\infty} \|g^{-1}\|^{-1/2}_{L_\infty}.
\eqno(8.5)
$$
Inequality (7.11) takes the form
$$
\|\wh{X}_1(\bt)\| \le \alpha_1^{1/2} \|g\|^{1/2}_{L_\infty}.
\eqno(8.6)
$$

\smallskip\noindent\textbf{8.2. The germ of the operator $\wh{A}(t,\bt)$.}
According to [BSu1, Chapter~3, \S 1], the spectral germ $\wh{S}(\bt)$ of the family $\wh{A}(t,\bt)$ acting in $\wh{\NN}$ is represented as
$$
\wh{S}(\bt) = b(\bt)^* g^0 b(\bt),\quad \bt \in {\mathbb S}^{d-1},
$$
where $b(\bt)$ is the symbol of the operator $b(\D)$, and $g^0$  is the so called effective matrix. The constant positive $(m\times m)$-matrix $g^0$
is defined as follows. Suppose that a $\Gamma$-periodic $(n \times m)$-matrix-valued function $\Lambda\in \wt{H}^1(\Omega)$ is the weak solution of the problem
$$
b(\D)^* g(\x) (b(\D)\Lambda(\x)+ \1_m) =0, \quad \intop_\Omega \Lambda(\x)\, d\x =0.
\eqno(8.7)
$$
We put
$$
\wt{g}(\x) := g(\x) (b(\D)\Lambda(\x)+ \1_m).
\eqno(8.8)
$$
Then
$$
g^0 = |\Omega|^{-1} \intop_\Omega \wt{g}(\x) \, d\x.
\eqno(8.9)
$$
It turns out that the matrix $g^0$ is positive definite.

\smallskip\noindent\textbf{8.3. The effective operator.}
Consider the symbol
$$
\wh{S}(\k) := t^2 \wh{S}(\bt) = b(\k)^* g^0 b(\k), \quad \k \in \R^d.
\eqno(8.10)
$$
Expression (8.10) is the symbol of the DO
$$
\wh{\mathcal A}^0 = b(\D)^* g^0 b(\D)
\eqno(8.11)
$$
acting in $L_2(\R^d;\C^n)$ and called the \textit{effective operator} for the operator~$\wh{\mathcal A}$.

Let $\wh{\mathcal A}^0(\k)$ be the operator family in $L_2(\Omega;\C^n)$ corresponding to $\wh{\A}^0$.
Then $\wh{\mathcal A}^0(\k)$ is given by the expression $b(\D+\k)^* g^0 b(\D+\k)$ with periodic boundary conditions.
By (8.3) and (8.10), we have
$$
\wh{S}(\k) \wh{P} = \wh{\mathcal A}^0(\k) \wh{P}.
\eqno(8.12)
$$

\smallskip\noindent\textbf{8.4. Properties of the effective matrix.}
The following properties of the matrix $g^0$ were checked in [BSu1, Chapter~3, Theorem~1.5].

 \smallskip\noindent\textbf{Proposition 8.1.} \textit{The effective matrix satisfies the estimates}
 $$
 \underline{g} \le g^0 \le \overline{g},
 \eqno(8.13)
 $$
\textit{where}
$$
\overline{g}:= |\Omega|^{-1} \intop_\Omega g(\x)\, d\x,\quad
\underline{g}:= \biggl( |\Omega|^{-1} \intop_\Omega g(\x)\, d\x \biggr)^{-1}.
$$
\textit{If $m=n$, then} $g^0 = \underline{g}$.

\smallskip
For specific DOs, estimates (8.13) are known in homogenization theory as the Voigt-Reuss bracketing.
Now we distinguish the cases where one of the inequalities in (8.13) becomes an identity.
The following statements were obtained in [BSu1, Chapter~3, Propositions 1.6, 1.7].

 \smallskip\noindent\textbf{Proposition 8.2.}
 \textit{The identity $g^0 = \overline{g}$ is equivalent to the relations}
  $$
  b(\D)^* {\mathbf g}_k(\x) =0,\quad k=1,\dots,m,
  \eqno(8.14)
  $$
\textit{where ${\mathbf g}_k(\x)$, $k=1,\dots, m$, are the columns of the matrix $g(\x)$.}

 \smallskip\noindent\textbf{Proposition 8.3.} \textit{The identity $g^0 = \underline{g}$ is
 equivalent to the representations}
  $$
  {\mathbf l}_k(\x) =  {\mathbf l}_k^0 + b(\D) {\mathbf w}_k(\x),\quad {\mathbf l}_k^0 \in \C^m, \quad {\mathbf w}_k \in \wt{H}^1(\Omega;\C^n), \quad k=1,\dots,m,
  \eqno(8.15)
  $$
\textit{where ${\mathbf l}_k(\x)$, $k=1,\dots, m$, are the columns of the matrix $g(\x)^{-1}$.}

\smallskip\noindent\textbf{8.5. The analytic branches of eigenvalues and eigenvectors.}
The analytic (in $t$) branches of the eigenvalues $\wh{\lambda}_l(t,\bt)$ and the analytic branches of
the eigenvectors $\wh{\varphi}_l(t,\bt)$ of $\wh{A}(t,\bt)$ admit the power series expansions of the form (1.6) and (1.7)  with the coefficients depending on $\bt$:
$$
\wh{\lambda}_l(t,\bt) = \wh{\gamma}_l(\bt) t^2 + \wh{\mu}_l(\bt) t^3 + \dots,\quad l=1,\dots,n,
\eqno(8.16)
$$
$$
\wh{\varphi}_l(t,\bt) = \wh{\omega}_l(\bt)  + t \wh{\psi}_l^{(1)}(\bt) + \dots,\quad l=1,\dots,n.
\eqno(8.17)
$$
(However,  we do not control the interval of convergence $t=|\k|\le t_*(\bt)$.)
According to (1.9), the numbers $\wh{\gamma}_l(\bt)$ and the elements $\wh{\omega}_l(\bt)$
are eigenvalues and eigenvectors of the germ:
$$
b(\bt)^* g^0 b(\bt) \wh{\omega}_l(\bt) = \wh{\gamma}_l(\bt)\wh{\omega}_l(\bt),\quad l=1,\dots,n.
$$

\smallskip\noindent\textbf{8.6. The operator $\wh{N}(\bt)$.}
We need to describe the operator~$N$ (that in abstract terms is defined in Theorem~1.4).
According to [BSu3,\S 4], for the family $\wh{A}(t,\bt)$ this operator takes the form
$$
\wh{N}(\bt) = b(\bt)^* L(\bt) b(\bt) \wh{P},
\eqno(8.18)
$$
where the  $(m\times m)$-matrix $L(\bt)$ is given by
$$
L(\bt) = |\Omega|^{-1} \intop_\Omega \left( \Lambda(\x)^* b(\bt)^* \wt{g}(\x) + \wt{g}(\x)^* b(\bt) \Lambda(\x)\right) \,d\x.
\eqno(8.19)
$$
Here $\Lambda(\x)$ is the $\Gamma$-periodic solution of problem (8.7), and $\wt{g}(\x)$ is given by~(8.8).

Observe that  $L(\k):=t L(\bt)$, $\k\in \R^d$,  is a Hermitian matrix-valued function first order homogeneous in $\k$.
 We put $\wh{N}(\k):=t^3 \wh{N}(\bt)$, $\k\in \R^d$. Then
$\wh{N}(\k) = b(\k)^* L(\k) b(\k) \wh{P}$. The matrix-valued function $b(\k)^* L(\k)b(\k)$ is a homogeneous third order polynomial of $\k \in \R^d$.
Therefore, either $\wh{N}(\bt)=0$ for all $\bt \in {\mathbb S}^{d-1}$, or $\wh{N}(\bt)\ne 0$  at most points $\bt$
(except for the zeroes of this polynomial).

Some cases where the operator (8.18) is equal to zero were distinguished in [BSu3, \S 4].

\smallskip\noindent\textbf{Proposition 8.4.}
\textit{Suppose that at least one of the following conditions is fulfilled}:

\noindent
$1^\circ$. \textit{The operator $\wh{\mathcal A}$ has the form
 $\wh{\mathcal A}= \D^* g(\x)\D$, where $g(\x)$ is a symmetric matrix with real entries}.

\noindent
$2^\circ$. \textit{Relations} (8.14) \textit{are satisfied, i.~e.} $g^0 = \overline{g}$.

\noindent
$3^\circ$. \textit{Relations} (8.15) \textit{are satisfied, i.~e.} $g^0 = \underline{g}$.
(\textit{In particular, this is true if} $m=n$.)

\noindent\textit{Then $\wh{N}(\bt)=0$ for all} $\bt \in {\mathbb S}^{d-1}$.

\smallskip
On the other hand, there are examples (see [BSu3, Subsections 10.4, 13.2, 14.6])
showing that, in general, the operator $\wh{N}(\bt)$ is not equal to zero for
the scalar elliptic operator $\D^* g(\x) \D$, where $g(\x)$ is a Hermitian matrix with complex entries,
as well as for matrix operators even with real-valued coefficients; see also Example 8.7 below.

Recall (see Remark 1.5) that $\wh{N}(\bt)= \wh{N}_0(\bt) + \wh{N}_*(\bt)$,
where the operator $\wh{N}_0(\bt)$ is diagonal in the basis $\wh{\omega}_1(\bt),\dots, \wh{\omega}_n(\bt)$ (see (8.17)),
while the diagonal elements of the operator $\wh{N}_*(\bt)$ are equal to zero. We have
$$
(\wh{N}(\bt) \wh{\omega}_l(\bt), \wh{\omega}_l(\bt))_{L_2(\Omega)} =(\wh{N}_0(\bt) \wh{\omega}_l(\bt), \wh{\omega}_l(\bt))_{L_2(\Omega)} = \wh{\mu}_l(\bt),
\quad l=1,\dots,n.
\eqno(8.20)
$$

In [BSu3, Subsection 4.3], the following argument is given. Suppose that $b(\bt)$ and $g(\x)$ are matrices with \textit{real entries}.
Then the matrix $\Lambda(\x)$ (see (8.7)) has purely imaginary entries, while $\wt{g}(\x)$ and $g^0$ are matrices with real entries.
In this case $L(\bt)$ (see (8.19)) and $b(\bt)^* L(\bt) b(\bt)$ are Hermitian matrices with purely imaginary entries.
Hence, for any \textit{real} vector ${\mathbf q} \in \wh{\NN}$ we have
$(\wh{N}(\bt) {\mathbf q}, {\mathbf q})=0$. If the analytic branches of the eigenvalues $\wh{\lambda}_l(t,\bt)$
and the analytic branches of the eigenvectors $\wh{\varphi}_l(t,\bt)$ of $\wh{A}(t,\bt)$ can be chosen so that the vectors $\wh{\omega}_1(\bt), \dots, \wh{\omega}_n(\bt)$ are real,
then, by  (8.20), we have  $\wh{\mu}_l(\bt)=0$, $l=1,\dots,n$, i.~e., $\wh{N}_0(\bt)=0$. We arrive at the following statement.

\smallskip\noindent\textbf{Proposition 8.5.} \textit{Suppose that $b(\bt)$ and $g(\x)$ have real entries. Suppose that in the expansions} (8.17)
\textit{for the analytic branches of the eigenvectors of $\wh{A}(t,\bt)$ the ``embrios''  $\wh{\omega}_l(\bt)$, $l=1,\dots,n$, can be chosen to be real.}
\textit{Then in} (8.16) \textit{we have $\wh{\mu}_l(\bt)=0$, $l=1,\dots,n$, i.~e., $\wh{N}_0(\bt)=0$ for all $\bt \in {\mathbb S}^{d-1}$.}

\smallskip
In the  ``real'' case under consideration, the germ $\wh{S}(\bt)$ is a symmetric matrix with real entries.
Clearly, if the eigenvalue $\wh{\gamma}_j(\bt)$ of the germ is simple,
then the embrio $\wh{\omega}_j(\bt)$ is defined uniquely up to a phase factor, and we can always choose
$\wh{\omega}_j(\bt)$ to be real.  We arrive at the following corollary.

\smallskip\noindent\textbf{Corollary 8.6.} \textit{Suppose that $b(\bt)$ and $g(\x)$ have real entries. Suppose that the spectrum of the germ $\wh{S}(\bt)$ is simple. Then $\wh{N}_0(\bt)=0$ for all $\bt \in {\mathbb S}^{d-1}$.}

\smallskip
However, as is seen from Example 8.7 considered below, even in the ``real'' case it is not always possible to choose the vectors $\wh{\omega}_l(\bt)$ to be real.
Moreover, it can happen that $\wh{N}_0(\bt) \ne 0$ at some isolated points $\bt$.

\smallskip\noindent\textbf{8.7. Multiplicities of the eigenvalues of the germ. Example.}
Considerations of this subsection concern the case where $n\ge 2$.
Now we return to the notation of  \S 2, tracing the multiplicities of the eigenvalues of the spectral germ $\wh{S}(\bt)$.
In general, the number $p(\bt)$ of different eigenvalues $\wh{\gamma}^\circ_1(\bt), \dots \wh{\gamma}^\circ_{p(\bt)}(\bt)$ of the spectral germ $\wh{S}(\bt)$
and their multiplicities $k_1(\bt),\dots, k_{p(\bt)}(\bt)$ depend on the parameter $\bt \in {\mathbb S}^{d-1}$.
For a fixed $\bt$ denote by $\wh{P}_j(\bt)$ the orthogonal projection of $L_2(\Omega;\C^n)$ onto the eigenspace
of the germ $\wh{S}(\bt)$ corresponding to the eigenvalue $\wh{\gamma}^\circ_j(\bt)$.
According to (3.20), the operators $\wh{N}_0(\bt)$ and $\wh{N}_*(\bt)$ admit the following invariant representations:
$$
\wh{N}_0(\bt) = \sum_{j=1}^{p(\bt)} \wh{P}_j(\bt) \wh{N}(\bt) \wh{P}_j(\bt),
\eqno(8.21)
$$
$$
\wh{N}_*(\bt) = \sum_{1\le j,l  \le p(\bt):\, j\ne l} \wh{P}_j(\bt) \wh{N}(\bt) \wh{P}_l(\bt).
\eqno(8.22)
$$

In conclusion of this section, we consider the example which shows that for matrix operators even with real-valued coefficients the eigenvalues of the germ may be multiple, and the coefficients $\wh{\mu}_l(\bt)$ in  (8.16)
may be nonzero.

\smallskip\noindent\textbf{Example 8.7.}
In this example, the matrices $b(\bt)$ and $g(\x)$ have real entries. Let $d=2$, $n=2$, and $m=3$.
For simplicity, assume that $\Gamma= (2\pi \Z)^2$.
Suppose that the operator $b(\D)$ and the matrix $g(\x)$ are given by
$$
b(\D) = \begin{pmatrix}
D_1 & 0 \cr
\frac{1}{2} D_2 & \frac{1}{2} D_1 \cr
0 & D_2
\end{pmatrix},
\quad
g(\x) = \begin{pmatrix}
1 & 0 & 0 \cr
0  &  g_2(x_1) & 0 \cr
0 & 0 & g_3(x_1)
\end{pmatrix},
$$
where $g_2(x_1)$ and $g_3(x_1)$ are $(2\pi)$-periodic bounded and positive definite functions of $x_1$, and  $\overline{g_3}=1$.
It is easy to find the  $\Gamma$-periodic solution of problem (8.7):
$$
\Lambda(\x) =
\begin{pmatrix}
0 & 0 & 0 \cr
0  &  \Lambda_{22}(x_1) & 0
\end{pmatrix}.
$$
Here $\Lambda_{22}(x_1)$ is the $(2\pi)$-periodic solution of the problem
$$
\frac{1}{2} D_1 \Lambda_{22}(x_1) +1 = \underline{g_2} (g_2(x_1))^{-1}, \quad \intop_0^{2\pi} \Lambda_{22}(x_1)\, dx_1 =0.
$$
Obviously, $\Lambda_{22}(x_1)$ is purely imaginary. Then
$\wt{g}(\x) = {\rm diag}\, \{ 1, \underline{g_2}, g_3(x_1)\}$, and
$g^0 = {\rm diag}\, \{ 1, \underline{g_2}, 1 \}$. The spectral germ $\wh{S}(\bt) = b(\bt)^* g^0 b(\bt)$ is given by
$$
\wh{S}(\bt)=
\begin{pmatrix}
\theta_1^2 + \frac{1}{4} \theta_2^2 \underline{g_2} & \frac{1}{4} \theta_1 \theta_2 \underline{g_2} \cr
 \frac{1}{4} \theta_1 \theta_2 \underline{g_2} & \frac{1}{4} \theta_1^2 \underline{g_2} + \theta_2^2
\end{pmatrix},\quad \bt = (\theta_1,\theta_2) \in {\mathbb S}^1.
\eqno(8.23)
$$
It is easily seen that the matrix (8.23) has a multiple eigenvalue (for some $\bt$) only if  $\underline{g_2}=4$.

  So, let $\underline{g_2}=4$.  Then the eigenvalues of the germ
$$
\wh{S}(\bt)=
\begin{pmatrix}
1  &  \theta_1 \theta_2  \cr
 \theta_1 \theta_2 & 1
\end{pmatrix}
$$
are $\wh{\gamma}_1(\bt) = 1+ \theta_1 \theta_2$ and $\wh{\gamma}_2(\bt) = 1- \theta_1 \theta_2$.
They coincide at four points $\bt^{(1)}=(0,1)$, $\bt^{(2)}=(0,-1)$, $\bt^{(3)}=(1,0)$, $\bt^{(4)}=(-1,0)$.

Next, we calculate the $(3\times 3)$-matrix $L(\bt)$ (see (8.19)):
$$
L(\bt) = \begin{pmatrix}
0 & 0 & 0 \cr
0 & 0 & \theta_2 \overline{\Lambda_{22}^* g_3}
\cr
0 & \theta_2 \overline{\Lambda_{22} g_3}  & 0
\end{pmatrix}.
$$
Hence,
$$
\wh{N}(\bt) = b(\bt)^* L(\bt) b(\bt)
= \frac{1}{2} \theta_2^3 \begin{pmatrix}
0 &  \overline{\Lambda_{22}^* g_3}
\cr
\overline{\Lambda_{22}  g_3} & 0
\end{pmatrix}.
$$
Below we assume that  $\overline{\Lambda_{22} g_3} \ne 0$.
(It is easy to give a concrete example: if $g_2(x_1)= 4 (1+ \frac{1}{2} \sin x_1)^{-1}$ and
$g_3(x_1)= 1+ \frac{1}{2} \cos x_1$, all the conditions are fulfilled.)

 For $\bt \ne \bt^{(j)}$, $j=1,2,3,4,$ we have $\wh{\gamma}_1(\bt) \ne \wh{\gamma}_2(\bt)$
and then $\wh{N}(\bt) = \wh{N}_*(\bt) \ne 0$. At the points $\bt^{(1)}$ and $\bt^{(2)}$ we have
$$
\wh{\gamma}_1(\bt^{(j)}) = \wh{\gamma}_2(\bt^{(j)}) =1, \quad \wh{N}(\bt^{(j)}) = \wh{N}_0(\bt^{(j)}) \ne 0,\quad j=1,2.
$$
Obviously, the numbers $\pm \mu$, where $\mu = \frac{1}{2}|\overline{\Lambda_{22} g_3}|$,
 are the eigenvalues of the operator $\wh{N}(\bt^{(j)})$ for $j=1,2$.
 In the exansions (8.16) there are nonzero coefficients at $t^3$:
$$
\wh{\lambda}_1(t, \bt^{(j)})= t^2 + \mu t^3 +\dots, \quad
\wh{\lambda}_2(t, \bt^{(j)})= t^2 - \mu t^3 +\dots,\quad j=1,2.
$$
In this case, the embrios $\wh{\omega}_1(\bt^{(j)})$, $\wh{\omega}_2(\bt^{(j)})$ in the expansions  (8.17) can not be real
(see Proposition 8.5).

 At the points $\bt^{(3)}$ and $\bt^{(4)}$ the situation is different. We have
$$
\wh{\gamma}_1(\bt^{(j)}) = \wh{\gamma}_2(\bt^{(j)}) =1, \quad \wh{N}(\bt^{(j)}) = 0,\quad j=3,4.
$$

This example also shows that, though the operator $\wh{N}(\bt)$ is always continuous in $\bt$
(it is a polynomial of the third degree),
its  ``blocks'' $\wh{N}_0(\bt)$ and $\wh{N}_*(\bt)$ can be discontinuous:
at the points where the branches of the eigenvalues of the germ intersect,
$\wh{N}_0(\bt)$ and $\wh{N}_*(\bt)$ may have jumps. Moreover, it may
happen that $\wh{N}_0(\bt)$ is not equal to zero only at some isolated points.

\section*{§9. Approximation of the smoothed operator $e^{-i \eps^{-2} \tau \wh{\mathcal A}(\k)}$}

\smallskip\noindent\textbf{9.1. Approximation of the smoothed operator $e^{-i \eps^{-2} \tau \wh{\mathcal A}(\k)}$ for \hbox{$|\k| \le \wh{t}^0$}.}
Consider the operator ${\mathcal H}_0= - \Delta$ in $L_2(\R^d;\C^n)$.
Under the Gelfand transformation, this operator expands in the direct integral of the operators
${\mathcal H}_0(\k)$ acting in $L_2(\Omega;\C^n)$. The operator ${\mathcal H}_0(\k)$ is given by the differential expression $|\D + \k|^2$ with periodic boundary conditions.
Denote
$$
{\mathcal R}(\k,\eps):= \eps^2 ({\mathcal H}_0(\k) + \eps^2 I)^{-1}.
\eqno(9.1)
$$
Obviously,
$$
{\mathcal R}(\k,\eps)^{s/2} \wh{P} = \eps^s (t^2 + \eps^2 )^{-s/2} \wh{P}, \quad s>0.
\eqno(9.2)
$$

We will apply theorems of \S 4 to the operator $\wh{A}(t,\bt)= \wh{\mathcal A}(\k)$. We start with Theorem~4.1.
First, we need to specify the constants. By (8.4) and (8.6),  instead of the precise values of the constants $\wh{C}_1(\bt)= \beta_1 \wh{\delta}^{-1/2} \| \wh{X}_1(\bt)\|$,
$\wh{C}_2(\bt) = \beta_2 \wh{\delta}^{-1/2} \| \wh{X}_1(\bt)\|^3$ (which depend on $\bt$) we can take
$$
\begin{aligned}
\wh{C}_1 = 2 \beta_1 r_0^{-1} \alpha_1^{1/2} \alpha_0^{-1/2} \|g\|^{1/2}_{L_\infty} \|g^{-1}\|^{1/2}_{L_\infty},
\cr
\wh{C}_2 = 2 \beta_2 r_0^{-1} \alpha_1^{3/2} \alpha_0^{-1/2} \|g\|^{3/2}_{L_\infty} \|g^{-1}\|^{1/2}_{L_\infty}.
\end{aligned}
\eqno(9.3)
$$

Combining (4.1), (8.12), and (9.2), we arrive at the inequality
$$
\begin{aligned}
&\| \left( e^{- i \eps^{-2} \tau \wh{\mathcal A}(\k)}- e^{- i \eps^{-2} \tau \wh{\mathcal A}^0(\k)} \right) {\mathcal R}(\k,\eps)^{3/2} \wh{P}\|_{L_2(\Omega) \to L_2(\Omega)}
\le (\wh{C}_1 + \wh{C}_2 |\tau|) \eps,
\cr
&\tau \in \R,\quad \eps >0,\quad |\k| \le \wh{t}^{\,0}.
\end{aligned}
\eqno(9.4)
$$

\subsection*{9.2. Estimate for $|\k| > \wh{t}^{\,0}$.}
For $\k \in \wt{\Omega}$ and $|\k| > \wh{t}^{\,0}$ estimates are trivial. By (9.2), we have
$$
\| {\mathcal R}(\k,\eps)^{1/2} \wh{P} \|_{L_2(\Omega) \to L_2(\Omega)} \le (\wh{t}^{\,0})^{-1} \eps,\quad \eps>0,\quad \k \in \wt{\Omega}, \quad |\k| > \wh{t}^{\,0}.
\eqno(9.5)
$$
Therefore,
$$
\begin{aligned}
&\| \left( e^{- i \eps^{-2} \tau \wh{\mathcal A}(\k)}- e^{- i \eps^{-2} \tau \wh{\mathcal A}^0(\k)} \right) {\mathcal R}(\k,\eps)^{1/2} \wh{P}\|_{L_2(\Omega) \to L_2(\Omega)}
\le 2 (\wh{t}^{\,0})^{-1} \eps,
\cr
&\tau \in \R,\quad \eps >0,\quad \k \in \wt{\Omega}, \quad |\k| > \wh{t}^{\,0}.
\end{aligned}
\eqno(9.6)
$$
Note that here the smoothing operator is ${\mathcal R}(\k,\eps)^{1/2}$ (i.~e.,  $s=1$).
Of course, the left-hand side of (9.4) also satisfies the same estimate.

From (9.4) and (9.6), using expressions for $\wh{t}^{\,0}$ and $\wh{C}_1$ (see (8.5) and (9.3)), we obtain
$$
\begin{aligned}
&\| \left( e^{- i \eps^{-2} \tau \wh{\mathcal A}(\k)}- e^{- i \eps^{-2} \tau \wh{\mathcal A}^0(\k)} \right) {\mathcal R}(\k,\eps)^{3/2} \wh{P}\|_{L_2(\Omega) \to L_2(\Omega)}
\le (\wh{C}^*_1 + \wh{C}_2 |\tau|) \eps,
\cr
&\tau \in \R,\quad \eps >0,\quad \k \in \wt{\Omega},
\end{aligned}
\eqno(9.7)
$$
where
$\wh{C}^*_1 = \max \{ \wh{C}_1, 2(\wh{t}^{\,0})^{-1}\} = \beta_1^* r_0^{-1} \alpha_1^{1/2} \alpha_0^{-1/2} \|g\|^{1/2}_{L_\infty} \|g^{-1}\|^{1/2}_{L_\infty}$.

\subsection*{9.3. Removal of the operator $\wh{P}$}
Now we show that, up to an admissible error,
the projection $\wh{P}$ can be replaced by the identity operator  under the norm sign in (9.7).
For this, we estimate the norm of the operator ${\mathcal R}(\k,\eps)^{s/2} (I - \wh{P})$.
Under the discrete Fourier transformation (see (6.6)),  the operator ${\mathcal R}(\k,\eps)^{s/2}$ turns into multiplication of the Fourier coefficients by the symbol
$\eps^{s} (|\b+\k|^2 + \eps^2)^{-s/2}$. The operator $I - \wh{P}$ makes the zero Fourier coefficient
equal to zero. Therefore,
$$
\begin{aligned}
\| {\mathcal R}(\k,\eps)^{s/2} (I - \wh{P}) \|_{L_2(\Omega) \to L_2(\Omega)} \le
\sup_{0 \ne \b \in \wt{\Gamma}} \eps^{s} (|\b+\k|^2 + \eps^2)^{-s/2} \le r_0^{-s} \eps^s,
\cr
\eps>0,\quad \k \in \wt{\Omega}.
\end{aligned}
\eqno(9.8)
$$
Hence,
$$
\begin{aligned}
&\| \left( e^{- i \eps^{-2} \tau \wh{\mathcal A}(\k)}- e^{- i \eps^{-2} \tau \wh{\mathcal A}^0(\k)} \right) {\mathcal R}(\k,\eps)^{1/2} (I- \wh{P})\|_{L_2(\Omega) \to L_2(\Omega)}
\le 2 r_0^{-1} \eps,
\cr
&\tau \in \R,\quad \eps >0,\quad \k \in \wt{\Omega}.
\end{aligned}
\eqno(9.9)
$$
Note that here the smoothing operator is ${\mathcal R}(\k,\eps)^{1/2}$ (i.~e., $s=1$).

Finally, from (9.7) and (9.9), using the obvious inequality $\|{\mathcal R}(\k,\eps)\| \le 1$ and expressions for the constants,
we obtain the following result which has been proved before in  [BSu5, Theorem 7.1].

\smallskip\noindent\textbf{Theorem 9.1.} \textit{For $\tau \in \R$, $\eps >0$, and $\k \in \wt{\Omega}$ we have}
$$
\| \left( e^{- i \eps^{-2} \tau \wh{\mathcal A}(\k)}- e^{- i \eps^{-2} \tau \wh{\mathcal A}^0(\k)} \right) {\mathcal R}(\k,\eps)^{3/2} \|_{L_2(\Omega) \to L_2(\Omega)}
\le (\wh{\mathcal C}_1 + \wh{\mathcal{C}}_2 |\tau|) \eps,
$$
\textit{where}
$$
\begin{aligned}
\wh{\mathcal C}_1 &= \wh{\beta}_1 r_0^{-1} \alpha_1^{1/2} \alpha_0^{-1/2} \|g\|^{1/2}_{L_\infty} \|g^{-1}\|^{1/2}_{L_\infty},
\cr
 \wh{\mathcal C}_2 &= \wh{C}_2 = 2 \beta_2 r_0^{-1} \alpha_1^{3/2} \alpha_0^{-1/2} \|g\|^{3/2}_{L_\infty} \|g^{-1}\|^{1/2}_{L_\infty}.
\end{aligned}
\eqno(9.10)
$$

\subsection*{9.4. Refinement of approximation of the smoothed operator $e^{-i \eps^{-2} \tau \wh{\mathcal A}(\k)}$
in the case where $\wh{N}(\bt)=0$}
Now we apply Theorem~4.2, assuming that  $\wh{N}(\bt)=0$ for all $\bt \in {\mathbb S}^{d-1}$.
Taking (8.12) and (9.2) into account, we have
$$
\begin{aligned}
&\| \left( e^{- i \eps^{-2} \tau \wh{\mathcal A}(\k)}- e^{- i \eps^{-2} \tau \wh{\mathcal A}^0(\k)} \right) {\mathcal R}(\k,\eps) \wh{P}\|_{L_2(\Omega) \to L_2(\Omega)}
\le  (\wh{C}'_1 + \wh{C}_5 |\tau|) \eps,
\cr
&\tau \in \R,\quad \eps >0,\quad |\k| \le \wh{t}^{\,0}.
\end{aligned}
$$
Here $\wh{C}'_1  = \max\{2, \wh{C}_1\}$, and the constant $\wh{C}_5$ is given by
$\wh{C}_5 = 4 \beta_5 r_0^{-2} \alpha_1^2 \alpha_0^{-1} \|g\|^2_{L_\infty} \| g^{-1}\|_{\infty}$.

Together with (9.6) and (9.9) this implies the following result.

\smallskip\noindent\textbf{Theorem 9.2.} \textit{Let $\wh{N}(\bt)$ be the operator defined by} (8.18), (8.19).
\textit{Suppose that $\wh{N}(\bt)=0$ for all $\bt \in {\mathbb S}^{d-1}$.
Then for $\tau \in \R$, $\eps >0$, and $\k \in \wt{\Omega}$ we have}
$$
\| \left( e^{- i \eps^{-2} \tau \wh{\mathcal A}(\k)}- e^{- i \eps^{-2} \tau \wh{\mathcal A}^0(\k)} \right) {\mathcal R}(\k,\eps) \|_{L_2(\Omega) \to L_2(\Omega)}
\le (\wh{\mathcal C}_3 + \wh{\mathcal{C}}_4 |\tau|) \eps,
$$
\textit{where
$\wh{\mathcal C}_3 = \max \{ 2+ 2r_0^{-1}, \wh{\mathcal C}_1\}$, $\wh{\mathcal C}_1$ is defined by} (9.10), \textit{and}
$\wh{\mathcal{C}}_4 = \wh{C}_5=4 \beta_5 r_0^{-2} \alpha_1^2 \alpha_0^{-1} \|g\|^2_{L_\infty} \| g^{-1}\|_{\infty}$.

\smallskip
Recall that some sufficient conditions ensuring that
$\wh{N}(\bt)=0$ for all $\bt \in {\mathbb S}^{d-1}$ are given in Proposition 8.4.

\subsection*{9.5. Refinement of approximation of the smoothed operator $e^{-i \eps^{-2} \tau \wh{\mathcal A}(\k)}$
in the case where $\wh{N}_0(\bt)=0$}
Now, we reject the assumption of Theorem 9.2, but instead we assume that $\wh{N}_0(\bt)=0$ for all $\bt$.
We may also assume that $\wh{N}(\bt)= \wh{N}_*(\bt) \ne 0$ for some $\bt$, and then at most points $\bt$
(otherwise, one can apply Theorem~9.2.) We would like to apply the ``abstract'' result, namely, Theorem~4.3.
However, there is an additional difficulty related to the fact that
the multiplicities of the eigenvalues of the germ $\wh{S}(\bt)$ may change at some points $\bt$.
 Near such points the distance between some pair of different eigenvalues
 tends to zero, and we are not able to choose the parameters (2.3) and (2.4) to be independent of $\bt$.
 Therefore, we are forced to impose an additional condition.
 We have to take care only about those pairs of eigenvalues for which the corresponding term in (8.22) is not zero.
 Since the number of different eigenvalues of the germ and their multiplicities may depend on $\bt$,
now it is more convenient to use the initial enumeration of the eigenvalues
$\wh{\gamma}_1(\bt),\dots, \wh{\gamma}_n(\bt)$ of $\wh{S}(\bt)$
(each eigenvalue is repeated according to its multiplicity). We enumerate them in the nondecreasing order:
$$
\wh{\gamma}_1(\bt) \le \wh{\gamma}_2(\bt) \le \dots \le \wh{\gamma}_n(\bt).
$$
 For each $\bt$ denote by $\wh{P}^{(k)}(\bt)$ the orthogonal projection of $L_2(\Omega;\C^n)$ onto the eigenspace
 of $\wh{S}(\bt)$ corresponding to the eigenvalue $\wh{\gamma}_k(\bt)$.
      Clearly, for each $\bt$ the operator $\wh{P}^{(k)}(\bt)$ coincides with one of the projections $\wh{P}_j(\bt)$
      introduced in Subsection~8.7 (but the number $j$
      may depend on $\bt$).

\smallskip\noindent\textbf{Condition 9.3.}
{$1^\circ$.} \textit{The operator $\wh{N}_0(\bt)$ defined by} (8.21) \textit{is equal to zero}: $\wh{N}_0(\bt)=0$ \textit{for all}
$\bt \in {\mathbb S}^{d-1}$.

\noindent{$2^\circ$.} \textit{For any pair of indices $(k,r)$, $1\le k,r \le n$, $k\ne r$,
such that $\wh{\gamma}_k(\bt_0)=\wh{\gamma}_r(\bt_0)$ for some $\bt_0 \in {\mathbb S}^{d-1}$, we have}
$\wh{P}^{(k)}(\bt) \wh{N}(\bt) \wh{P}^{(r)}(\bt) =0$ \textit{for any} $\bt\in {\mathbb S}^{d-1}$.

    \smallskip
    Note that $\wh{P}^{(k)}(\bt_0)=\wh{P}^{(r)}(\bt_0)$
    at the points $\bt_0$ such that $\wh{\gamma}_k(\bt_0)=\wh{\gamma}_r(\bt_0)$.
    Therefore, the identity $\wh{P}^{(k)}(\bt_0) \wh{N}(\bt_0) \wh{P}^{(r)}(\bt_0) =0$
    holds automatically in virtue of condition $1^\circ$. Condition $2^\circ$ can be reformulated
    as follows: we assume that, for the ``blocks''
    $\wh{P}^{(k)}(\bt) \wh{N}(\bt) \wh{P}^{(r)}(\bt)$ of the operator $\wh{N}(\bt)$ that are not
    identically zero, the corresponding branches of the eigenvalues $\wh{\gamma}_k(\bt)$ and $\wh{\gamma}_r(\bt)$
    do not intersect.

    Obviously, Condition 9.3 is ensured by the following more restrictive condition.

\smallskip\noindent\textbf{Condition 9.4.}
{$1^\circ$.} \textit{The operator $\wh{N}_0(\bt)$ defined by} (8.21) \textit{is equal to zero}:
$\wh{N}_0(\bt)=0$ \textit{for all} $\bt \in {\mathbb S}^{d-1}$.

\noindent$2^\circ$. \textit{Assume that the number $p$ of different eigenvalues of the spectral germ $\wh{S}(\bt)$
does not depend on $\bt \in {\mathbb S}^{d-1}$. Denote different eigenvalues of the germ enumerated in the increasing order
 by $\wh{\gamma}_1^\circ(\bt),\dots, \wh{\gamma}_p^\circ(\bt)$, and assume that their multiplicities $k_1,\dots, k_p$
 do not depend on $\bt \in {\mathbb S}^{d-1}$.}

\smallskip\noindent\textbf{Remark 9.5.}
Assumption $2^\circ$ of Condition~9.4 is a fortiori satisfied, if the spectrum of the germ $\wh{S}(\bt)$
is simple for any $\bt \in {\mathbb S}^{d-1}$.

\smallskip

So, we assume that Condition~9.3 is satisfied.
We are interested only in the pairs of indices from the set
$$
\wh{\mathcal K}  := \{(k,r): \ 1\le k,r \le n,\ k\ne r,\ \wh{P}^{(k)}(\bt) \wh{N}(\bt) \wh{P}^{(r)}(\bt) \not\equiv  0 \}.
$$
Denote (cf. (2.3))
$$
\wh{c}^{\,\circ}_{kr}(\bt) := \min \{ \wh{c}_*,  {n}^{-1} |\wh{\gamma}_k(\bt) - \wh{\gamma}_{r}(\bt)|\},\quad (k, r) \in \wh{\mathcal K}.
$$
Since $\wh{S}(\bt)$ is continuous in $\bt\in {\mathbb S}^{d-1}$ (this is a polynomial of the second degree),
then the perturbation theory of discrete spectrum implies that the functions
$\wh{\gamma}_j(\bt)$ are continuous  on the sphere ${\mathbb S}^{d-1}$.
By Condition~9.3($2^\circ$), for $(k,r)\in \wh{\mathcal K}$ we have $|\wh{\gamma}_k (\bt) - \wh{\gamma}_{r}(\bt)| >0$ for all $\bt$,
whence
$$
\wh{c}^{\,\circ}_{kr} := \min_{\bt \in {\mathbb S}^{d-1}} \wh{c}_{kr}^{\,\circ}(\bt) >0, \quad (k,r) \in \wh{\mathcal K}.
$$
We put
$$
\wh{c}^{\,\circ} := \min_{(k,r) \in \wh{\mathcal K}} \wh{c}_{kr}^{\,\circ}.
\eqno(9.11)
$$
Clearly, the number (9.11) is a realization of (3.33) chosen independently of $\bt$.

Under Condition 9.3, the number $t^{00}$ subject to (3.34) also can be chosen independently of $\bt \in {\mathbb S}^{d-1}$.
Taking (8.4) and (8.6) into account, we put
$$
\wh{t}^{\,00} = (8 \beta_2)^{-1} r_0 \alpha_1^{-3/2} \alpha_0^{1/2}  \|g\|_{L_\infty}^{-3/2} \|g^{-1} \|_{L_\infty}^{-1/2} \wh{c}^{\,\circ},
\eqno(9.12)
$$
where $\wh{c}^{\,\circ}$ is defined by (9.11).
(The condition $\wh{t}^{\,00} \le \wh{t}^{\,0}$ is valid automatically since $\wh{c}^{\,\circ} \le \| \wh{S} (\bt)\| \le \alpha_1 \|g\|_{L_\infty}$.)

\smallskip\noindent\textbf{Remark 9.6.} 1.
Unlike $\wh{t}^{\,0}$ (see (8.5)) that is controlled only in terms of
$r_0$, $\alpha_0$, $\alpha_1$, $\|g\|_{L_\infty}$, and $\|g^{-1}\|_{L_\infty}$,
 the number $\wh{t}^{\,00}$ depends on the spectral characteristics of the germ, namely, on the
  minimal distance between its different eigenvalues
$\wh{\gamma}_k(\bt)$ and $\wh{\gamma}_r(\bt)$ (where $(k,r)$ runs through $\wh{\mathcal K}$).
2. If we reject Condition~9.3 and admit intersection of the branches
$\wh{\gamma}_k(\bt)$ and $\wh{\gamma}_r(\bt)$ (for some $(k,r)\in \wh{\mathcal K}$),
then $\wh{c}^{\,\circ}_{kr}(\bt)$ will be not positive definite, and we
will be not able to choose the number $\wh{t}^{\,00}$ independently of $\bt$.

\smallskip
Under Condition 9.3, we apply Theorem 4.3 and obtain
$$
\begin{aligned}
&\| \left( e^{- i \eps^{-2} \tau \wh{\mathcal A}(\k)}- e^{- i \eps^{-2} \tau \wh{\mathcal A}^0(\k)} \right) {\mathcal R}(\k,\eps) \wh{P}\|_{L_2(\Omega) \to L_2(\Omega)}
\le ( \wh{C}'_{9} + \wh{C}_{10} |\tau|) \eps,
\cr
&\tau \in \R,\quad \eps >0,\quad |\k| \le \wh{t}^{\,00},
\end{aligned}
\eqno(9.13)
$$
where
$$
\begin{aligned}
 \wh{C}'_{9} &= \max \{ 2,  \beta_{9} r_0^{-1} \alpha_1^{1/2} \alpha_0^{-1/2} \|g\|_{L_\infty}^{1/2}\|g^{-1}\|_{L_\infty}^{1/2}  \left( 1+ n^2 \alpha_1 \|g\|_{L_\infty} (\wh{c}^{\,\circ})^{-1}\right)\},
\cr
 \wh{C}_{10} &=  \beta_{10} r_0^{-2} \alpha_1^{2} \alpha_0^{-1} \|g\|_{L_\infty}^{2}\|g^{-1}\|_{L_\infty}  \left( 1+ n^2 \alpha^2_1 \|g\|^2_{L_\infty} (\wh{c}^{\,\circ})^{-2}\right).
 \end{aligned}
\eqno(9.14)
$$
Similarly to (9.6), we have
$$
\begin{aligned}
&\| \left( e^{- i \eps^{-2} \tau \wh{\mathcal A}(\k)}- e^{- i \eps^{-2} \tau \wh{\mathcal A}^0(\k)} \right) {\mathcal R}(\k,\eps)^{1/2} \wh{P}\|_{L_2(\Omega) \to L_2(\Omega)}
\le 2 (\wh{t}^{\,00})^{-1} \eps,
\cr
&\tau \in \R,\quad \eps >0,\quad \k \in \wt{\Omega}, \quad |\k| > \wh{t}^{\,00}.
\end{aligned}
\eqno(9.15)
$$

Now, relations (9.9), (9.13), and (9.15) directly imply the following result.

\smallskip\noindent\textbf{Theorem 9.7.} \textit{Suppose that Condition} 9.3 (\textit{or more restrictive Condition} 9.4)
\textit{is satisfied. Then for any $\tau \in \R$, $\eps >0$, and $\k \in \wt{\Omega}$ we have}
$$
\| \left( e^{- i \eps^{-2} \tau \wh{\mathcal A}(\k)}- e^{- i \eps^{-2} \tau \wh{\mathcal A}^0(\k)} \right) {\mathcal R}(\k,\eps) \|_{L_2(\Omega) \to L_2(\Omega)}
\le (\wh{\mathcal C}_5 + \wh{\mathcal{C}}_6 |\tau|) \eps,
\eqno(9.16)
$$
\textit{where
$\wh{\mathcal C}_5 = \max \{ \wh{C}_{9}' , 2 (\wh{t}^{\,00})^{-1}\} + 2 r_0^{-1}$, $\wh{\mathcal{C}}_{6} = \wh{C}_{10}$,
and the constants $\wh{C}_{9}'$, $\wh{C}_{10}$, and $\wh{t}^{\,00}$ are defined by} (9.14)
\textit{and} (9.12).

\smallskip
The assumptions of Theorem 9.7 are a fortiori satisfied in the ``real'' case, if the
spectrum of the germ is simple (see Corollary~8.6 and Remark~9.5).
We arrive at the following corollary.

\smallskip\noindent\textbf{Corollary 9.8.}
\textit{Suppose that the matrices $b(\bt)$ and $g(\x)$ have real entries.
Suppose that the spectrum of the germ $\wh{S}(\bt)$ is simple for any $\bt \in {\mathbb S}^{d-1}$.
Then estimate} (9.16) \textit{holds for any $\tau \in \R$, $\eps >0$, and $\k \in \wt{\Omega}$}.

\smallskip\noindent\textbf{9.6. The sharpness of the result in the general case.}
Application of Theorem~4.4 allows us to confirm the sharpness of the result of Theorem~9.1 in the general case.

\smallskip\noindent\textbf{Theorem 9.9.}  \textit{Let $\wh{N}_0(\bt)$ be the operator defined by} (8.21).
\textit{Suppose that $\wh{N}_0(\bt_0) \ne 0$ at some point $\bt_0 \in {\mathbb S}^{d-1}$.
Let $0 \ne \tau \in \R$. Then for any $1 \le s <3$ it is impossible that the estimate}
$$
\| \left( e^{- i \eps^{-2} \tau \wh{\mathcal A}(\k)}- e^{- i \eps^{-2} \tau \wh{\mathcal A}^0(\k)} \right) {\mathcal R}(\k,\eps)^{s/2} \|_{L_2(\Omega) \to L_2(\Omega)}
\le {\mathcal C}(\tau) \eps
\eqno(9.17)
$$
\textit{holds for almost all $\k = t \bt \in \wt{\Omega}$ and sufficiently small $\eps >0$.}

\smallskip
For the proof we need the following lemma.

\smallskip\noindent\textbf{Lemma 9.10.}  \textit{Let $\wh{\delta}$ and $\wh{t}^{\,0}$ be given by} (8.4) \textit{and} (8.5), \textit{respectively.
Let $\wh{F}(\k)= \wh{F} (t,\bt)$ be the spectral projection of the operator $\wh{\A}(\k)$
for the interval $[0,\wh{\delta}]$.}
\textit{Then for $|\k|\le \wh{t}^{\,0}$ and $|\k_0| \le \wh{t}^{\,0}$ we have}
$$
\begin{aligned}
&\| \wh{F}(\k) - \wh{F}(\k_0)\|_{L_2(\Omega) \to L_2(\Omega)} \le \wh{C}' |\k-\k_0|,
\cr
&\| \wh{\A}(\k) \wh{F}(\k) -  \wh{\A}(\k_0) \wh{F}(\k_0)\|_{L_2(\Omega) \to L_2(\Omega)} \le \wh{C}'' |\k-\k_0|,
\cr
&\|e^{-i \tau \wh{\A}(\k)} \wh{F}(\k) - e^{-i \tau \wh{\A}(\k_0)} \wh{F}(\k_0)\|_{L_2(\Omega) \to L_2(\Omega)} \le (2 \wh{C}' + \wh{C}''|\tau|)|\k-\k_0|.
\end{aligned}
\eqno(9.18)
$$

\smallskip\noindent\textbf{Proof.}
First, we estimate the difference of the resolvents of the operators $\wh{\A}(\k)$ and $\wh{\A}(\k_0)$.
Consider the difference of the corresponding sesquilinear forms on the elements
$\u,\v \in \wt{H}^1(\Omega;\C^n)$:
$$
\begin{aligned}
&\wh{a}(\k)[\u,\v] - \wh{a}(\k_0)[\u,\v]
\cr
&=
\intop_\Omega \left( \langle g(\x) b(\k - \k_0) \u, b(\D+\k) \v \rangle
+ \langle g(\x) b(\D + \k_0) \u, b(\k - \k_0)  \v \rangle\right) \, d\x.
\end{aligned}
$$
Let $z\in \C$ be a common regular point of $\wh{\A}(\k)$ and $\wh{\A}(\k_0)$. Substituting
$\u = (\wh{\A}(\k) - zI)^{-1} {\boldsymbol \varphi}$ and \hbox{$\v = (\wh{\A}(\k_0) - z^* I)^{-1} {\boldsymbol \psi}$}, where
${\boldsymbol \varphi}, {\boldsymbol \psi} \in L_2(\Omega;\C^n)$, it is easy to see that
$$
\begin{aligned}
&\left| ((\wh{\A}(\k) - zI)^{-1} {\boldsymbol \varphi} - (\wh{\A}(\k_0) - z I)^{-1} {\boldsymbol \varphi}, {\boldsymbol \psi})_{L_2(\Omega)} \right|
\cr
&\le \|g\|_{L_\infty}^{1/2} \alpha_1^{1/2} | \k-\k_0 |  \| (\wh{\A}(\k) - zI)^{-1} {\boldsymbol \varphi} \|_{L_2}
\| \wh{\A}(\k_0)^{1/2} (\wh{\A}(\k_0) - z^* I)^{-1} {\boldsymbol \psi} \|_{L_2}
\cr
&+\|g\|_{L_\infty}^{1/2} \alpha_1^{1/2} | \k-\k_0 |  \| (\wh{\A}(\k_0) - z^* I)^{-1} {\boldsymbol \psi} \|_{L_2}
\| \wh{\A}(\k)^{1/2} (\wh{\A}(\k) - zI)^{-1} {\boldsymbol \varphi} \|_{L_2}
\cr
&+ 2 \|g\|_{L_\infty} \alpha_1 | \k-\k_0 |^2  \| (\wh{\A}(\k) - zI)^{-1} {\boldsymbol \varphi} \|_{L_2}
\| (\wh{\A}(\k_0) - z^* I)^{-1} {\boldsymbol \psi} \|_{L_2}.
\end{aligned}
\eqno(9.19)
$$

As follows from the results of \S 1, for $|\k|\le \wh{t}^{\,0}$ the first $n$ eigenvalues of the operator $\wh{\A}(\k)$
lie on the interval $[0,\wh{\delta} ]$, and the rest of the spectrum lies on the semiaxis $[3 \wh{\delta},\infty)$.
 Consider the contour $\Gamma_{\wh{\delta}} \subset \C$ which encloses the interval $[0, \wh{\delta}]$
 equidistantly at the distance $\wh{\delta}$.
 Then
 $$
 \wh{F}(\k) = - \frac{1}{2\pi i} \intop_{\Gamma_{\wh{\delta}}} (\wh{\A}(\k) - zI)^{-1} \,dz,
 \eqno(9.20)
 $$
where we integrate  in the positive direction.
If $z\in \Gamma_{\wh{\delta}}$, then
$$
\begin{aligned}
&\| (\wh{\A}(\k) - zI)^{-1} {\boldsymbol \varphi}\|_{L_2(\Omega)} \le {\wh{\delta}}^{-1} \|{\boldsymbol \varphi}\|_{L_2(\Omega)},
\cr
&\| (\wh{\A}(\k_0) - z^* I)^{-1} {\boldsymbol \psi}\|_{L_2(\Omega)} \le {\wh{\delta}}^{-1} \|{\boldsymbol \psi}\|_{L_2(\Omega)},
\cr
&\| \wh{\A}(\k)^{1/2}(\wh{\A}(\k) - zI)^{-1} {\boldsymbol \varphi}\|_{L_2(\Omega)} \le \sqrt{3} \,{\wh{\delta}}^{-1/2}\|{\boldsymbol \varphi}\|_{L_2(\Omega)},
\cr
&\| \wh{\A}(\k_0)^{1/2}(\wh{\A}(\k_0) - z^* I)^{-1} {\boldsymbol \psi}\|_{L_2(\Omega)} \le \sqrt{3}\,{\wh{\delta}}^{-1/2}  \|{\boldsymbol \psi}\|_{L_2(\Omega)}.
\end{aligned}
\eqno(9.21)
$$

From (9.19) and (9.21) it follows that
$$
\|  (\wh{\A}(\k) - zI)^{-1}  - (\wh{\A}(\k_0) - zI)^{-1} \|_{L_2(\Omega) \to L_2(\Omega)} \le C |\k-\k_0|
\eqno(9.22)
 $$
 for $z \in \Gamma_{\wh{\delta}}$ and $|\k|, |\k_0| \le \wh{t}^{\,0}$.
 Now representation (9.20) (at the points $\k$ and $\k_0$) and (9.22) imply the first
 inequality in (9.18).

The second estimate in (9.18) is deduced from the representation
 $$
 \wh{\A}(\k) \wh{F}(\k) = - \frac{1}{2\pi i} \intop_{\Gamma_{\wh{\delta}}} z (\wh{\A}(\k) - zI)^{-1} \,dz
 \eqno(9.23)
 $$
 at the points $\k$ and $\k_0$  with the help of (9.22).

Let us prove the third inequality. We have
$$
\begin{aligned}
&e^{-i \tau \wh{\A}(\k)} \wh{F}(\k) - e^{-i \tau \wh{\A}(\k_0)} \wh{F}(\k_0)
= e^{-i \tau \wh{\A}(\k)} \wh{F}(\k)(\wh{F}(\k) - \wh{F}(\k_0))
\cr
&+ (\wh{F}(\k) - \wh{F}(\k_0))  e^{-i \tau \wh{\A}(\k_0)} \wh{F}(\k_0)
+ \Xi(\tau,\k, \k_0),
\end{aligned}
\eqno(9.24)
$$
where
$$
\Xi(\tau,\k, \k_0)= e^{-i \tau \wh{\A}(\k)} \wh{F}(\k) \wh{F}(\k_0) - \wh{F}(\k)e^{-i \tau \wh{\A}(\k_0)} \wh{F}(\k_0).
$$
The sum of the first two terms in (9.24) does not exceed $2 \wh{C}' |\k - \k_0 |$, in view of the first estimate in (9.18).
The third term can be written as
$$
\begin{aligned}
&\Xi(\tau,\k, \k_0) = e^{-i \tau \wh{\A}(\k)} \Sigma(\tau,\k,\k_0),
\cr
&\Sigma(\tau,\k,\k_0) =
\wh{F}(\k) \wh{F}(\k_0) - e^{i \tau \wh{\A}(\k)} \wh{F}(\k)e^{-i \tau \wh{\A}(\k_0)} \wh{F}(\k_0).
\end{aligned}
$$
Obviously, $\Sigma(0,\k,\k_0)=0$, and
$$
\frac{d\Sigma(\tau,\k,\k_0)}{d\tau} =  -i \wh{F}(\k) e^{i \tau \wh{\A}(\k)} ( \wh{\A}(\k)\wh{F}(\k) - \wh{\A}(\k_0)\wh{F}(\k_0)) e^{-i \tau \wh{\A}(\k_0)} \wh{F}(\k_0).
$$
Integrating over the interval $[0,\tau]$ and using the second estimate in (9.18), we obtain
$$
\|\Xi (\tau,\k,\k_0)\| = \|\Sigma(\tau,\k,\k_0)\| \le \wh{C}'' |\tau| |\k-\k_0|.
$$
We arrive at  the third estimate in (9.18). $\bullet$

\smallskip\noindent\textbf{Proof of Theorem 9.9.} We prove by contradiction.
Let us fix $\tau \ne 0$. Assume that for some $1\le s<3$ there exists a constant ${\mathcal C}(\tau)>0$ such that
estimate (9.17) holds for almost every $\k\in \wt{\Omega}$ and sufficiently small $\eps>0$.
By (9.9) and (9.2), it follows that there exists a constant $\wt{\mathcal C}(\tau)>0$ such that
$$
\| \left( e^{- i \eps^{-2} \tau \wh{\mathcal A}(\k)}- e^{- i \eps^{-2} \tau \wh{\mathcal A}^0(\k)} \right) \wh{P}  \|_{L_2(\Omega) \to L_2(\Omega)} \eps^s (|\k|^2 + \eps^2)^{-s/2}
\le \wt{\mathcal C}(\tau) \eps
\eqno(9.25)
$$
for almost every $\k\in \wt{\Omega}$ and sufficiently small $\eps$.

Now, let $|\k| \le \wh{t}^{\,0}$. By (1.13),
$$
\| \wh{F}(\k) - \wh{P} \|_{L_2(\Omega) \to L_2(\Omega)} \le \wh{C}_1 |\k|,\quad |\k| \le \wh{t}^{\,0}.
\eqno(9.26)
$$
From (9.25) and (9.26) it follows that there exists a constant $\check{\mathcal C}(\tau)>0$ such that
$$
\|  e^{- i \eps^{-2} \tau \wh{\mathcal A}(\k)} \wh{F}(\k) - e^{- i \eps^{-2} \tau \wh{\mathcal A}^0(\k)}  \wh{P}  \|_{L_2(\Omega) \to L_2(\Omega)} \eps^s (|\k|^2 + \eps^2)^{-s/2}
\le \check{\mathcal C}(\tau) \eps
\eqno(9.27)
$$
for almost every $\k$ in the ball $|\k|\le \wh{t}^{\,0}$ and sufficiently small $\eps$.

Observe that $\wh{P}$ is the spectral projection of the operator $\wh{\A}^0(\k)$ for the interval $[0,\wh{\delta}]$.
Applying Lemma 9.10 to $\wh{\A}(\k)$ and $\wh{\A}^0(\k)$, we conclude that
for fixed $\tau$ and $\eps$ the operator under the norm sign in (9.27) is continuous with respect to $\k$  in the ball $|\k|\le \wh{t}^{\,0}$.
Consequently, estimate (9.27) holds for all $\k$ in that ball.
In particular, it holds at the point $\k=t \bt_0$ if $t \le \wh{t}^{\,0}$.
Applying (9.26) once more, we see that
$$
\| \left( e^{- i \eps^{-2} \tau \wh{\mathcal A}(t\bt_0)}  - e^{- i \eps^{-2} \tau \wh{\mathcal A}^0(t\bt_0)}  \right) \wh{P}  \|_{L_2(\Omega) \to L_2(\Omega)} \eps^s (t^2 + \eps^2)^{-s/2}
\le \check{\mathcal C}'(\tau) \eps
\eqno(9.28)
$$
for all $t \le \wh{t}^{\,0}$ and sufficiently small $\eps$.

Estimate (9.28) corresponds to the abstract estimate (4.3).
Since \hbox{$\wh{N}_0(\bt_0) \ne 0$}, applying Theorem~4.4, we arrive at a contradiction.  $\bullet$

\section*{§10. The operator ${\mathcal A}(\k)$. Application of the scheme of \S 5}

\smallskip\noindent\textbf{10.1. Application of the scheme of \S 5 to the operator ${\mathcal A}(\k)$.}
We apply the scheme of \S 5 to study the operator ${\mathcal A}(\k) = f^* \wh{\mathcal A}(\k)f$.
Now $\H = \wh{\H}=L_2(\Omega;\C^n)$, $\H_* = L_2(\Omega;\C^m)$, the role of $A(t)$ is played by $A(t,\bt)= {\mathcal A}(\k)$,
the role of  $\wh{A}(t)$ is played by $\wh{A}(t,\bt) = \wh{\mathcal A}(\k)$.
Next, the  isomorphism $M$ is the operator of multiplication by the matrix-valued function $f(\x)$.
The operator $Q$ is the operator of multiplication by the matrix-valued function
$$
Q(\x) = (f(\x) f(\x)^*)^{-1}.
$$
The block of the operator $Q$ in the subspace $\wh{\NN}$ (see (8.2)) is the operator of multiplication by the constant matrix
$$
\overline{Q} = (\underline{f f^*} )^{-1} = |\Omega|^{-1} \intop_{\Omega}( f(\x) f(\x)^*)^{-1} \,d\x.
$$
Next, $M_0$ is the operator of multiplication by the constant matrix
$$
f_0 = (\overline{Q})^{-1/2}= (\underline{f f^*} )^{1/2}.
\eqno(10.1)
$$
Note that
$$
|f_0| \le \|f \|_{L_\infty},\quad |f_0^{-1}| \le \|f^{-1} \|_{L_\infty}.
\eqno(10.2)
$$

In $L_2(\R^d;\C^n)$, we define the operator
$$
{\mathcal A}^0:= f_0 \wh{\mathcal A}^0 f_0 = f_0 b(\D)^* g^0 b(\D) f_0.
\eqno(10.3)
$$
Let ${\mathcal A}^0(\k)$ be the corresponding family of operators in $L_2(\Omega; \C^n)$. Then ${\mathcal A}^0(\k) = f_0 \wh{\mathcal A}^0(\k) f_0$.
By (8.3) and (8.10), we have
$$
f_0 \wh{S}(\k) f_0 \wh{P} = {\mathcal A}^0(\k) \wh{P}.
\eqno(10.4)
$$

\smallskip\noindent\textbf{10.2. The analytic branches of eigenvalues and eigenvectors.}
According to (5.3), the spectral germ $S(\bt)$ of the operator $A(t,\bt)$ acting in the subspace $\NN$ (see (7.3)) is represented as
$$
S(\bt) = P f^* b(\bt)^* g^0 b(\bt) f\vert_{\NN},
$$
where $P$ is the orthogonal projection of $L_2(\Omega;\C^n)$ onto $\NN$.

The analytic (in $t$) branches of the eigenvalues ${\lambda}_l(t,\bt)$ and the branches of the eigenvectors ${\varphi}_l(t,\bt)$
of ${A}(t,\bt)$ admit the power series expansions of the form (1.6), (1.7) with the coefficients depending on $\bt$:
$$
{\lambda}_l(t,\bt) = {\gamma}_l(\bt) t^2 + {\mu}_l(\bt) t^3 + \dots,\quad l=1,\dots,n,
\eqno(10.5)
$$
$$
{\varphi}_l(t,\bt) = {\omega}_l(\bt) t^2 + t {\psi}_l^{(1)}(\bt) + \dots,\quad l=1,\dots,n.
\eqno(10.6)
$$
The vectors  $\omega_1(\bt),\dots, \omega_n(\bt)$ form an orthonormal basis in the subspace $\NN$ (see (7.3)),
and the vectors
$$
\zeta_l(\bt) := f \omega_l(\bt),\quad l=1,\dots,n,
$$
form a basis in  $\wh{\NN}$ (see (8.2)) orthonormal with the weight $\overline{Q}$, i.~e., $(\overline{Q} \zeta_l(\bt), \zeta_j(\bt)) = \delta_{jl}$, $j,l=1,\dots,n$.

The numbers $\gamma_l(\bt)$ and the elements $\omega_l(\bt)$ are eigenvalues and eigenvectors of the spectral germ $S(\bt)$.
However, it is more convenient to turn to the generalized spectral problem for $\wh{S}(\bt)$.
According to (5.7), the numbers  ${\gamma}_l(\bt)$ and the elements ${\zeta}_l(\bt)$ are eigenvalues and eigenvectors
of the following generalized spectral problem:
$$
b(\bt)^* g^0 b(\bt) {\zeta}_l(\bt) = {\gamma}_l(\bt) \overline{Q} {\zeta}_l(\bt),\quad l=1,\dots,n.
\eqno(10.7)
$$

\subsection*{10.3. The operator $\wh{N}_Q(\bt)$}
We need to describe the operator $\wh{N}_Q(\bt)$ (in abstract terms it was defined in Subsection~5.2).
Let $\Lambda_Q(\x)$ be the $\Gamma$-periodic solution of the problem
$$
b(\D)^* g(\x)(b(\D) \Lambda_Q(\x) + \1_m)=0,\quad \intop_\Omega Q(\x) \Lambda_Q(\x)\, d\x =0.
\eqno(10.8)
$$
Clearly, $\Lambda_Q(\x)$ differs from the periodic solution $\Lambda(\x)$ of the problem (8.7) by a constant summand:
$$
\Lambda_Q(\x) = \Lambda(\x) + \Lambda_Q^0,\quad \Lambda_Q^0 = - (\overline{Q})^{-1} (\overline{Q \Lambda}).
\eqno(10.9)
$$

As shown in  [BSu3,\S 5],  the operator $\wh{N}_Q(\bt)$ takes the form
$$
\wh{N}_Q(\bt) = b(\bt)^* L_Q(\bt) b(\bt) \wh{P},
\eqno(10.10)
$$
where $L_Q(\bt)$ is an $(m\times m)$-matrix given by
$$
L_Q(\bt) = |\Omega|^{-1} \intop_\Omega \left( \Lambda_Q(\x)^* b(\bt)^* \wt{g}(\x) + \wt{g}(\x)^* b(\bt) \Lambda_Q(\x)\right) \,d\x.
\eqno(10.11)
$$
Combining (10.9), (10.11), and (8.19), we see that
$$
L_Q(\bt) = L(\bt) + L_Q^0(\bt), \quad  L_Q^0(\bt)= (\Lambda_Q^0)^* b(\bt)^* g^0 + g^0 b(\bt) \Lambda_Q^0.
$$

Observe that $L_Q(\k):=t L_Q(\bt)$, $\k\in \R^d$, is  a Hermitian matrix-valued function first order homogeneous in $\k$.
We put $\wh{N}_Q(\k):=t^3 \wh{N}_Q(\bt)$, $\k\in \R^d$. Then
$\wh{N}_Q(\k) = b(\k)^* L_Q(\k) b(\k) \wh{P}$. The matrix-valued function $b(\k)^* L_Q(\k)b(\k)$ is a homogeneous polynomial of the third degree  in $\k \in \R^d$.
It follows that either $\wh{N}_Q(\bt)=0$ identically for $\bt \in {\mathbb S}^{d-1}$, or $\wh{N}_Q(\bt)\ne 0$ for most points $\bt$
(except for the zeroes of this polynomial).

Some cases where the operator (10.10) is equal to zero were distinguished in [BSu3, \S 5].

\smallskip\noindent\textbf{Proposition 10.1.}
\textit{Suppose that at least one of the following conditions is satisfied}:

\noindent
$1^\circ$. \textit{The operator ${\mathcal A}$ has the form
 ${\mathcal A}= f(\x)^*\D^* g(\x)\D f(\x)$, where $g(\x)$ is a symmetric matrix with real entries.}

\noindent
$2^\circ$. \textit{Relations} (8.14) \textit{are satisfied, i.~e.} $g^0 = \overline{g}$.

\noindent\textit{Then $\wh{N}_Q(\bt)=0$ for all} $\bt \in {\mathbb S}^{d-1}$.

\smallskip
Recall that (see Subsections 5.2, 5.3) $\wh{N}_Q(\bt)= \wh{N}_{0,Q}(\bt) + \wh{N}_{*,Q}(\bt)$. By (5.9),
$$
\wh{N}_{0,Q}(\bt) = \sum_{l=1}^n \mu_l(\bt) (\cdot, \overline{Q}\zeta_l(\bt))_{L_2(\Omega)} \overline{Q}\zeta_l(\bt).
$$
We have
$$
(\wh{N}_Q(\bt) {\zeta}_l(\bt), {\zeta}_l(\bt))_{L_2(\Omega)} =(\wh{N}_{0,Q}(\bt) {\zeta}_l(\bt), {\zeta}_l(\bt))_{L_2(\Omega)} = {\mu}_l(\bt),
\quad l=1,\dots,n.
\eqno(10.12)
$$

Now we assume that the matrices $b(\bt)$, $g(\x)$, and $Q(\x)$ \textit{have real entries}. Then the matrix $\Lambda_Q(\x)$ (see (10.8))
has purely imaginary entries, and $\wt{g}(\x)$ and $g^0$ have real entries. In this case $L_Q(\bt)$ (see (10.11)) and $b(\bt)^* L_Q(\bt) b(\bt)$ are Hermitian matrices with
purely imaginary entries. If the analytic branches of the eigenvalues ${\lambda}_l(t,\bt)$ and the analytic branches of the eigenvectors ${\varphi}_l(t,\bt)$
of the operator ${A}(t,\bt)$ can be chosen so that the vectors $\zeta_l(\bt) = f {\omega}_l(\bt)$, $l=1,\dots,n,$ are real, then, by (10.12), we have 
 ${\mu}_l(\bt)=0$, $l=1,\dots,n$, i.~e., $\wh{N}_{0,Q}(\bt)=0$. We arrive at the following statement.

\smallskip\noindent\textbf{Proposition 10.2.} \textit{Suppose that the matrices $b(\bt)$, $g(\x)$, and $Q(\x)$ have real entries. Suppose that in the expansions} (10.6)
\textit{for the analytic branches of the eigenvectors of ${A}(t,\bt)$ the ``embrios''  ${\omega}_l(\bt)$, $l=1,\dots,n$, can be chosen so that the vectors $\zeta_l(\bt)= f \omega_l(\bt)$ are real.}
\textit{Then in} (10.5) \textit{we have ${\mu}_l(\bt)=0$, $l=1,\dots,n$, i.~e., $\wh{N}_{0,Q}(\bt)=0$ for all $\bt \in {\mathbb S}^{d-1}$.}

\smallskip
 In the ``real'' case under consideration, the operator $\wh{S}(\bt)$ is a symmetric matrix with real entries; $\overline{Q}$
is also a symmetric matrix with real entries.
 Clearly,  if the eigenvalue ${\gamma}_j(\bt)$ of the generalized problem (10.7) is simple, then the eigenvector $\zeta_j(\bt) = f {\omega}_j(\bt)$
 is defined uniquely up to a phase factor, and we always can choose it to be real.  We arrive at the following corollary.

\smallskip\noindent\textbf{Corollary 10.3.} \textit{Suppose that the matrices $b(\bt)$, $g(\x)$, and $Q(\x)$ have real entries.
Suppose that the spectrum of the generalized spectral problem} (10.7) \textit{is simple. Then $\wh{N}_{0,Q}(\bt)=0$ for all $\bt \in {\mathbb S}^{d-1}$.}

\smallskip\noindent\textbf{10.4. Multiplicities of the eigenvalues of the germ.}
This subsection concerns the case where $n\ge 2$.
We return to the notation of \S 2, tracing the multiplicities of the eigenvalues of the spectral germ  ${S}(\bt)$.
From what was said in Subsection 10.2 it follows that these eigenvalues are also the eigenvalues of the generalized problem (10.7).
In general, the number $p(\bt)$ of different eigenvalues $\gamma_1^\circ(\bt),\dots, \gamma^\circ_{p(\bt)}(\bt)$ of this problem and their multiplicities
$k_1(\bt),\dots, k_{p(\bt)}(\bt)$ depend on the parameter $\bt \in {\mathbb S}^{d-1}$.
 For a fixed $\bt$, let $\NN_j(\bt)$ be the eigenspace  of the germ $S(\bt)$ corresponding to the eigenvalue $\gamma_j^\circ(\bt)$.
Then $f \NN_j(\bt)$ is the eigenspace of the problem (10.7) corresponding to the same eigenvalue  $\gamma_j^\circ(\bt)$.
Let ${\mathcal P}_j(\bt)$ denote the ``skew'' projection of $L_2(\Omega;\C^n)$ onto the subspace $f \NN_j(\bt)$;
${\mathcal P}_j(\bt)$ is orthogonal with respect to the inner product with the weight $\overline{Q}$.
Then, by (5.11), we have the following invariant representations for the operators $\wh{N}_{0,Q}(\bt)$ and $\wh{N}_{*,Q}(\bt)$:
$$
\wh{N}_{0,Q}(\bt) = \sum_{j=1}^{p(\bt)} {\mathcal P}_j(\bt)^* \wh{N}_{Q}(\bt) {\mathcal P}_j (\bt),
\eqno(10.13)
$$
$$
\wh{N}_{*,Q}(\bt) = \sum_{1\le j,l \le p(\bt):\, j\ne l} {\mathcal P}_j(\bt)^* \wh{N}_{Q}(\bt) {\mathcal P}_l (\bt).
$$

\section*{§11. Approximation of the smoothed sandwiched \\ operator $e^{-i \eps^{-2} \tau {\mathcal A}(\k)}$}

\smallskip\noindent\textbf{11.1. Approximation of the smoothed sandwiched operator $e^{-i \eps^{-2} \tau {\mathcal A}(\k)}$ in the general case.}
We apply theorems of \S 5 to the operator ${\mathcal A}(\k)$. First we apply Theorem 5.7.
The constant $t^0$ is given by (7.12). Taking (7.10) and (7.11) into account, instead of the precise values of the constants $C_1$ and $C_2$ which now depend on $\bt$
we take the larger values
$$
\begin{aligned}
C_1 &= 2 \beta_1 r_0^{-1} \alpha_1^{1/2} \alpha_0^{-1/2} \|g\|_{L_\infty}^{1/2}\|g^{-1}\|_{L_\infty}^{1/2} \|f\|_{L_\infty} \|f^{-1}\|_{L_\infty},
\cr
C_2 &= 2 \beta_2 r_0^{-1} \alpha_1^{3/2} \alpha_0^{-1/2} \|g\|_{L_\infty}^{3/2}\|g^{-1}\|_{L_\infty}^{1/2} \|f\|^3_{L_\infty} \|f^{-1}\|_{L_\infty}.
 \end{aligned}
 \eqno(11.1)
$$
Applying (5.16) for the operator ${\mathcal A}(\k)$ and using (9.2) and (10.4), we obtain
$$
\begin{aligned}
&\| \left( f e^{- i \eps^{-2} \tau {\mathcal A}(\k)} f^{-1} - f_0 e^{- i \eps^{-2} \tau {\mathcal A}^0(\k)} f_0^{-1} \right) {\mathcal R}(\k,\eps)^{3/2} \wh{P}\|_{L_2(\Omega) \to L_2(\Omega)}
\cr
&\le \|f\|^2_{L_\infty} \|f^{-1}\|^2_{L_\infty} ({C}_1 + {C}_2 |\tau|) \eps,
\quad
\tau \in \R,\quad \eps >0,\quad |\k| \le {t}^0.
\end{aligned}
\eqno(11.2)
$$

 For $|\k|> t^0$ estimates are trivial. Using the analog of (9.5) for $|\k|>t^0$ and (10.2), we have
$$
\begin{aligned}
&\| \left( f e^{- i \eps^{-2} \tau {\mathcal A}(\k)} f^{-1} - f_0 e^{- i \eps^{-2} \tau {\mathcal A}^0(\k)} f_0^{-1} \right) {\mathcal R}(\k,\eps)^{1/2} \wh{P}\|_{L_2(\Omega) \to L_2(\Omega)}
\cr
&\le 2 \|f\|_{L_\infty} \|f^{-1}\|_{L_\infty}  (t^0 )^{-1} \eps,
\quad
\tau \in \R,\quad \eps >0,\quad |\k| > {t}^0.
\end{aligned}
\eqno(11.3)
$$

Next, by (9.8),
$$
\begin{aligned}
&\| \left( f e^{- i \eps^{-2} \tau {\mathcal A}(\k)} f^{-1} - f_0 e^{- i \eps^{-2} \tau {\mathcal A}^0(\k)} f_0^{-1} \right) {\mathcal R}(\k,\eps)^{1/2} (I- \wh{P}) \|_{L_2(\Omega) \to L_2(\Omega)}
\cr
&\le 2 \|f\|_{L_\infty} \|f^{-1}\|_{L_\infty}  r_0^{-1}  \eps,
\quad
\tau \in \R,\quad \eps >0,\quad \k \in \wt{\Omega}.
\end{aligned}
\eqno(11.4)
$$

Finally, relations (11.2)--(11.4) (and expressions for the constants) imply the following result.

\smallskip\noindent\textbf{Theorem 11.1.} \textit{For $\tau \in \R$, $\eps >0$, and $\k \in \wt{\Omega}$ we have}
$$
\begin{aligned}
&\| \left( f e^{- i \eps^{-2} \tau {\mathcal A}(\k)} f^{-1} - f_0 e^{- i \eps^{-2} \tau {\mathcal A}^0(\k)} f_0^{-1} \right) {\mathcal R}(\k,\eps)^{3/2} \|_{L_2(\Omega) \to L_2(\Omega)}
\cr
&\le  ({\mathcal C}_1 + {\mathcal C}_2 |\tau|) \eps,
\end{aligned}
$$
\textit{where}
$$
\begin{aligned}
{\mathcal C}_1 &=
\wt{\beta}_1 r_0^{-1} \alpha_1^{1/2} \alpha_0^{-1/2} \|g\|^{1/2}_{L_\infty} \|g^{-1}\|^{1/2}_{L_\infty} \|f\|^3_{L_\infty} \|f^{-1}\|^3_{L_\infty},
\cr
 {\mathcal C}_2 &=
 2 \beta_2 r_0^{-1} \alpha_1^{3/2} \alpha_0^{-1/2} \|g\|_{L_\infty}^{3/2}\|g^{-1}\|_{L_\infty}^{1/2} \|f\|^5_{L_\infty} \|f^{-1}\|^3_{L_\infty}.
 \end{aligned}
\eqno(11.5)
$$

\subsection*{11.2. Refinement of approximation of the smoothed sandwiched operator $e^{-i \eps^{-2} \tau {\mathcal A}(\k)}$
in the case where $\wh{N}_Q(\bt) =0$}
Now we assume that $\wh{N}_Q(\bt)=0$ and apply Theorem~5.8. By (9.2) and (10.4), we have
$$
\begin{aligned}
&\| \left( f e^{- i \eps^{-2} \tau {\mathcal A}(\k)} f^{-1} - f_0 e^{- i \eps^{-2} \tau {\mathcal A}^0(\k)} f_0^{-1} \right) {\mathcal R}(\k,\eps) \wh{P}\|_{L_2(\Omega) \to L_2(\Omega)}
\cr
&\le  \|f\|^2_{L_\infty} \|f^{-1}\|^2_{L_\infty} ({ C}_1' + { C}_5 |\tau|) \eps,\quad \tau \in \R, \quad \eps>0,\quad |\k|\le t^0.
\end{aligned}
$$
Here $C_1'= \max\{2, C_1\}$, the constant $C_1$ is defined by (11.1), and the constant $C_5$ is given by
${C}_5 = 4 \beta_5 r_0^{-2} \alpha_1^2 \alpha_0^{-1} \|g\|^2_{L_\infty} \| g^{-1}\|_{\infty}\|f\|^4_{L_\infty} \| f^{-1}\|^2_{\infty}$ .

Together with (11.3) and (11.4) this implies the following result.

\smallskip\noindent\textbf{Theorem 11.2.} \textit{Let $\wh{N}_Q(\bt)$ be the operator defined by} (10.10), (10.11). \textit{Suppose that
$\wh{N}_Q(\bt)=0$ for all $\bt \in {\mathbb S}^{d-1}$.
Then for $\tau \in \R$, $\eps >0$, and $\k \in \wt{\Omega}$ we have}
$$
\| \left( f e^{- i \eps^{-2} \tau {\mathcal A}(\k)} f^{-1}-  f_0 e^{- i \eps^{-2} \tau {\mathcal A}^0(\k)} f_0^{-1} \right) {\mathcal R}(\k,\eps) \|_{L_2(\Omega) \to L_2(\Omega)}
\le ({\mathcal C}_3 + {\mathcal{C}}_4 |\tau|) \eps,
$$
\textit{where
$$
\begin{aligned}
{\mathcal C}_3 &= \max \{ 2\|f\|_{L_\infty} \|f^{-1}\|_{L_\infty}(\|f\|_{L_\infty} \|f^{-1}\|_{L_\infty} + r_0^{-1}), {\mathcal C}_1\},
\cr
\mathcal{C}_4 &= 4 \beta_5 r_0^{-2} \alpha_1^2 \alpha_0^{-1} \|g\|^2_{L_\infty} \| g^{-1}\|_{L_\infty}\|f\|^6_{L_\infty} \| f^{-1}\|^4_{L_\infty},
\end{aligned}
\eqno(11.6)
$$
 and ${\mathcal C}_1$ is defined by} (11.5).

\smallskip
Recall that some sufficient conditions ensuring that $\wh{N}_Q(\bt)=0$ for all $\bt \in {\mathbb S}^{d-1}$ are given in Proposition~10.1.

\subsection*{11.3. Refinement of approximation of the smoothed sandwiched operator $e^{-i \eps^{-2} \tau {\mathcal A}(\k)}$
in the case where $\wh{N}_{0,Q}(\bt)=0$}
Now we reject the assumption of Theorem~11.2, but instead we assume that $\wh{N}_{0,Q}(\bt)=0$ for all $\bt$.
We may also assume that $\wh{N}_Q(\bt)= \wh{N}_{*,Q}(\bt) \ne 0$ for some $\bt$, and then for most points $\bt$ (otherwise, one
can apply Theorem~11.2.) As in Subsection~9.5, in order to apply ``abstract'' Theorem~5.9, we have to
impose an additional condition. We use the initial enumeration of the eigenvalues
${\gamma}_1(\bt),\dots, {\gamma}_n(\bt)$ of the germ ${S}(\bt)$ (each eigenvalue is repeated corresponding to its multiplicity)
and enumerate them in the nondecreasing order:
$$
{\gamma}_1(\bt) \le {\gamma}_2(\bt) \le \dots \le {\gamma}_n(\bt).
\eqno(11.7)
$$
As has been already mentioned, the numbers (11.7) are also the eigenvalues of the generalized spectral problem (10.7).
 For each $\bt$, let ${\mathcal P}^{(k)}(\bt)$ be the ``skew'' projection (orthogonal with the weight $\overline{Q}$) of $L_2(\Omega;\C^n)$
 onto the eigenspace of the problem (10.7) corresponding to the eigenvalue ${\gamma}_k(\bt)$.
      Clearly, ${\mathcal P}^{(k)}(\bt)$ coincides with one of the projections ${\mathcal P}_j(\bt)$ introduced in Subsection 10.4
      (but the number $j$ may depend on $\bt$).

\smallskip\noindent\textbf{Condition 11.3.}
{$1^\circ$.} \textit{The operator $\wh{N}_{0,Q}(\bt)$ defined by} (10.13) \textit{is equal to zero}: $\wh{N}_{0,Q}(\bt)=0$ \textit{for all} $\bt \in {\mathbb S}^{d-1}$.

\noindent{$2^\circ$.} \textit{For any pair of indices $(k,r)$, $1\le k,r \le n$, $k\ne r$,
such that ${\gamma}_k(\bt_0) = {\gamma}_r(\bt_0)$ for some $\bt_0 \in {\mathbb S}^{d-1}$, we have}
$({\mathcal P}^{(k)}(\bt))^* \wh{N}_Q(\bt) {\mathcal P}^{(r)}(\bt) =0$ \textit{for any} $\bt\in {\mathbb S}^{d-1}$.

    \smallskip
    Note that, if   ${\gamma}_k(\bt_0)={\gamma}_r(\bt_0)$,
    then the projections ${\mathcal P}^{(k)}(\bt_0)$ and ${\mathcal P}^{(r)}(\bt_0)$ coincide, and the identity
    $({\mathcal P}^{(k)}(\bt_0))^* \wh{N}_Q(\bt_0) {\mathcal P}^{(r)}(\bt_0) =0$ is valid automatically due to condition $1^\circ$.
    Condition~$2^\circ$ can be reformulated as follows: it is assumed that, for the ``blocks''
    $({\mathcal P}^{(k)}(\bt))^* \wh{N}_Q(\bt) {\mathcal P}^{(r)}(\bt)$ of the operator $\wh{N}_Q(\bt)$ that are not identically zero,
    the corresponding branches of the eigenvalues ${\gamma}_k(\bt)$ and ${\gamma}_r(\bt)$ do not intersect.

    Condition 11.3 is ensured by the following more restrictive condition.

\smallskip\noindent\textbf{Condition 11.4.}
{$1^\circ$.} \textit{The operator $\wh{N}_{0,Q}(\bt)$ defined by} (10.13) \textit{is equal to zero}: $\wh{N}_{0,Q}(\bt)=0$ \textit{for all} $\bt \in {\mathbb S}^{d-1}$.

\noindent$2^\circ$. \textit{Suppose that the number $p$ of different eigenvalues of the generalized spectral problem} (10.7) \textit{does not depend on $\bt \in {\mathbb S}^{d-1}$.
Denote the different eigenvalues of this problem enumerated in the increasing order by ${\gamma}_1^\circ(\bt),\dots, {\gamma}_p^\circ(\bt)$, and assume that their
multiplicities $k_1,\dots, k_p$ do not depend on $\bt \in {\mathbb S}^{d-1}$.}

\smallskip\noindent\textbf{Remark 11.5.}
Assumption $2^\circ$ of Condition 11.4 is a fortiori valid if the spectrum of the problem (10.7) is simple for any $\bt \in {\mathbb S}^{d-1}$.

\smallskip

So, we assume that Condition 11.3 is satisfied. We have to take care only about the pairs of indices from the set
$$
{\mathcal K}  := \{(k,r): \ 1\le k,r \le n,\ k\ne r,\ ({\mathcal P}^{(k)}(\bt))^* \wh{N}_Q(\bt) {\mathcal P}^{(r)}(\bt) \not\equiv  0 \}.
$$
Denote (cf. (2.3))
$$
{c}^{\,\circ}_{k r}(\bt) := \min \{ {c}_*,  {n}^{-1} |{\gamma}_k(\bt) - {\gamma}_{r}(\bt)|\},\quad (k, r) \in {\mathcal K}.
$$
Since the operator-valued function ${S}(\bt)$ is continuous with respect to $\bt\in {\mathbb S}^{d-1}$, then
${\gamma}_j(\bt)$ are continuous functions on the sphere ${\mathbb S}^{d-1}$.
By Condition~11.3($2^\circ$), for $(k,r)\in {\mathcal K}$ we have $|{\gamma}_k (\bt) - {\gamma}_{r}(\bt)| >0$ for all $\bt$, whence
$$
{c}^{\,\circ}_{kr} := \min_{\bt \in {\mathbb S}^{d-1}} {c}_{kr}^{\,\circ}(\bt) >0, \quad (k,r) \in {\mathcal K}.
$$
We put
$$
{c}^{\,\circ} := \min_{(k,r) \in {\mathcal K}} {c}_{kr}^\circ.
\eqno(11.8)
$$
Clearly, the number (11.8) plays the role of the number (3.33); it is important that, due to Condition~11.3($2^\circ$), we  choose $c^\circ$ idependently of $\bt$.

The number $t^{00}$ subject to (3.34) also can be chosen independently of $\bt \in {\mathbb S}^{d-1}$. Taking (7.10) and (7.11) into account, we put
$$
{t}^{\,00} = (8 \beta_2)^{-1} r_0 \alpha_1^{-3/2} \alpha_0^{1/2}  \|g\|_{L_\infty}^{-3/2} \|g^{-1} \|_{L_\infty}^{-1/2} \|f \|_{L_\infty}^{- 3} \| f^{-1} \|^{-1}_{L_\infty} {c}^{\,\circ},
\eqno(11.9)
$$
where ${c}^{\,\circ}$ is given by (11.8).
(The condition ${t}^{\,00} \le {t}^{\,0}$ is valid automatically, since ${c}^{\,\circ} \le \| {S} (\bt)\| \le \alpha_1 \|g\|_{L_\infty} \|f\|^2_{L_\infty}$.)

Under Condition~11.3, we apply Theorem~5.9 and obtain
$$
\begin{aligned}
& \| \left( f e^{- i \eps^{-2} \tau {\mathcal A}(\k)} f^{-1} - f_0 e^{- i \eps^{-2} \tau {\mathcal A}^0(\k)} f_0^{-1} \right) {\mathcal R}(\k,\eps) \wh{P}\|_{L_2(\Omega) \to L_2(\Omega)}
\cr
&\le  \|f\|_{L_\infty}^2 \|f^{-1}\|_{L_\infty}^2  ( {C}'_{9} + {C}_{10} |\tau|) \eps,
\quad \tau \in \R,\quad \eps >0,\quad |\k| \le {t}^{00},
\end{aligned}
\eqno(11.10)
$$
where
$$
\begin{aligned}
 {C}'_{9} &= \max \{ 2,  \beta_{9} r_0^{-1} \alpha_1^{1/2} \alpha_0^{-1/2} \|g\|_{L_\infty}^{1/2}\|g^{-1}\|_{L_\infty}^{1/2}\|f\|_{L_\infty}\|f^{-1}\|_{L_\infty}
 \cr
 &\times \left( 1+ n^2 \alpha_1 \|g\|_{L_\infty}\|f\|^2_{L_\infty} (c^\circ)^{-1}\right)\},
\cr
 {C}_{10} &= \beta_{10} r_0^{-2} \alpha_1^{2} \alpha_0^{-1} \|g\|_{L_\infty}^{2}\|g^{-1}\|_{L_\infty} \|f\|_{L_\infty}^{4}\|f^{-1}\|_{L_\infty}^2
  \cr
  &\times \left( 1+ n^2 \alpha_1^2 \|g\|_{L_\infty}^2  \|f\|^4_{L_\infty} (c^\circ)^{-2}\right).
 \end{aligned}
\eqno(11.11)
$$
Similarly to (11.3), we have
$$
\begin{aligned}
& \| \left( f e^{- i \eps^{-2} \tau {\mathcal A}(\k)} f^{-1} - f_0 e^{- i \eps^{-2} \tau {\mathcal A}^0(\k)} f_0^{-1} \right) {\mathcal R}(\k,\eps)^{1/2} \wh{P}\|_{L_2(\Omega) \to L_2(\Omega)}
\cr
&\le 2 \|f\|_{L_\infty} \|f^{-1}\|_{L_\infty}  ({t}^{00})^{-1} \eps,
\quad  \tau \in \R,\quad \eps >0,\quad |\k| > {t}^{00}.
\end{aligned}
\eqno(11.12)
$$

Now relations  (11.4), (11.10), and (11.12) directly imply the following result.

\smallskip\noindent\textbf{Theorem 11.6.} \textit{Suppose that Condition}~11.3 (\textit{or more restrictive Condition}~11.4)
\textit{is satisfied. Then for $\tau \in \R$, $\eps >0$, and $\k \in \wt{\Omega}$ we have}
$$
\begin{aligned}
\|& \left( f e^{- i \eps^{-2} \tau {\mathcal A}(\k)} f^{-1} - f_0 e^{- i \eps^{-2} \tau {\mathcal A}^0(\k)} f_0^{-1} \right) {\mathcal R}(\k,\eps) \|_{L_2(\Omega) \to L_2(\Omega)}
\cr
&\le ({\mathcal C}_5 + {\mathcal{C}}_6 |\tau|) \eps,
\end{aligned}
\eqno(11.13)
$$
\textit{where
${\mathcal C}_5 = \|f\|_{L_\infty} \|f^{-1}\|_{L_\infty} \left(\max \{ \|f\|_{L_\infty} \|f^{-1}\|_{L_\infty} {C}_{9}' , 2 ({t}^{00})^{-1}\} + 2 r_0^{-1}\right)$,
${\mathcal{C}}_{6} = \|f\|^2_{L_\infty} \|f^{-1}\|^2_{L_\infty} {C}_{10}$, and the constants ${C}_{9}'$, ${C}_{10}$, and ${t}^{00}$ are defined by} (11.11) \textit{and} (11.9).

\smallskip
The assumptions of Theorem~11.6 are a fortiori satisfied in the  ``real'' case, if the spectrum of the problem~(10.7) is simple (see Corollary~10.3 and Remark~11.5).
We arrive at the following corollary.

\smallskip\noindent\textbf{Corollary 11.7.}
\textit{Suppose that the matrices $b(\bt)$, $g(\x)$, and $Q(\x)$ have real entries. Suppose that the spectrum of the problem} (10.7) \textit{is simple for any $\bt \in {\mathbb S}^{d-1}$.
Then estimate} (11.13) \textit{holds for $\tau \in \R$, $\eps >0$, and $\k \in \wt{\Omega}$.}

\smallskip\noindent\textbf{11.4.  The sharpness of the result in the general case.}
Application of Theorem 5.10 allows us to confirm the sharpness of the result of Theorem~11.1 in the general case.

\smallskip\noindent\textbf{Theorem 11.8.}  \textit{Let $\wh{N}_{0,Q}(\bt)$ be the operator defined by} (10.13).
\textit{Suppose that $\wh{N}_{0,Q}(\bt_0) \ne 0$ at some point $\bt_0 \in {\mathbb S}^{d-1}$.
Let $0 \ne \tau \in \R$. Then for any $1 \le s <3$ it is impossible that the estimate}
$$
\| \left( f e^{- i \eps^{-2} \tau {\mathcal A}(\k)} f^{-1} - f_0 e^{- i \eps^{-2} \tau {\mathcal A}^0(\k)} f_0^{-1} \right) {\mathcal R}(\k,\eps)^{s/2} \|_{L_2(\Omega) \to L_2(\Omega)}
\le {\mathcal C}(\tau) \eps
\eqno(11.14)
$$
\textit{holds for almost every $\k = t \bt \in \wt{\Omega}$ and sufficiently small  $\eps >0$.}

\smallskip
For the proof we need the following lemma which can be easily checked by analogy with the proof of  Lemma~9.10.

\smallskip\noindent\textbf{Lemma 11.9.}  \textit{Let ${\delta}$ and $t^0$ be given  by} (7.10) \textit{and} (7.12), \textit{respectively.
Let ${F}(\k)= {F} (t,\bt)$ be the spectral projection of the operator ${\A}(\k)$ for the interval $[0,\delta]$.}
\textit{Then for $|\k|\le {t}^0$ and $|\k_0| \le {t}^0$ we have}
$$
\begin{aligned}
&\| {F}(\k) - {F}(\k_0)\|_{L_2(\Omega) \to L_2(\Omega)} \le {C}' |\k-\k_0|,
\cr
&\| {\A}(\k) {F}(\k) -  {\A}(\k_0) {F}(\k_0)\|_{L_2(\Omega) \to L_2(\Omega)} \le {C}'' |\k-\k_0|,
\cr
&\|e^{-i \tau {\A}(\k)} {F}(\k) - e^{-i \tau {\A}(\k_0)} {F}(\k_0)\|_{L_2(\Omega) \to L_2(\Omega)} \le (2 {C}' + {C}''|\tau|)|\k-\k_0|.
\end{aligned}
$$

\smallskip\noindent\textbf{Proof of Theorem 11.8.}
We prove by contradiction. Let us fix $\tau \ne 0$. Suppose that for some $1\le s <3$ there exists a constant ${\mathcal C}(\tau)>0$
such that estimate~(11.14) holds for almost every $\k\in \wt{\Omega}$ and sufficiently small $\eps>0$.
By (11.4) and (9.2), it follows that  there exists a constant $\wt{\mathcal C}(\tau)>0$ such that
$$
\| \left( f e^{- i \eps^{-2} \tau {\mathcal A}(\k)} f^{-1} - f_0 e^{- i \eps^{-2} \tau {\mathcal A}^0(\k)} f_0^{-1} \right) \wh{P} \|
 \eps^s (|\k|^2 + \eps^2)^{-s/2}
\le \wt{\mathcal C}(\tau) \eps
\eqno(11.15)
$$
for almost every $\k\in \wt{\Omega}$ and sufficiently small $\eps$.

By (5.2), we have $f^{-1} \wh{P} = P f^* \overline{Q}$, where $P$ is the orthogonal projection of  $L_2(\Omega;\C^n)$ onto the subspace $\NN$ (see (7.3)).
Then the operator under the norm sign in (11.15) can be written as
$ f e^{- i \eps^{-2} \tau {\mathcal A}(\k)} P f^{*} \overline{Q} - f_0 e^{- i \eps^{-2} \tau {\mathcal A}^0(\k)} f_0^{-1}  \wh{P}$.

Now, let $|\k| \le {t}^0$. By (1.13),  $\| {F}(\k) - {P} \| \le {C}_1 |\k|$ for $|\k| \le {t}^0$.
Together with  (11.15)  this implies that there exists a constant $\check{\mathcal C}(\tau)$ such that
$$
\| f e^{- i \eps^{-2} \tau {\mathcal A}(\k)} F(\k) f^{*} \overline{Q} - f_0 e^{- i \eps^{-2} \tau {\mathcal A}^0(\k)} f_0^{-1}  \wh{P} \|
 \eps^s (|\k|^2 + \eps^2)^{-s/2}
\le \check{\mathcal C}(\tau) \eps
\eqno(11.16)
$$
for almost every  $\k$ in the ball $|\k|\le {t}^0$ and sufficiently small $\eps$.

Observe that  $\wh{P}$ is the spectral projection of the operator ${\A}^0(\k)$ for the interval $[0,\delta]$.
Therefore, Lemma 11.9 (applied to  ${\A}(\k)$ and ${\A}^0(\k)$) implies that for fixed $\tau$ and $\eps$ the operator under
the norm sign in (11.16) is continuous with respect to $\k$ in the ball $|\k|\le {t}^0$.
Hence, estimate (11.16) holds for all $\k$ in this ball.
In particular, it is satisfied at the point $\k=t \bt_0$ if $t \le {t}^0$.
Applying once more the inequality $\| {F}(\k) - {P} \| \le {C}_1 |\k|$ and the identity   $P f^* \overline{Q} = f^{-1}\wh{P}$,
we see that the inequality
$$
\| \left( f e^{- i \eps^{-2} \tau {\mathcal A}(t \bt_0)} f^{-1} - f_0 e^{- i \eps^{-2} \tau {\mathcal A}^0(t \bt_0)} f_0^{-1} \right) \wh{P} \|
 \eps^s (t^2 + \eps^2)^{-s/2}
\le \check{\mathcal C}'(\tau) \eps
\eqno(11.17)
$$
holds for all $t \le t^0$ and sufficiently small $\eps$.

In abstract terms, estimate (11.17) corresponds to the inequality (5.17). Since we assume that
$\wh{N}_{0,Q}(\bt_0) \ne 0$, applying Theorem 5.10, we arrive at a contradiction.  $\bullet$

\section*{§12. Approximation of the smoothed operator $e^{-i \eps^{-2} \tau {\mathcal A}}$}

\smallskip\noindent\textbf{12.1. Approximation of the smoothed operator $e^{-i \eps^{-2} \tau \wh{\mathcal A}}$.}
 In $L_2(\R^d;\C^n)$, we consider the operator  $\wh{\mathcal A}= b(\D)^* g(\x) b(\D)$ (see (8.1)). Let
$\wh{\mathcal A}^0= b(\D)^* g^0 b(\D)$ be the effective operator (see (8.11)). Recall the notation ${\mathcal H}_0 = -\Delta$ and put
$$
{\mathcal R}(\eps) := \eps^2 ({\mathcal H}_0 + \eps^2 I)^{-1}.
\eqno(12.1)
$$
Expansion (7.7) for $\wh{\mathcal A}$ yields
$$
e^{- i \tau \eps^{-2} \wh{\mathcal A}} = {\mathcal U}^{-1}\biggl(\int_{\wt{\Omega}}\oplus e^{- i \tau \eps^{-2} \wh{\mathcal A}(\k)}\,d\k \biggr) {\mathcal U}.
$$
The operator $e^{- i \tau \eps^{-2} \wh{\mathcal A}^0}$ admits a similar representation.
The operator (12.1) expands in the direct integral of the operators~(9.1):
$$
{\mathcal R}(\eps) = {\mathcal U}^{-1}\biggl(\int_{\wt{\Omega}}\oplus {\mathcal R}(\k,\eps) \,d\k \biggr) {\mathcal U}.
\eqno(12.2)
$$
It follows that the operator $\left( e^{- i \tau \eps^{-2} \wh{\mathcal A}} - e^{- i \tau \eps^{-2} \wh{\mathcal A}^0}\right) {\mathcal R}(\eps)^{s/2}$
expands in the direct integral of the operators
$\left( e^{- i \tau \eps^{-2} \wh{\mathcal A}(\k)} - e^{- i \tau \eps^{-2} \wh{\mathcal A}^0(\k)}\right) {\mathcal R}(\k,\eps)^{s/2}$.
Hence,
$$
\begin{aligned}
&\|\left( e^{- i \tau \eps^{-2} \wh{\mathcal A}} - e^{- i \tau \eps^{-2} \wh{\mathcal A}^0}\right) {\mathcal R}(\eps)^{s/2}\|_{L_2(\R^d) \to L_2(\R^d)}
\cr
&=
\textrm{ess-}\!\sup_{\k \in \wt{\Omega}}
\|\left( e^{- i \tau \eps^{-2} \wh{\mathcal A}(\k)} - e^{- i \tau \eps^{-2} \wh{\mathcal A}^0(\k)}\right) {\mathcal R}(\k,\eps)^{s/2}\|_{L_2(\Omega) \to L_2(\Omega)}.
\end{aligned}
\eqno(12.3)
$$
Therefore, Theorem 9.1 directly implies the following statement.

\smallskip\noindent\textbf{Theorem 12.1.}  \textit{Let  $\wh{\mathcal A}$ be the operator in $L_2(\R^d;\C^n)$ given by $\wh{\mathcal A}= b(\D)^* g(\x) b(\D)$, where $g(\x)$ and $b(\D)$ satisfy the assumptions of Subsection {\rm 6.1}.
Let $\wh{\mathcal A}^0= b(\D)^* g^0 b(\D)$ be the effective operator, where $g^0$ is given by {\rm (8.9)}. Let ${\mathcal R}(\eps)$ be defined by} (12.1). \textit{Then for $\tau \in \R$ and $\eps>0$ we have}
$$
\|\left( e^{- i \tau \eps^{-2} \wh{\mathcal A}} - e^{- i \tau \eps^{-2} \wh{\mathcal A}^0}\right) {\mathcal R}(\eps)^{3/2}\|_{L_2(\R^d) \to L_2(\R^d)}
\le (\wh{\mathcal C}_1 + \wh{\mathcal C}_2 |\tau|) \eps.
$$
\textit{The constants $\wh{\mathcal C}_1$ and $\wh{\mathcal C}_2$ are defined by} (9.10)
\textit{and depend only on} $r_0$, $\alpha_0$, $\alpha_1$, $\|g\|_{L_\infty}$, \textit{and} $\| g^{-1}\|_{L_\infty}$.

\smallskip
Similarly,  Theorem~9.2 implies the following result.

\smallskip\noindent\textbf{Theorem 12.2.}  \textit{Suppose that the assumptions of Theorem} 12.1
\textit{are satisfied. Let $\wh{N}(\bt)$ be the operator defined by} (8.18), (8.19).  \textit{Suppose that}
$\wh{N}(\bt)=0$ \textit{for all} $\bt \in {\mathbb S}^{d-1}$. \textit{Then for $\tau \in \R$ and $\eps>0$ we have}
$$
\|\left( e^{- i \tau \eps^{-2} \wh{\mathcal A}} - e^{- i \tau \eps^{-2} \wh{\mathcal A}^0}\right) {\mathcal R}(\eps) \|_{L_2(\R^d) \to L_2(\R^d)}
\le (\wh{\mathcal C}_3 + \wh{\mathcal C}_4 |\tau|) \eps.
$$
\textit{The constants  $\wh{\mathcal C}_3$ and $\wh{\mathcal C}_4$ are defined in Theorem}~9.2
\textit{and depend only on} $r_0$, $\alpha_0$, $\alpha_1$, $\|g\|_{L_\infty}$, \textit{and} $\| g^{-1}\|_{L_\infty}$.

\smallskip
Recall that some sufficient conditions ensuring  that the assumptions of Theorem 12.2 are satisfied are given in Proposition~8.4.

Finally, applying  Theorem 9.7 and using the direct integral expansion, we obtain the following statement.

\smallskip\noindent\textbf{Theorem 12.3.}  \textit{Suppose that the assumptions of Theorem} 12.1
\textit{are satisfied. Suppose also that Condition}~9.3 (\textit{or more restrictive Condition}~9.4)
 \textit{is satisfied. Then for $\tau \in \R$ and $\eps>0$ we have}
$$
\|\left( e^{- i \tau \eps^{-2} \wh{\mathcal A}} - e^{- i \tau \eps^{-2} \wh{\mathcal A}^0}\right) {\mathcal R}(\eps) \|_{L_2(\R^d) \to L_2(\R^d)}
\le (\wh{\mathcal C}_5 + \wh{\mathcal C}_6 |\tau|) \eps.
$$
\textit{The constants $\wh{\mathcal C}_5$ and $\wh{\mathcal C}_6$ are defined in Theorem}~9.7
\textit{and depend only on} $r_0$, $\alpha_0$, $\alpha_1$, $\|g\|_{L_\infty}$, $\| g^{-1}\|_{L_\infty}$, \textit{and also on the number $\wh{c}^{\,\circ}$ defined by} (9.11).

\smallskip
Recall that some sufficient conditions ensuring
that the assumptions of Theorem~12.3 are satisfied are given in Corollary 9.8.

Applying Theorem 9.9, we confirm the sharpness of the result of Theorem~12.1.

\smallskip\noindent\textbf{Theorem 12.4.}  \textit{Suppose that the assumptions of Theorem} 12.1
\textit{are satisfied. Let $\wh{N}_0(\bt)$ be the operator defined by} (8.21).  \textit{Suppose that}
$\wh{N}_0(\bt_0)\ne 0$ \textit{for some $\bt_0 \in {\mathbb S}^{d-1}$. Let $0\ne \tau \in \R$.
Then for any $1\le s<3$ it is impossible that the estimate}
$$
\|\left( e^{- i \tau \eps^{-2} \wh{\mathcal A}} - e^{- i \tau \eps^{-2} \wh{\mathcal A}^0}\right) {\mathcal R}(\eps)^{s/2} \|_{L_2(\R^d) \to L_2(\R^d)}
\le {\mathcal C}(\tau)\eps
\eqno(12.4)
$$
\textit{holds for all sufficiently small $\eps >0$.}

\smallskip\noindent\textbf{Proof.} We prove by contradiction. Let us fix $\tau \ne 0$. Suppose that for some $1\le s<3$ there exists a constant ${\mathcal C}(\tau)>0$
such that  (12.4) holds for all sufficiently small $\eps$.
By (12.3), this means that for almost every $\k \in \wt{\Omega}$ and sufficiently small $\eps$ estimate (9.17) holds.
But this contradicts to the statement of Theorem~9.9. $\bullet$

\smallskip\noindent\textbf{12.2. Approximation of the smoothed sandwiched operator $e^{-i \eps^{-2} \tau {\mathcal A}}$.}
 In $L_2(\R^d;\C^n)$, consider the operator ${\mathcal A}= f(\x)^* \wh{\mathcal A}f(\x)= f(\x)^* b(\D)^* g(\x)b(\D) f(\x)$  (see (6.4)).
Let $\wh{\mathcal A}^0$ be the operator~(8.11), and let  $f_0$ be the matrix~(10.1). Let ${\mathcal A}^0=f_0 \wh{\mathcal A}^0 f_0=
f_0 b(\D)^* g^0 b(\D) f_0$ (see (10.3)).

Similarly to (12.3), using (7.7) and (12.2),  we have
$$
\begin{aligned}
&\|\left( f e^{- i \tau \eps^{-2} {\mathcal A}} f^{-1} - f_0 e^{- i \tau \eps^{-2} {\mathcal A}^0} f_0^{-1}\right) {\mathcal R}(\eps)^{s/2}\|_{L_2(\R^d) \to L_2(\R^d)}
\cr
&=
\textrm{ess-}\!\sup_{\k \in \wt{\Omega}}
\|\left( f e^{- i \tau \eps^{-2} {\mathcal A}(\k)} f^{-1} -f_0 e^{- i \tau \eps^{-2} {\mathcal A}^0(\k)} f_0^{-1} \right) {\mathcal R}(\k,\eps)^{s/2}\|_{L_2(\Omega) \to L_2(\Omega)}.
\end{aligned}
\eqno(12.5)
$$

Theorem~11.1 together with (12.5) implies the following result.

\smallskip\noindent\textbf{Theorem 12.5.}  \textit{Let
${\mathcal A}$ be the operator in $L_2(\R^d;\C^n)$ given by ${\mathcal A}= f(\x)^* b(\D)^* g(\x) b(\D) f(\x)$, where $g(\x)$, $f(\x)$
and $b(\D)$ satisfy the assumptions of Subsection {\rm 6.1}. Let
${\mathcal A}^0= f_0 b(\D)^* g^0 b(\D) f_0$, where  $g^0$ is the effective matrix} (8.9) \textit{and $f_0 = (\underline{ff^*})^{1/2}$.
Let ${\mathcal R}(\eps)$ be defined by} (12.1).
\textit{Then for $\tau \in \R$ and $\eps>0$ we have}
$$
\|\left( f e^{- i \tau \eps^{-2} {\mathcal A}} f^{-1} -f_0 e^{- i \tau \eps^{-2} {\mathcal A}^0} f_0^{-1}\right) {\mathcal R}(\eps)^{3/2}\|_{L_2(\R^d) \to L_2(\R^d)}
\le ({\mathcal C}_1 + {\mathcal C}_2 |\tau|) \eps.
$$
\textit{The constants ${\mathcal C}_1$ and ${\mathcal C}_2$ are defined by} (11.5)
\textit{and depend only on} $r_0$, $\alpha_0$, $\alpha_1$, $\|g\|_{L_\infty}$, $\| g^{-1}\|_{L_\infty}$, $\|f\|_{L_\infty}$, \textit{and} $\| f^{-1}\|_{L_\infty}$.

\smallskip
Similarly, Theorem 11.2 leads to the following statement.

\smallskip\noindent\textbf{Theorem 12.6.}  \textit{Suppose that the assumptions of Theorem} 12.5 \textit{are satisfied. Let $\wh{N}_Q(\bt)$ be the operator defined by} (10.10), (10.11).
\textit{Suppose that} \hbox{$\wh{N}_Q(\bt)=0$} \textit{for all} $\bt \in {\mathbb S}^{d-1}$. \textit{Then for $\tau \in \R$ and $\eps>0$ we have}
$$
\|\left( f e^{- i \tau \eps^{-2} {\mathcal A}} f^{-1} -f_0 e^{- i \tau \eps^{-2} {\mathcal A}^0} f_0^{-1}\right) {\mathcal R}(\eps) \|_{L_2(\R^d) \to L_2(\R^d)}
\le ({\mathcal C}_3 + {\mathcal C}_4 |\tau|) \eps.
$$
\textit{The constants ${\mathcal C}_3$ and ${\mathcal C}_4$ are given by} (11.6)
\textit{and depend only on} $r_0$, $\alpha_0$, $\alpha_1$, $\|g\|_{L_\infty}$, $\| g^{-1}\|_{L_\infty}$, $\|f\|_{L_\infty}$, \textit{and} $\| f^{-1}\|_{L_\infty}$.

\smallskip
Recall that some sufficient conditions ensuring that the assumptions of Theorem~12.6 are satisfiied are given in Proposition~10.1.

Finally, from Theorem 11.6 and the direct integral expansion  we deduce the following result.

\smallskip\noindent\textbf{Theorem 12.7.}  \textit{Suppose that the assumptions of Theorem}~12.5 \textit{are satisfied.
 Suppose also that Condition}~11.3 (\textit{or more restrictive Condition}~11.4)
\textit{is satisfied. Then for $\tau \in \R$ and  $\eps>0$ we have}
$$
\|\left( f e^{- i \tau \eps^{-2} {\mathcal A}} f^{-1} - f_0 e^{- i \tau \eps^{-2} {\mathcal A}^0} f_0^{-1} \right) {\mathcal R}(\eps) \|_{L_2(\R^d) \to L_2(\R^d)}
\le ({\mathcal C}_5 + {\mathcal C}_6 |\tau|) \eps.
$$
\textit{The constants ${\mathcal C}_5$ and ${\mathcal C}_6$ are defined in  Theorem}~11.6
\textit{and depend only on} $r_0$, $\alpha_0$, $\alpha_1$, $\|g\|_{L_\infty}$, $\| g^{-1}\|_{L_\infty}$, $\|f\|_{L_\infty}$, $\| f^{-1}\|_{L_\infty}$,
\textit{and also on the number ${c}^\circ$ defined by} (11.8).

\smallskip
Recall that some sufficient conditions ensuring that the assumptions of Theorem~12.7 are satisfied are given in Corollary~11.7.

By analogy with the proof of Theorem~12.4,  we deduce the following result from Theorem~11.8; this confirms that the result of Theorem~12.5 is sharp.

\smallskip\noindent\textbf{Theorem 12.8.}  \textit{Suppose that the assumptions of Theorem}~12.5 \textit{are satisfied. Let $\wh{N}_{0,Q}(\bt)$ be the operator defined by} (10.13).
\textit{Suppose that} $\wh{N}_{0,Q}(\bt_0)\ne 0$ \textit{for some $\bt_0 \in {\mathbb S}^{d-1}$. Let $0\ne \tau \in \R$. Then for any  $1\le s<3$ it is impossible that the estimate}
$$
\|\left(f e^{- i \tau \eps^{-2} {\mathcal A}}f^{-1} -f_0 e^{- i \tau \eps^{-2} {\mathcal A}^0} f_0^{-1} \right) {\mathcal R}(\eps)^{s/2} \|_{L_2(\R^d) \to L_2(\R^d)}
\le {\mathcal C}(\tau)\eps
\eqno(12.6)
$$
\textit{holds for all sufficiently small $\eps >0$.}

\section*{Chapter 3. Homogenization problems for  nonstationary Schr\"odinger type equation }

\section*{\S 13. Homogenization of the operator exponential $e^{-i \tau \A_\eps}$ }

\smallskip\noindent\textbf{13.1. The operators $\wh{\mathcal A}_\eps$ and ${\mathcal A}_\eps$. Statement of the problem.}
If $\phi(\x)$ is a  $\Gamma$-periodic function in~$\R^d$, we agree to denote $\phi^\eps(\x):= \phi(\eps^{-1}\x)$, $\eps>0$.
\textit{Our main objects} are  the operators $\wh{\mathcal A}_\eps$ and ${\mathcal A}_\eps$ acting in $L_2(\R^d;\C^n)$ and formally given by
$$
\wh{\mathcal A}_\eps = b(\D)^* g^\eps(\x)b(\D),
\quad
{\mathcal A}_\eps =(f^\eps(\x))^* b(\D)^* g^\eps(\x)b(\D) f^\eps(\x).
\eqno(13.1)
$$
The precise definitions are given in terms of the corresponding quadratic forms (cf. Subsection~6.1).
The coefficients of the operators~(13.1) oscillate rapidly as $\eps \to 0$.

\textit{Our goal} is to find approximations for small $\eps$ for the operator exponential $e^{-i \tau \wh{\mathcal A}_\eps}$
and for the sandwiched exponential $f^\eps e^{-i \tau {\mathcal A}_\eps } (f^\eps)^{-1}$
and to apply the results to homogenization of the Cauchy problem for the Schr\"odinger type equation.

\smallskip\noindent\textbf{13.2. The scaling transformation.}
Let $T_\eps$ be the unitary scaling transformation in $L_2(\R^d;\C^n)$ defined by
$(T_\eps \u)(\x) = \eps^{d/2}\u(\eps \x)$, $\eps>0$. Then $\A_\eps = \eps^{-2} T_\eps^* \A T_\eps$. Hence,
$$
e^{-i \tau \A_\eps} = T^*_\eps e^{- i \tau \eps^{-2} \A} T_\eps.
\eqno(13.2)
$$
The operator $\wh{\mathcal A}_\eps$ satisfies similar relations.

Applying the scaling transformation to the resolvent of the operator ${\mathcal H}_0 = -\Delta$, we obtain
$$
({\mathcal H}_0 +I)^{-1} = \eps^2 T^*_\eps ({\mathcal H}_0 + \eps^2 I)^{-1} T_\eps = T^*_\eps {\mathcal R}(\eps) T_\eps.
\eqno(13.3)
$$
Here we have used the notation (12.1).

Finally, if $\phi(\x)$ is a  $\Gamma$-periodic function, then, under the scaling transformation, the operator $[\phi^\eps]$ of multiplication by the function $\phi^\eps(\x)$
turns into the operator $[\phi]$ of multiplication by $\phi(\x)$:
$$
[\phi^\eps] = T_\eps^* [\phi] T_\eps.
\eqno(13.4)
$$

\smallskip\noindent\textbf{13.3. Approximation of the operator $e^{-i \tau \wh{\A}_\eps}$.}
We start with the simpler operator $\wh{\A}_\eps$. Let $\wh{\A}^0$ be the effective operator (8.11).
Using relations of the form (13.2) (for the operators  $\wh{\A}_\eps$ and $\wh{\A}^0$) and (13.3), we obtain
$$
\begin{aligned}
&\left( e^{-i \tau \wh{\A}_\eps} - e^{-i \tau \wh{\A}^0 }\right) ({\mathcal H}_0+I)^{-s/2}
\\
&= T_\eps^*
\left( e^{-i \tau \eps^{-2}\wh{\A}} - e^{-i \tau \eps^{-2}\wh{\A}^0 }\right) {\mathcal R}(\eps)^{s/2} T_\eps,\quad \eps>0.
\end{aligned}
\eqno(13.5)
$$

Since $T_\eps$ is unitary, combining this with Theorem~12.1, we deduce the following result (which has been proved before in [BSu5, Theorem~12.1]).

\smallskip\noindent\textbf{Theorem 13.1.}  \textit{Let $\wh{\mathcal A}_\eps= b(\D)^* g^\eps b(\D)$,
where $g(\x)$ and $b(\D)$ satisfy the assumptions of Subsection~{\rm 6.1}.
Let $\wh{\mathcal A}^0= b(\D)^* g^0 b(\D)$ be the effective operator, where $g^0$ is given by {\rm (8.9)}. Let ${\mathcal H}_0=-\Delta$.}
\textit{Then for $\tau \in \R$ and $\eps>0$ we have}
$$
\|\left( e^{- i \tau \wh{\mathcal A}_\eps} - e^{- i \tau  \wh{\mathcal A}^0}\right) ({\mathcal H}_0 +I)^{-3/2}\|_{L_2(\R^d) \to L_2(\R^d)}
\le (\wh{\mathcal C}_1 + \wh{\mathcal C}_2 |\tau|) \eps.
\eqno(13.6)
$$
\textit{The constants $\wh{\mathcal C}_1$ and $\wh{\mathcal C}_2$ are given by} (9.10)
\textit{and depend only on} $r_0$, $\alpha_0$, $\alpha_1$, $\|g\|_{L_\infty}$, \textit{and} $\| g^{-1}\|_{L_\infty}$.

\smallskip

Obviously,
$$
\| e^{- i \tau \wh{\mathcal A}_\eps} - e^{- i \tau  \wh{\mathcal A}^0}\|_{L_2(\R^d) \to L_2(\R^d)}
\le 2,\quad \tau \in \R,\quad \eps>0.
\eqno(13.7)
$$
Interpolating between (13.7) and (13.6), for $0\le s \le 3$ we obtain
$$
\begin{aligned}
&\|\left( e^{- i \tau \wh{\mathcal A}_\eps} - e^{- i \tau  \wh{\mathcal A}^0}\right) ({\mathcal H}_0 +I)^{-s/2}\|_{L_2(\R^d) \to L_2(\R^d)}
\cr
&\le 2^{1-s/3} (\wh{\mathcal C}_1 + \wh{\mathcal C}_2 |\tau|)^{s/3} \eps^{s/3},\quad \tau \in \R,\quad \eps>0.
\end{aligned}
\eqno(13.8)
$$
The operator $({\mathcal H}_0+I)^{s/2}$ is an isometric isomorphism of the Sobolev space $H^s(\R^d;\C^n)$ onto $L_2(\R^d;\C^n)$.
Therefore,  (13.8) is equivalent to
$$
\| e^{- i \tau \wh{\mathcal A}_\eps} - e^{- i \tau  \wh{\mathcal A}^0} \|_{H^s(\R^d) \to L_2(\R^d)}
\le 2^{1-s/3} (\wh{\mathcal C}_1 + \wh{\mathcal C}_2 |\tau|)^{s/3} \eps^{s/3},\quad \tau \in \R,\quad \eps>0.
\eqno(13.9)
$$
In particular, for $0< \eps \le 1$ estimate (13.9) allows us  to consider large values of time $\tau$, namely, we can consider $|\tau| = O(\eps^{-\alpha})$ for \hbox{$0< \alpha <1$}.
We arrive at the following theorem (which has been proved before in  [BSu5, Theorem~12.2]).

\smallskip\noindent\textbf{Theorem 13.2.}  \textit{Suppose that the assumptions of Theorem}~13.1
\textit{are satisfied. Then for \hbox{$0\le s \le 3$},  $\tau \in \R$, and $\eps>0$ we have}
$$
\| e^{- i \tau \wh{\mathcal A}_\eps} - e^{- i \tau  \wh{\mathcal A}^0} \|_{H^s(\R^d) \to L_2(\R^d)}
\le \wh{\mathfrak C}_1(s;\tau) \eps^{s/3},
\eqno(13.10)
$$
\textit{where}
$$
\wh{\mathfrak C}_1(s;\tau)= 2^{1-s/3} (\wh{\mathcal C}_1 + \wh{\mathcal C}_2 |\tau|)^{s/3}.
\eqno(13.11)
$$
\textit{In particular, for $0< \eps \le 1$ and $|\tau|= \eps^{-\alpha}$, $0< \alpha <1$, we have}
$$
\| e^{- i \tau \wh{\mathcal A}_\eps} - e^{- i \tau  \wh{\mathcal A}^0} \|_{H^s(\R^d) \to L_2(\R^d)}
\le \wh{\mathfrak C}_1(s;1) \eps^{s(1-\alpha)/3},\quad
0< \eps \le 1, \quad |\tau|=\eps^{-\alpha}.
\eqno(13.12)
$$

\smallskip\noindent\textbf{13.4. Refinement of approximation of the operator $e^{-i \tau \wh{\A}_\eps}$ under the additional assumptions.}
Similarly, from (13.5) and Theorem~12.2 we deduce the following result.

\smallskip\noindent\textbf{Theorem 13.3.}  \textit{Suppose that the assumptions of Theorem} 13.1
\textit{are satisfied. Let $\wh{N}(\bt)$ be the operator defined by} (8.18), (8.19). \textit{Suppose that} $\wh{N}(\bt)=0$ \textit{for all} $\bt \in {\mathbb S}^{d-1}$.
\textit{Then for $\tau \in \R$ and $\eps>0$ we have}
$$
\|\left( e^{- i \tau \wh{\mathcal A}_\eps} - e^{- i \tau  \wh{\mathcal A}^0}\right) ({\mathcal H}_0 +I)^{-1}\|_{L_2(\R^d) \to L_2(\R^d)}
\le (\wh{\mathcal C}_3 + \wh{\mathcal C}_4 |\tau|) \eps.
\eqno(13.13)
$$
\textit{The constants $\wh{\mathcal C}_3$ and $\wh{\mathcal C}_4$ are defined in Theorem}~9.2
\textit{and depend only on} $r_0$, $\alpha_0$, $\alpha_1$, $\|g\|_{L_\infty}$, \textit{and} $\| g^{-1}\|_{L_\infty}$.

\smallskip

Interpolating between (13.7) and (13.13), we obtain the following statement.

\smallskip\noindent\textbf{Theorem 13.4.}  \textit{Suppose that the assumptions of Theorem}~13.3
\textit{are satisfied. Then for \hbox{$0\le s \le 2$},  $\tau \in \R$, and $\eps>0$ we have}
$$
\| e^{- i \tau \wh{\mathcal A}_\eps} - e^{- i \tau  \wh{\mathcal A}^0} \|_{H^s(\R^d) \to L_2(\R^d)}
\le \wh{\mathfrak C}_2(s;\tau) \eps^{s/2},
$$
\textit{where}
$$
\wh{\mathfrak C}_2(s;\tau)= 2^{1-s/2} (\wh{\mathcal C}_3 + \wh{\mathcal C}_4 |\tau|)^{s/2}.
\eqno(13.14)
$$
\textit{In particular, for $0< \eps \le 1$ and $|\tau|= \eps^{-\alpha}$, $0< \alpha <1$, we have}
$$
\| e^{- i \tau \wh{\mathcal A}_\eps} - e^{- i \tau  \wh{\mathcal A}^0} \|_{H^s(\R^d) \to L_2(\R^d)}
\le \wh{\mathfrak C}_2(s;1) \eps^{s(1-\alpha)/2},\quad
0< \eps \le 1, \quad |\tau|=\eps^{-\alpha}.
$$

\smallskip
Theorem 13.4 and  Proposition 8.4 imply the following statement.

\smallskip\noindent\textbf{Corollary 13.5.}
\textit{Suppose that at least one of the following conditions is fulfilled}:

\noindent
$1^\circ$. \textit{The operator $\wh{\mathcal A}_\eps$ has the form
 $\wh{\mathcal A}_\eps= \D^* g^\eps(\x)\D$, where $g(\x)$ is a symmetric matrix with real entries.}

\noindent
$2^\circ$. \textit{Relations} (8.14) \textit{are satisfied, i.~e.,} $g^0 = \overline{g}$.

\noindent
$3^\circ$. \textit{Relations} (8.15) \textit{are satisfied, i.~e.,} $g^0 = \underline{g}$. (\textit{In particular, this is true if} $m=n$.)

\noindent\textit{Then the statements of Theorem} 13.4 \textit{are valid}.

\smallskip
Finally, Theorem 12.3 and (13.5) lead to the following result.

\smallskip\noindent\textbf{Theorem 13.6.}  \textit{Suppose that the assumptions of Theorem}~13.1
\textit{are satisfied. Suppose that Condition}~9.3 (\textit{or more restrictive Condition}~9.4)
\textit{is satisfied. Then for $\tau \in \R$ and $\eps>0$ we have}
$$
\|\left( e^{- i \tau \wh{\mathcal A}_\eps} - e^{- i \tau  \wh{\mathcal A}^0}\right) ({\mathcal H}_0 +I)^{-1}\|_{L_2(\R^d) \to L_2(\R^d)}
\le (\wh{\mathcal C}_5 + \wh{\mathcal C}_6 |\tau|) \eps.
\eqno(13.15)
$$
\textit{The constants $\wh{\mathcal C}_5$ and $\wh{\mathcal C}_6$ are defined in Theorem}~9.7
\textit{and depend only on} $r_0$, $\alpha_0$, $\alpha_1$, $\|g\|_{L_\infty}$, $\| g^{-1}\|_{L_\infty}$, \textit{and on the number $\wh{c}^\circ$ given by}~(9.11).

\smallskip

Interpolating between (13.7) and (13.15), we obtain the following theorem.

\smallskip\noindent\textbf{Theorem 13.7.}  \textit{Suppose that the assumptions of Theorem} 13.6
\textit{are satisfied. Then for \hbox{$0\le s \le 2$},  $\tau \in \R$, and $\eps>0$ we have}
$$
\| e^{- i \tau \wh{\mathcal A}_\eps} - e^{- i \tau  \wh{\mathcal A}^0} \|_{H^s(\R^d) \to L_2(\R^d)}
\le \wh{\mathfrak C}_3(s;\tau) \eps^{s/2},
$$
\textit{where}
$$
\wh{\mathfrak C}_3(s;\tau)= 2^{1-s/2} (\wh{\mathcal C}_5 + \wh{\mathcal C}_6 |\tau|)^{s/2}.
\eqno(13.16)
$$
\textit{In particular, for $0< \eps \le 1$ and $|\tau|= \eps^{-\alpha}$, $0< \alpha <1$, we have}
$$
\| e^{- i \tau \wh{\mathcal A}_\eps} - e^{- i \tau  \wh{\mathcal A}^0} \|_{H^s(\R^d) \to L_2(\R^d)}
\le \wh{\mathfrak C}_3(s;1) \eps^{s(1-\alpha)/2},\quad
0< \eps \le 1, \quad |\tau|=\eps^{-\alpha}.
$$

\smallskip

Theorem 13.7 and Corollary 9.8 imply the following statement.

\smallskip\noindent\textbf{Corollary 13.8.}
\textit{Suppose that the matrices $b(\bt)$ and $g(\x)$ have real entries. Suppose that the spectrum of the germ  $\wh{S}(\bt)$ is simple for all
$\bt\in {\mathbb S}^{d-1}$. Then the statements of Theorem} 13.7 \textit{are valid}.

\smallskip\noindent\textbf{13.5. The sharpness of the result.}
Applying Theorem~12.4, we confirm the sharpness of the result of Theorem~13.1 in the general case.

\smallskip\noindent\textbf{Theorem 13.9.}  \textit{Suppose that the assumptions of Theorem} 13.1
\textit{are satisfied. Let $\wh{N}_0(\bt)$ be the operator defined by} (8.21).  \textit{Suppose that}
$\wh{N}_0(\bt_0)\ne 0$ \textit{for some $\bt_0 \in {\mathbb S}^{d-1}$. Let $0\ne \tau \in \R$. Then for any $1\le s<3$
it is impossible that the estimate}
$$
\|\left( e^{- i \tau  \wh{\mathcal A}_\eps} - e^{- i \tau \wh{\mathcal A}^0}\right) ({\mathcal H}_0+I)^{-s/2} \|_{L_2(\R^d) \to L_2(\R^d)}
\le {\mathcal C}(\tau)\eps
\eqno(13.17)
$$
\textit{holds for all sufficiently small $\eps >0$.}

\smallskip\noindent\textbf{Proof.} We prove by contradiction. Let us fix $\tau \ne 0$.
Suppose that for some $1\le s < 3$ there exists a constant ${\mathcal C}(\tau)$ such that estimate (13.17) holds for all sufficiently small $\eps$.
Then, by (13.5), estimate (12.4) is also satisfied. But this contradicts to the statement of Theorem~12.4. $\bullet$

\smallskip\noindent\textbf{13.6. Approximation of the operator $f^\eps e^{-i \tau {\A}_\eps} (f^\eps)^{-1}$.}
Now we proceed to the study of the operator ${\A}_\eps$ (see (13.1)).  Let ${\A}^0$ be defined by (10.3).
Using (13.2) for the operators ${\A}_\eps$ and ${\A}^0$ and taking (13.3) and (13.4) into account,
we obtain the following analog of identity (13.5):
$$
\begin{aligned}
&\left( f^\eps e^{-i \tau {\A}_\eps} (f^\eps)^{-1}- f_0 e^{-i \tau {\A}^0 } f_0^{-1} \right) ({\mathcal H}_0+I)^{-s/2}
\cr
&= T_\eps^*
\left( f e^{-i \tau \eps^{-2}{\A}} f^{-1} - f_0 e^{-i \tau \eps^{-2}{\A}^0 } f_0^{-1} \right) {\mathcal R}(\eps)^{s/2} T_\eps,\quad \eps>0.
\end{aligned}
\eqno(13.18)
$$

Since $T_\eps$ is unitary, combining this with Theorem~12.5, we obtain the following result (which has been proved before in [BSu5, Theorem 12.3]).

\smallskip\noindent\textbf{Theorem 13.10.}  \textit{Let
${\mathcal A}_\eps= (f^\eps(\x))^*b(\D)^* g^\eps b(\D) f^\eps(\x)$,
where $g(\x)$, $f(\x)$ and $b(\D)$ satisfy the assumptions of Subsection {\rm 6.1}.
Let ${\mathcal A}^0= f_0 b(\D)^* g^0 b(\D) f_0$, where $g^0$ is the effective matrix {\rm (8.9)} and $f_0 = (\underline{ff^*})^{1/2}$. Let ${\mathcal H}_0=-\Delta$.}
\textit{Then for $\tau \in \R$ and $\eps>0$ we have}
$$
\|\left( f^\eps e^{- i \tau {\mathcal A}_\eps} (f^\eps)^{-1} - f_0 e^{- i \tau  {\mathcal A}^0} f_0^{-1} \right) ({\mathcal H}_0 +I)^{-3/2}\|_{L_2(\R^d) \to L_2(\R^d)}
\le ({\mathcal C}_1 + {\mathcal C}_2 |\tau|) \eps.
\eqno(13.19)
$$
\textit{The constants ${\mathcal C}_1$ and ${\mathcal C}_2$ are defined by} (11.5)
\textit{and depend only on} $r_0$, $\alpha_0$, $\alpha_1$, $\|g\|_{L_\infty}$, $\| g^{-1}\|_{L_\infty}$, $\|f\|_{L_\infty}$, \textit{and} $\| f^{-1}\|_{L_\infty}$.

\smallskip

Obviously, by (10.2),
$$
\begin{aligned}
\| f^\eps e^{- i \tau {\mathcal A}_\eps} (f^\eps)^{-1} - f_0 e^{- i \tau  {\mathcal A}^0} f_0^{-1}\|_{L_2(\R^d) \to L_2(\R^d)}
\le 2 \|f\|_{L_\infty} \|f^{-1}\|_{L_\infty},
\cr \tau \in \R,\quad \eps>0.
\end{aligned}
\eqno(13.20)
$$
Interpolating between (13.20) and (13.19), we arrive at the following result which has been proved before in  [BSu5, Theorem 12.4].

\smallskip\noindent\textbf{Theorem 13.11.}  \textit{Suppose that the assumptions of Theorem} 13.10
\textit{are satisfied. Then for $0\le s \le 3$,  $\tau \in \R$, and $\eps>0$ we have}
$$
\| f^\eps e^{- i \tau {\mathcal A}_\eps} (f^\eps)^{-1}-  f_0 e^{- i \tau  {\mathcal A}^0} f_0^{-1} \|_{H^s(\R^d) \to L_2(\R^d)}
\le {\mathfrak C}_1(s;\tau) \eps^{s/3},
$$
\textit{where}
$$
{\mathfrak C}_1(s;\tau)= (2\|f\|_{L_\infty} \|f^{-1}\|_{L_\infty})^{1-s/3} ({\mathcal C}_1 + {\mathcal C}_2 |\tau|)^{s/3}.
\eqno(13.21)
$$
\textit{In particular, for $0< \eps \le 1$ and $|\tau|= \eps^{-\alpha}$, $0< \alpha <1$, we have}
$$
\begin{aligned}
\| f^\eps e^{- i \tau {\mathcal A}_\eps}(f^\eps)^{-1} - f_0 e^{- i \tau  {\mathcal A}^0} f_0^{-1} \|_{H^s(\R^d) \to L_2(\R^d)}
\le {\mathfrak C}_1(s;1) \eps^{s(1-\alpha)/3},
\cr
0< \eps \le 1, \quad |\tau|=\eps^{-\alpha}.
\end{aligned}
$$

\smallskip\noindent\textbf{13.7. Refinement of approximation of the operator $f^\eps e^{-i \tau {\A}_\eps} (f^\eps)^{-1}$ under the additional assumptions.}
Using (13.18) and Theorem~12.6, we deduce  the following result.

\smallskip\noindent\textbf{Theorem 13.12.}  \textit{Suppose that the assumptions of Theorem} 13.10
\textit{are satisfied. Let $\wh{N}_Q(\bt)$ be the operator defined by} (10.10), (10.11). \textit{Suppose that} $\wh{N}_Q(\bt)=0$ \textit{for all} $\bt \in {\mathbb S}^{d-1}$.
\textit{Then for $\tau \in \R$ and $\eps>0$ we have}
$$
\|\left( f^\eps e^{- i \tau {\mathcal A}_\eps} (f^\eps)^{-1} - f_0 e^{- i \tau  {\mathcal A}^0} f_0^{-1} \right) ({\mathcal H}_0 +I)^{-1}\|_{L_2(\R^d) \to L_2(\R^d)}
\le ({\mathcal C}_3 + {\mathcal C}_4 |\tau|) \eps.
\eqno(13.22)
$$
\textit{The constants ${\mathcal C}_3$ and ${\mathcal C}_4$ are defined by} (11.6)
\textit{and depend only on} $r_0$, $\alpha_0$, $\alpha_1$, $\|g\|_{L_\infty}$, $\| g^{-1}\|_{L_\infty}$, $\|f\|_{L_\infty}$, \textit{and} $\| f^{-1}\|_{L_\infty}$.

\smallskip

Interpolating between (13.20) and (13.22), we obtain the following statement.

\smallskip\noindent\textbf{Theorem 13.13.}  \textit{Suppose that the assumptions of Theorem} 13.12
\textit{are satisfied. Then for $0\le s \le 2$, $\tau \in \R$, and $\eps>0$ we have}
$$
\| f^\eps e^{- i \tau {\mathcal A}_\eps} (f^\eps)^{-1}-  f_0 e^{- i \tau  {\mathcal A}^0} f_0^{-1} \|_{H^s(\R^d) \to L_2(\R^d)}
\le {\mathfrak C}_2(s;\tau) \eps^{s/2},
$$
\textit{where}
$$
{\mathfrak C}_2(s;\tau)= (2 \|f\|_{L_\infty} \|f^{-1}\|_{L_\infty})^{1-s/2} ({\mathcal C}_3 + {\mathcal C}_4 |\tau|)^{s/2}.
\eqno(13.23)
$$
\textit{In particular, for $0< \eps \le 1$ and $|\tau|= \eps^{-\alpha}$, $0< \alpha <1$, we have}
$$
\begin{aligned}
\| f^\eps e^{- i \tau {\mathcal A}_\eps} (f^\eps)^{-1} - f_0 e^{- i \tau  {\mathcal A}^0} f_0^{-1} \|_{H^s(\R^d) \to L_2(\R^d)}
\le {\mathfrak C}_2(s;1) \eps^{s(1-\alpha)/2},
\cr
0< \eps\le 1, \quad |\tau|=\eps^{-\alpha}.
\end{aligned}
$$

\smallskip
From Theorem~13.13 and  Proposition~10.1 we deduce the following corollary.

\smallskip\noindent\textbf{Corollary 13.14.}
\textit{Suppose that at least one of the following conditions is fulfilled}:

\noindent
$1^\circ$. \textit{The operator ${\mathcal A}_\eps$ has the form
 $\wh{\mathcal A}_\eps= (f^\eps)^* \D^* g^\eps(\x)\D f^\eps(\x)$, where $g(\x)$ is a symmetric matrix
 with real entries. }

\noindent
$2^\circ$. \textit{Relations} (8.14) \textit{are satisfied, i.~e.} $g^0 = \overline{g}$.

\noindent\textit{Then the statements of Theorem} 13.13 \textit{are valid}.

\smallskip
Finally, Theorem 12.7 and (13.18) imply the following result.

\smallskip\noindent\textbf{Theorem 13.15.}  \textit{Suppose that the assumptions of Theorem} 13.10
\textit{are satisfied. Suppose that Condition} 11.3 (\textit{or more restrictive Condition} 11.4)
 \textit{is satisfied. Then for $\tau \in \R$ and $\eps>0$ we have}
$$
\|\left( f^\eps e^{- i \tau {\mathcal A}_\eps} (f^\eps)^{-1}-
f_0 e^{- i \tau  {\mathcal A}^0} f_0^{-1} \right) ({\mathcal H}_0 +I)^{-1}\|_{L_2(\R^d) \to L_2(\R^d)}
\le ({\mathcal C}_5 + {\mathcal C}_6 |\tau|) \eps.
\eqno(13.24)
$$
\textit{The constants ${\mathcal C}_5$ and ${\mathcal C}_6$ are defined in Theorem}~11.6
\textit{and depend only on} $r_0$, $\alpha_0$, $\alpha_1$, $\|g\|_{L_\infty}$, $\| g^{-1}\|_{L_\infty}$, $\|f\|_{L_\infty}$, $\| f^{-1}\|_{L_\infty}$,
\textit{and also on the number ${c}^\circ$ defined by} (11.8).

\smallskip

Interpolating between (13.20) and (13.24), we obtain the following theorem.

\smallskip\noindent\textbf{Theorem 13.16.}  \textit{Suppose that the assumptions of Theorem} 13.15
\textit{are satisfied. Then for $0\le s \le 2$,  $\tau \in \R$, and $\eps>0$ we have}
$$
\| f^\eps e^{- i \tau {\mathcal A}_\eps} (f^\eps)^{-1}- f_0 e^{- i \tau  {\mathcal A}^0} f_0^{-1} \|_{H^s(\R^d) \to L_2(\R^d)}
\le {\mathfrak C}_3(s;\tau) \eps^{s/2},
$$
\textit{where }
$$
{\mathfrak C}_3(s;\tau)= (2 \|f\|_{L_\infty} \|f^{-1}\|_{L_\infty})^{1-s/2}  ({\mathcal C}_5 + {\mathcal C}_6 |\tau|)^{s/2}.
\eqno(13.25)
$$
\textit{In particular, for $0< \eps \le 1$ and $|\tau|= \eps^{-\alpha}$, $0< \alpha <1$, we have}
$$
\begin{aligned}
\| f^\eps e^{- i \tau {\mathcal A}_\eps} (f^\eps)^{-1} -
f_0 e^{- i \tau  {\mathcal A}^0} f_0^{-1} \|_{H^s(\R^d) \to L_2(\R^d)}
\le {\mathfrak C}_3(s;1) \eps^{s(1-\alpha)/2},
\cr
0< \eps \le 1, \quad |\tau|=\eps^{-\alpha}.
\end{aligned}
$$

\smallskip

Theorem 13.16 and Corollary 11.7 imply the following statement.

\smallskip\noindent\textbf{Corollary 13.17.}
\textit{Suppose that the matrices $b(\bt)$, $g(\x)$, and $Q(\x)= (f(\x)f(\x)^*)^{-1}$ have real entries.
Suppose that the spectrum of the generalized spectral problem}~(10.7) \textit{is simple for all $\bt\in {\mathbb S}^{d-1}$. Then the statements of Theorem} 13.16
\textit{are valid}.

\smallskip\noindent\textbf{13.8. The sharpness of the result.}
Applying Theorem~12.8, we confirm the sharpness of the result of Theorem~13.10 in the general case.

\smallskip\noindent\textbf{Theorem 13.18.}  \textit{Suppose that the assumptions of Theorem} 13.10
\textit{are satisfied. Let $\wh{N}_{0,Q}(\bt)$ be defined in Subsection} 10.3.  \textit{Suppose that}
\hbox{$\wh{N}_{0,Q}(\bt_0)\ne 0$} \textit{for some  $\bt_0 \in {\mathbb S}^{d-1}$. Let $0\ne \tau \in \R$. Then for any
$1\le s<3$ it is impossible that the estimate}
$$
\|\left( f^\eps e^{- i \tau  {\mathcal A}_\eps} (f^\eps)^{-1} -f_0 e^{- i \tau {\mathcal A}^0} f_0^{-1}\right) ({\mathcal H}_0+I)^{-s/2} \|_{L_2(\R^d) \to L_2(\R^d)}
\le {\mathcal C}(\tau)\eps
\eqno(13.26)
$$
\textit{holds for all sufficiently small $\eps >0$.}

\smallskip\noindent\textbf{Proof.} We prove by contradiction. Let us fix $\tau \ne 0$.
Suppose that for some $1\le s <3$  there exists a constant ${\mathcal C}(\tau)$ such that estimate (13.26) holds for all sufficiently small $\eps$.
Then, by (13.18), estimate (12.6) also holds. But this contradicts the statement of Theorem 12.8. $\bullet$

\section*{§14. Homogenization of the Cauchy problem for the Schr\"odinger type equation}

\smallskip\noindent\textbf{14.1. The Cauchy problem for the homogeneous equation with the operator $\wh{\mathcal A}_\eps$.}
Let $\u_\eps(\x,\tau)$, $\x \in \R^d$, $\tau \in \R$, be the solution of the Cauchy problem
$$
i \frac{\partial \u_\eps(\x,\tau)}{\partial \tau} = b(\D)^* g^\eps(\x)b(\D) \u_\eps(\x,\tau),\quad \u_\eps(\x,0) = \bphi(\x),
\eqno(14.1)
$$
where $\bphi \in L_2(\R^d;\C^n)$ is a given function.
The solution can be represented as $\u_\eps(\cdot,\tau) = e^{-i \tau \wh{\A}_\eps} \bphi$.
Let $\u_0(\x,\tau)$ be the solution of the ``homogenized'' Cauchy problem
$$
i \frac{\partial \u_0 (\x,\tau)}{\partial \tau} = b(\D)^* g^0 b(\D) \u_0 (\x,\tau),\quad \u_0(\x,0) = \bphi(\x),
\eqno(14.2)
$$
where $g^0$ is the effective matrix. Then $\u_0 = e^{-i \tau \wh{\mathcal A}^0} \phi$.

Theorem 13.2 directly implies the following result which has been proved before in [BSu5, Theorem~14.1].

\smallskip\noindent\textbf{Theorem 14.1.}  \textit{Let $\u_\eps$ be the solution of problem} (14.1),
\textit{and let $\u_0$ be the solution of problem} (14.2).

\noindent $1^\circ$. \textit{If $\bphi \in H^s(\R^d;\C^n)$, $0 \le s \le 3$, then for $\tau \in \R$ and $\eps>0$ we have}
$$
\| \u_\eps(\cdot,\tau) - \u_0(\cdot,\tau) \|_{L_2(\R^d)} \le \eps^{s/3} \wh{\mathfrak C}_1(s;\tau) \|\bphi \|_{H^s(\R^d)}.
$$
\textit{In particular, for $0< \eps \le 1$ and} $\tau = \pm \eps^{-\alpha}$, $0< \alpha <1$,
$$
\| \u_\eps(\cdot,\pm\eps^{-\alpha}) - \u_0(\cdot,\pm\eps^{-\alpha}) \|_{L_2(\R^d)} \le \eps^{s(1-\alpha)/3} \wh{\mathfrak C}_1(s;1) \|\bphi \|_{H^s(\R^d)}.
\eqno(14.3)
$$
\textit{The constant $\wh{\mathfrak C}_1(s;\tau)$ is defined by} (13.11).

\noindent $2^\circ$. \textit{If $\bphi \in L_2(\R^d;\C^n)$,  then}
$$
\begin{aligned}
&\lim_{\eps \to 0} \| \u_\eps(\cdot,\tau) - \u_0(\cdot,\tau)\|_{L_2(\R^d)} =0,\quad \tau \in \R;
\cr
&\lim_{\eps \to 0} \| \u_\eps(\cdot, \pm\eps^{-\alpha}) - \u_0(\cdot,\pm \eps^{-\alpha})\|_{L_2(\R^d)} =0,\quad 0< \alpha <1.
\end{aligned}
$$

\smallskip
Statement $2^\circ$ follows directly from statement $1^\circ$ and the Banach-Steinhaus theorem.

Statement $1^\circ$ can be refined under the additional assumptions.
Theorem 13.4 implies the following statement.

\smallskip\noindent\textbf{Theorem 14.2.}  \textit{Suppose that the assumptions of Theorem} 14.1
\textit{are satisfied. Let $\wh{N}(\bt)$ be the operator defined by} (8.18), (8.19). \textit{Suppose that} $\wh{N}(\bt)=0$ \textit{for all} $\bt \in {\mathbb S}^{d-1}$.
\textit{If $\bphi \in H^s(\R^d;\C^n)$, $0 \le s \le 2$, then for $\tau \in \R$ and $\eps>0$ we have}
$$
\| \u_\eps(\cdot,\tau) - \u_0(\cdot,\tau) \|_{L_2(\R^d)} \le \eps^{s/2} \wh{\mathfrak C}_2(s;\tau) \|\bphi \|_{H^s(\R^d)}.
$$
\textit{In particular, for $0< \eps \le 1$ and} $\tau = \pm \eps^{-\alpha}$, $0< \alpha <1$,
$$
\| \u_\eps(\cdot,\pm\eps^{-\alpha}) - \u_0(\cdot,\pm\eps^{-\alpha}) \|_{L_2(\R^d)} \le \eps^{s(1-\alpha)/2} \wh{\mathfrak C}_2(s;1) \|\bphi \|_{H^s(\R^d)}.
\eqno(14.4)
$$
\textit{The constant $\wh{\mathfrak C}_2(s;\tau)$ is given by} (13.14).

\smallskip

Finally, Theorem~13.7 implies the following result.

\smallskip\noindent\textbf{Theorem 14.3.}  \textit{Suppose that the assumptions of Theorem} 14.1
\textit{are satisfied. Suppose also that Condition} 9.3 (\textit{or more restrictive Condition} 9.4)
\textit{is satisfied. If $\bphi \in H^s(\R^d;\C^n)$, $0 \le s \le 2$, then for $\tau \in \R$ and $\eps>0$ we have}
$$
\| \u_\eps(\cdot,\tau) - \u_0(\cdot,\tau) \|_{L_2(\R^d)} \le \eps^{s/2} \wh{\mathfrak C}_3(s;\tau) \|\bphi \|_{H^s(\R^d)}.
$$
\textit{In particular, for $0< \eps \le 1$ and} $\tau = \pm \eps^{-\alpha}$, $0< \alpha <1$,
$$
\| \u_\eps(\cdot,\pm\eps^{-\alpha}) - \u_0(\cdot,\pm\eps^{-\alpha}) \|_{L_2(\R^d)} \le \eps^{s(1-\alpha)/2} \wh{\mathfrak C}_3(s;1) \|\bphi \|_{H^s(\R^d)}.
$$
\textit{The constant $\wh{\mathfrak C}_3(s;\tau)$ is defined by} (13.16).

\smallskip

\smallskip\noindent\textbf{14.2. The Cauchy problem for the nonhomogeneous equation with the operator $\wh{\mathcal A}_\eps$.}
Now we consider the Cauchy problem for the nonhomogeneous equation
$$
i \frac{\partial \u_\eps(\x,\tau)}{\partial \tau} = b(\D)^* g^\eps(\x)b(\D) \u_\eps(\x,\tau) + \FF(\x,\tau),\quad \u_\eps(\x,0) = \bphi(\x),
\eqno(14.5)
$$
 where $\bphi \in L_2(\R^d;\C^n)$ and $\FF \in L_{1,\text{loc}}(\R; L_2(\R^d;\C^n))$ are given functions.
 The solution of problem (14.5) can be represented as
$$
\u_\eps(\cdot,\tau) = e^{-i \tau \wh{\A}_\eps} \bphi - i \intop_0^\tau e^{-i (\tau- \wt{\tau}) \wh{\A}_\eps} \FF(\cdot, \wt{\tau})\, d\wt{\tau}.
\eqno(14.6)
$$
Let $\u_0(\x,\tau)$ be the solution of the homogenized problem
$$
i \frac{\partial \u_0(\x,\tau)}{\partial \tau} = b(\D)^* g^0 b(\D) \u_0(\x,\tau) + \FF(\x,\tau),\quad \u_0(\x,0) = \bphi(\x).
\eqno(14.7)
$$
Then
$$
\u_0(\cdot,\tau) = e^{-i \tau \wh{\A}^0} \bphi - i \intop_0^\tau e^{-i (\tau- \wt{\tau}) \wh{\A}^0} \FF(\cdot, \wt{\tau})\, d\wt{\tau}.
\eqno(14.8)
$$

The following theorem has been proved before in [BSu5, Theorem~14.2].

\smallskip\noindent\textbf{Theorem 14.4.}  \textit{Let $\u_\eps$ be the solution of problem} (14.5),
\textit{and let $\u_0$ be the solution of problem}~(14.7).

\noindent $1^\circ$.
\textit{If $\bphi \in H^s(\R^d;\C^n)$ and} $\FF \in L_{1,\text{loc}}(\R; H^s(\R^d;\C^n))$ \textit{with some \hbox{$0 \le s \le 3$}, then for $\tau \in \R$ and $\eps>0$ we have}
$$
\| \u_\eps(\cdot,\tau) - \u_0(\cdot,\tau) \|_{L_2(\R^d)} \le \eps^{s/3} \wh{\mathfrak C}_1(s;\tau) \left( \|\bphi \|_{H^s(\R^d)} +
\|\FF\|_{L_1((0,\tau);H^s(\R^d))} \right).
\eqno(14.9)
$$
\textit{Under the additional assumption that $\FF \in  L_{p}(\R_\pm; H^s(\R^d;\C^n))$, where $p \in [1,\infty]$, for  $0< \eps \le 1$
and $\tau = \pm \eps^{-\alpha}$,
$0< \alpha < s (s+3/p')^{-1}$, we have}
$$
\begin{aligned}
&\| \u_\eps(\cdot,\pm\eps^{-\alpha}) - \u_0(\cdot,\pm\eps^{-\alpha}) \|_{L_2(\R^d)}
\cr
&\le \eps^{s(1-\alpha)/3} \wh{\mathfrak C}_1(s;1) \left(\|\bphi \|_{H^s(\R^d)} + \eps^{-\alpha / p'} \|\FF\|_{L_p(\R_\pm; H^s(\R^d))} \right).
\end{aligned}
\eqno(14.10)
$$
\textit{The constant $\wh{\mathfrak C}_1(s;\tau)$ is defined by} (13.11). \textit{Here $p^{-1}+ (p')^{-1}=1$.}

\noindent $2^\circ$.
\textit{If $\bphi \in L_2(\R^d;\C^n)$ and} $\FF \in L_{1,\text{loc}}(\R; L_2(\R^d;\C^n))$, \textit{then}
$$
\lim_{\eps \to 0} \| \u_\eps(\cdot,\tau) - \u_0(\cdot,\tau) \|_{L_2(\R^d)} =0,\quad \tau \in \R.
$$
\textit{Under the additional assumption that $\FF \in L_{1}(\R_\pm; L_2(\R^d;\C^n))$, we have}
$$
\lim_{\eps \to 0} \| \u_\eps(\cdot, \pm \eps^{-\alpha}) - \u_0(\cdot, \pm \eps^{-\alpha}) \|_{L_2(\R^d)} =0,\quad 0< \alpha <1.
$$

\smallskip\noindent\textbf{Proof.}
If $\bphi \in H^s(\R^d;\C^n)$ and $\FF \in  L_{1,\text{loc}}(\R; H^s(\R^d;\C^n))$ with some \hbox{$0 \le s \le 3$}, then relations (13.10), (14.6), and (14.8) imply (14.9).

If  $\FF \in  L_{p}(\R_\pm; H^s(\R^d;\C^n))$ with some $p \in [1,\infty]$, then for
$\tau = \pm \eps^{-\alpha}$,  $0< \alpha < s (s+3/p')^{-1}$, relations (13.12), (14.6), and (14.8) yield (14.10).

If it is known only that $\bphi \in L_2(\R^d;\C^n)$ and $\FF \in  L_{1,\text{loc}}(\R; L_2(\R^d;\C^n))$, then from the obvious estimate
$$
\| \u_\eps(\cdot,\tau) - \u_0(\cdot,\tau) \|_{L_2(\R^d)} \le 2\|\bphi \|_{L_2(\R^d)}
+ 2 \|\FF\|_{L_1((0,\tau);L_2(\R^d))}
$$
and (14.9), by the Banach-Steinhaus theorem, it follows that $\u_\eps(\cdot,\tau)$ tends to $\u_0(\cdot,\tau)$ in  $L_2(\R^d;\C^n)$
for a fixed $\tau \in \R$.

If $\FF \in  L_{1}(\R_\pm; L_2(\R^d;\C^n))$, then for $\tau = \pm \eps^{-\alpha}$, $0< \alpha < 1$, from the obvious estimate
$$
\| \u_\eps(\cdot,\pm \eps^{-\alpha}) - \u_0(\cdot, \pm \eps^{-\alpha})  \|_{L_2(\R^d)} \le 2\|\bphi \|_{L_2(\R^d)}
+ 2 \|\FF\|_{L_1(\R_\pm;L_2(\R^d))}
\eqno(14.11)
$$
and (14.10) (with $p=1$), by the Banach-Steinhaus theorem, it follows that the left-hand side of (14.11) tends to zero.
$\bullet$

\smallskip

Statement $1^\circ$ of Theorem 14.4 can be refined under the additional assumptions.
Theorem~13.4 implies the following result.

\smallskip\noindent\textbf{Theorem 14.5.}  \textit{Suppose that the assumptions of Theorem}~14.4
\textit{are satisfied. Let $\wh{N}(\bt)$ be the operator defined by} (8.18), (8.19). \textit{Suppose that} \hbox{$\wh{N}(\bt)=0$} \textit{for all $\bt \in {\mathbb S}^{d-1}$.
If $\bphi \in H^s(\R^d;\C^n)$ and} $\FF \in L_{1,\text{loc}}(\R; H^s(\R^d;\C^n))$ \textit{with some $0 \le s \le 2$, then for $\tau \in \R$ and $\eps>0$ we have}
$$
\| \u_\eps(\cdot,\tau) - \u_0(\cdot,\tau) \|_{L_2(\R^d)} \le \eps^{s/2} \wh{\mathfrak C}_2(s;\tau) \left( \|\bphi \|_{H^s(\R^d)} +
\|\FF\|_{L_1((0,\tau);H^s(\R^d))} \right).
$$
\textit{Under the additional assumption that $\FF \in  L_{p}(\R_\pm; H^s(\R^d;\C^n))$, where $p \in [1,\infty]$, for
 $0< \eps \le 1$ and $\tau = \pm \eps^{-\alpha}$, $0< \alpha < s (s+2/p')^{-1}$, we have}
$$
\begin{aligned}
&\| \u_\eps(\cdot,\pm\eps^{-\alpha}) - \u_0(\cdot,\pm\eps^{-\alpha}) \|_{L_2(\R^d)}
\cr
&\le \eps^{s(1-\alpha)/2} \wh{\mathfrak C}_2(s;1) \left(\|\bphi \|_{H^s(\R^d)} + \eps^{-\alpha / p'} \|\FF\|_{L_p(\R_\pm; H^s(\R^d))} \right).
\end{aligned}
$$
\textit{The constant $\wh{\mathfrak C}_2(s;\tau)$ is defined by} (13.14).

\smallskip

Finally, Theorem 13.7 implies the following result.

\smallskip\noindent\textbf{Theorem 14.6.}  \textit{Suppose that the assumptions of Theorem} 14.4
\textit{are satisfied. Suppose that Condition} 9.3 (\textit{or more restrictive Condition}~9.4)
\textit{is satisfied. If $\bphi \in H^s(\R^d;\C^n)$ and} $\FF \in L_{1,\text{loc}}(\R; H^s(\R^d;\C^n))$ \textit{with some $0 \le s \le 2$, then for $\tau \in \R$ and $\eps>0$ we have}
$$
\| \u_\eps(\cdot,\tau) - \u_0(\cdot,\tau) \|_{L_2(\R^d)} \le \eps^{s/2} \wh{\mathfrak C}_3(s;\tau) \left( \|\bphi \|_{H^s(\R^d)} +
\|\FF\|_{L_1((0,\tau);H^s(\R^d))} \right).
$$
\textit{Under the additional assumption that $\FF \in  L_{p}(\R_\pm; H^s(\R^d;\C^n))$, where $p \in [1,\infty]$, for  $0< \eps \le 1$ and $\tau = \pm \eps^{-\alpha}$,
$0< \alpha < s (s+2/p')^{-1}$, we have}
$$
\begin{aligned}
&\| \u_\eps(\cdot,\pm\eps^{-\alpha}) - \u_0(\cdot,\pm\eps^{-\alpha}) \|_{L_2(\R^d)}
\cr
&\le \eps^{s(1-\alpha)/2} \wh{\mathfrak C}_3(s;1) \left(\|\bphi \|_{H^s(\R^d)} + \eps^{-\alpha / p'} \|\FF\|_{L_p(\R_\pm; H^s(\R^d))} \right).
\end{aligned}
$$
\textit{The constant $\wh{\mathfrak C}_3(s;\tau)$ is defined by} (13.16).

\smallskip\noindent\textbf{14.3. The Cauchy problem for the homogeneous equation with the operator ${\mathcal A}_\eps$.}
Let $\A_\eps$ be the operator (13.1). Consider the Cauchy problem
$$
i \frac{\partial \u_\eps(\x,\tau)}{\partial \tau} = (f^\eps(\x))^* b(\D)^* g^\eps(\x)b(\D) f^\eps(\x) \u_\eps(\x,\tau),\quad f^\eps(\x)\u_\eps(\x,0) = \bphi(\x),
\eqno(14.12)
$$
where $\bphi \in L_2(\R^d;\C^n)$.
The solution of problem (14.12) can be represented as $\u_\eps(\cdot,\tau) = e^{-i \tau {\A}_\eps} (f^\eps)^{-1} \bphi$.
Let $\u_0(\x,\tau)$ be the solution of the ``homogenized'' Cauchy problem
$$
i \frac{\partial \u_0 (\x,\tau)}{\partial \tau} = f_0 b(\D)^* g^0 b(\D) f_0 \u_0 (\x,\tau),\quad f_0 \u_0(\x,0) = \bphi(\x),
\eqno(14.13)
$$
 where $g^0$ is the effective matrix (8.9) and $f_0$ is defined by (10.1).
 Then $\u_0 = e^{-i \tau {\mathcal A}^0} f_0^{-1}\bphi$.

Theorem 13.11 implies the following result (which has been proved before in  [BSu5, Theorem 14.3]).

\smallskip\noindent\textbf{Theorem 14.7.}  \textit{Let $\u_\eps$ be the solution of problem} (14.12),
\textit{and let $\u_0$ be the solution of problem} (14.13).

\noindent $1^\circ$. \textit{If $\bphi \in H^s(\R^d;\C^n)$, $0 \le s \le 3$, then for $\tau \in \R$ and $\eps>0$ we have}
$$
\| f^\eps \u_\eps(\cdot,\tau) - f_0 \u_0(\cdot,\tau) \|_{L_2(\R^d)} \le \eps^{s/3} {\mathfrak C}_1(s;\tau) \|\bphi \|_{H^s(\R^d)}.
$$
\textit{In particular, for $0< \eps \le 1$ and} $\tau = \pm \eps^{-\alpha}$, $0< \alpha <1$,
$$
\| f^\eps \u_\eps(\cdot,\pm\eps^{-\alpha}) - f_0 \u_0(\cdot,\pm\eps^{-\alpha}) \|_{L_2(\R^d)} \le \eps^{s(1-\alpha)/3} {\mathfrak C}_1(s;1) \|\bphi \|_{H^s(\R^d)}.
$$
\textit{The constant ${\mathfrak C}_1(s;\tau)$ is defined by} (13.21).

\noindent $2^\circ$. \textit{If $\bphi \in L_2(\R^d;\C^n)$,  then}
$$
\begin{aligned}
&\lim_{\eps \to 0} \| f^\eps \u_\eps(\cdot,\tau) - f_0 \u_0(\cdot,\tau)\|_{L_2(\R^d)} =0,\quad \tau \in \R;
\cr
&\lim_{\eps \to 0} \| f^\eps \u_\eps(\cdot, \pm\eps^{-\alpha}) - f_0 \u_0(\cdot,\pm \eps^{-\alpha})\|_{L_2(\R^d)} =0,\quad 0< \alpha <1.
\end{aligned}
$$

\smallskip
Statement $1^\circ$ of Theorem 14.7 can be refined under the additional assumptions. Theorem~13.13 implies the following result.

\smallskip\noindent\textbf{Theorem 14.8.}  \textit{Suppose that the assumptions of Theorem} 14.7
\textit{are satisfied. Let $\wh{N}_Q(\bt)$ be the operator defined by} (10.10), (10.11). \textit{Suppose that} \hbox{$\wh{N}_Q(\bt)=0$}
\textit{for all} $\bt \in {\mathbb S}^{d-1}$.
\textit{If $\bphi \in H^s(\R^d;\C^n)$, $0 \le s \le 2$, then for $\tau \in \R$ and $\eps>0$ we have}
$$
\| f^\eps \u_\eps(\cdot,\tau) - f_0 \u_0(\cdot,\tau) \|_{L_2(\R^d)} \le \eps^{s/2} {\mathfrak C}_2(s;\tau) \|\bphi \|_{H^s(\R^d)}.
$$
\textit{In particular, for $0< \eps \le 1$ and} $\tau = \pm \eps^{-\alpha}$, $0< \alpha <1$,
$$
\| f^\eps \u_\eps(\cdot,\pm\eps^{-\alpha}) - f_0 \u_0(\cdot,\pm\eps^{-\alpha}) \|_{L_2(\R^d)} \le \eps^{s(1-\alpha)/2} {\mathfrak C}_2(s;1) \|\bphi \|_{H^s(\R^d)}.
$$
\textit{The constant ${\mathfrak C}_2(s;\tau)$ is defined by} (13.23).

\smallskip
Similarly, Theorem 13.16 implies the following result.

\smallskip\noindent\textbf{Theorem 14.9.}  \textit{Suppose that the assumptions of Theorem} 14.7
\textit{are satisfied. Suppose that Condition} 11.3 (\textit{or more restrictive Condition} 11.4)
\textit{is satisfied. If $\bphi \in H^s(\R^d;\C^n)$, $0 \le s \le 2$, then for $\tau \in \R$ and $\eps>0$ we have}
$$
\| f^\eps \u_\eps(\cdot,\tau) - f_0 \u_0(\cdot,\tau) \|_{L_2(\R^d)} \le \eps^{s/2} {\mathfrak C}_3(s;\tau) \|\bphi \|_{H^s(\R^d)}.
$$
\textit{In particular, for $0< \eps \le 1$ and} $\tau = \pm \eps^{-\alpha}$, $0< \alpha <1$,
$$
\| f^\eps \u_\eps(\cdot,\pm\eps^{-\alpha}) - f_0 \u_0(\cdot,\pm\eps^{-\alpha}) \|_{L_2(\R^d)} \le \eps^{s(1-\alpha)/2} {\mathfrak C}_3(s;1) \|\bphi \|_{H^s(\R^d)}.
$$
\textit{The constant ${\mathfrak C}_3(s;\tau)$ is defined by} (13.25).

\smallskip\noindent\textbf{14.4. The Cauchy problem for the nonhomogeneous equation with the operator ${\mathcal A}_\eps$.}
Now we consider the Cauchy problem for the nonhomogeneous  equation:
$$
\begin{aligned}
&i \frac{\partial \u_\eps(\x,\tau)}{\partial \tau} =(f^\eps(\x))^*  b(\D)^* g^\eps(\x)b(\D) f^\eps(\x) \u_\eps(\x,\tau) + (f^\eps(\x))^{-1}\FF(\x,\tau),
\cr
&f^\eps(\x)\u_\eps(\x,0) = \bphi(\x),
\end{aligned}
\eqno(14.14)
$$
where $\bphi \in L_2(\R^d;\C^n)$ and $\FF \in L_{1,\text{loc}}(\R; L_2(\R^d;\C^n))$.
The solution of problem (14.14) can be represented as
$$
\u_\eps(\cdot,\tau) = e^{-i \tau {\A}_\eps} ( f^\eps)^{-1} \bphi - i \intop_0^\tau e^{-i (\tau- \wt{\tau}) {\A}_\eps} (f^\eps)^{-1} \FF(\cdot, \wt{\tau})\, d\wt{\tau}.
\eqno(14.15)
$$
Let $\u_0(\x,\tau)$ be the solution of the homogenized problem
$$
i \frac{\partial \u_0(\x,\tau)}{\partial \tau} = f_0 b(\D)^* g^0 b(\D) f_0 \u_0(\x,\tau) + f_0^{-1}\FF(\x,\tau),\quad f_0 \u_0(\x,0) = \bphi(\x).
\eqno(14.16)
$$
Then
$$
\u_0(\cdot,\tau) = e^{-i \tau {\A}^0} f_0^{-1}\bphi - i \intop_0^\tau e^{-i (\tau- \wt{\tau}) {\A}^0} f_0^{-1}\FF(\cdot, \wt{\tau})\, d\wt{\tau}.
\eqno(14.17)
$$

By analogy with the proof of Theorem 14.4, from Theorem 13.11 and relations (14.15), (14.17) we deduce the following result (which has been proved before in [BSu5, Theorem~14.5]).

\smallskip\noindent\textbf{Theorem 14.10.}  \textit{Let $\u_\eps$ be the solution of problem} (14.14),
\textit{and let $\u_0$ be the solution of problem} (14.16).

\noindent $1^\circ$. \textit{If $\bphi \in H^s(\R^d;\C^n)$  and} $\FF \in L_{1,\text{loc}}(\R; H^s(\R^d;\C^n))$ \textit{with some \hbox{$0 \le s \le 3$}, then for $\tau \in \R$ and $\eps>0$ we have}
$$
\| f^\eps \u_\eps(\cdot,\tau) - f_0 \u_0(\cdot,\tau) \|_{L_2(\R^d)} \le \eps^{s/3} {\mathfrak C}_1(s;\tau) \left(\|\bphi \|_{H^s(\R^d)} + \|\FF\|_{L_1((0,\tau);H^s(\R^d))}\right).
$$
\textit{Under the additional assumption that $\FF \in L_{p}(\R_\pm; H^s(\R^d;\C^n))$, where $p \in [1,\infty]$,
for $0< \eps \le 1$ and} $\tau = \pm \eps^{-\alpha}$, $0< \alpha < s (s+3/p')^{-1}$, \textit{we have}
$$
\begin{aligned}
&\| f^\eps \u_\eps(\cdot,\pm\eps^{-\alpha})  - f_0 \u_0(\cdot,\pm\eps^{-\alpha}) \|_{L_2(\R^d)}
\cr
&\le \eps^{s(1-\alpha)/3} {\mathfrak C}_1(s;1)
\left( \|\bphi \|_{H^s(\R^d)} + \eps^{-\alpha/p'}\|\FF\|_{L_p(\R_\pm ;H^s(\R^d))}\right).
\end{aligned}
$$
\textit{The constant ${\mathfrak C}_1(s;\tau)$ is defined by} (13.21).

\noindent $2^\circ$. \textit{If $\bphi \in L_2(\R^d;\C^n)$  and} $\FF \in L_{1,\text{loc}}(\R; L_2(\R^d;\C^n))$, \textit{then}
$$
\lim_{\eps \to 0} \| f^\eps \u_\eps(\cdot,\tau) - f_0 \u_0(\cdot,\tau)\|_{L_2(\R^d)} =0,\quad \tau \in \R.
$$
\textit{Under the additional assumption that $\FF \in L_{1}(\R_\pm; L_2(\R^d;\C^n))$, we have}
$$
\lim_{\eps \to 0} \| f^\eps \u_\eps(\cdot, \pm\eps^{-\alpha}) - f_0 \u_0(\cdot,\pm \eps^{-\alpha})\|_{L_2(\R^d)} =0,\quad 0< \alpha <1.
$$

\smallskip

Statement $1^\circ$ of Theorem 14.10 can be refined under the additional assumptions. Theorem~13.13 implies the following result.

\smallskip\noindent\textbf{Theorem 14.11.}  \textit{Suppose that the assumptions of Theorem}~14.10 \textit{are satisfied. Let $\wh{N}_Q(\bt)$ be the operator defined by} (10.10), (10.11).
\textit{Suppose that} $\wh{N}_Q(\bt)=0$ \textit{for all $\bt \in {\mathbb S}^{d-1}$.}
\textit{If $\bphi \in H^s(\R^d;\C^n)$  and} $\FF \in L_{1,\text{loc}}(\R; H^s(\R^d;\C^n))$ \textit{with some $0 \le s \le 2$, then for $\tau \in \R$ and $\eps>0$ we have}
$$
\| f^\eps \u_\eps(\cdot,\tau) - f_0 \u_0(\cdot,\tau) \|_{L_2(\R^d)} \le \eps^{s/2} {\mathfrak C}_2(s;\tau) \left(\|\bphi \|_{H^s(\R^d)} + \|\FF\|_{L_1((0,\tau);H^s(\R^d))}\right).
$$
\textit{Under the additional assumption that $\FF \in L_{p}(\R_\pm; H^s(\R^d;\C^n))$, where $p \in [1,\infty]$,
for $0< \eps \le 1$ and} $\tau = \pm \eps^{-\alpha}$, $0< \alpha < s (s+2/p')^{-1}$, \textit{we have}
$$
\begin{aligned}
&\| f^\eps \u_\eps(\cdot,\pm\eps^{-\alpha})  - f_0 \u_0(\cdot,\pm\eps^{-\alpha}) \|_{L_2(\R^d)}
\cr
&\le \eps^{s(1-\alpha)/2} {\mathfrak C}_2(s;1)
\left( \|\bphi \|_{H^s(\R^d)} + \eps^{-\alpha/p'} \|\FF\|_{L_p(\R_\pm ;H^s(\R^d))}\right).
\end{aligned}
$$
\textit{The constant ${\mathfrak C}_2(s;\tau)$ is defined by} (13.23).

\smallskip
Similarly, applying Theorem 13.16, we deduce the following result.

\smallskip\noindent\textbf{Theorem 14.12.}  \textit{Suppose that the assumptions of Theorem} 14.10
\textit{are satisfied. Suppose also that Condition} 11.3  (\textit{or more restrictive Condition} 11.4)
\textit{is satisfied. If $\bphi \in H^s(\R^d;\C^n)$  and} $\FF \in L_{1,\text{loc}}(\R; H^s(\R^d;\C^n))$ \textit{with some $0 \le s \le 2$, then for $\tau \in \R$ and $\eps>0$ we have}
$$
\| f^\eps \u_\eps(\cdot,\tau) - f_0 \u_0(\cdot,\tau) \|_{L_2(\R^d)} \le \eps^{s/2} {\mathfrak C}_3(s;\tau) \left(\|\bphi \|_{H^s(\R^d)} + \|\FF\|_{L_1((0,\tau);H^s(\R^d))}\right).
$$
\textit{Under the additional assumption that $\FF \in L_{p}(\R_\pm; H^s(\R^d;\C^n))$, where $p \in [1,\infty]$,
for $0< \eps \le 1$ and} $\tau = \pm \eps^{-\alpha}$, $0< \alpha < s (s+2/p')^{-1}$, \textit{we have}
$$
\begin{aligned}
&\| f^\eps \u_\eps(\cdot,\pm\eps^{-\alpha})  - f_0 \u_0(\cdot,\pm\eps^{-\alpha}) \|_{L_2(\R^d)}
\cr
&\le \eps^{s(1-\alpha)/2} {\mathfrak C}_3(s;1)
\left( \|\bphi \|_{H^s(\R^d)} + \eps^{-\alpha/p'} \|\FF\|_{L_p(\R_\pm ;H^s(\R^d))}\right).
\end{aligned}
$$
\textit{The constant ${\mathfrak C}_3(s;\tau)$ is defined by} (13.25).

\section*{§15. Application of the general results: the nonstationary Schr\"odinger equation}

\smallskip\noindent\textbf{15.1. The model example: the Schr\"odinger type equation with the operator
$\wh{\mathcal A}_\eps = -{\rm div}\, g^\eps(\x) \nabla$.} In $L_2(\R^d)$, $d \ge 1$, we consider the operator
$$
\wh{\mathcal A} = \D^* g(\x)\D = -{\rm div}\, g(\x) \nabla.
\eqno(15.1)
$$
Here $g(\x)$ is a $\Gamma$-periodic Hermitian $(d \times d)$-matrix-valued function such that
$$
g(\x) >0;\quad g, g^{-1} \in L_\infty.
$$
The operator (15.1) is a particular case of the operator (8.1). We have $n=1$, $m=d$, and $b(\D) = \D$.
Obviously, condition (6.2) is satisfied with  $\alpha_0 = \alpha_1 =1$.
According to (8.11), the effective operator for the operator (15.1) is given by
$$
\wh{\mathcal A}^0 = \D^* g^0 \D =-{\rm div}\, g^0 \nabla.
$$
By the general rule, the effective matrix $g^0$ is defined as follows.
Let $\e_1,\dots, \e_d$ be the standard orthonormal basis in $\R^d$.
Let $\Phi_j \in \wt{H}^1(\Omega)$ be the weak $\Gamma$-periodic solution of the problem
$$
{\rm div}\, g(\x) (\nabla \Phi_j(\x) + \e_j) =0,\quad \intop_\Omega \Phi_j(\x)\,d\x=0.
\eqno(15.2)
$$
Then $g^0$ is the $(d\times d)$-matrix with the columns
$$
{\mathbf g}^0_j = |\Omega|^{-1} \intop_\Omega g(\x) (\nabla \Phi_j(\x) + \e_j)\,d\x,\quad j=1,\dots,d.
$$
If  $d=1$, then $m=n=1$, whence $g^0 = \underline{g}$.

If $g(\x)$ is a symmetric matrix with real entries, then, by Proposition~8.4($1^\circ$),  $\wh{N}(\bt)=0$ for all $\bt \in  {\mathbb S}^{d-1}$.
If $g(\x)$ is a Hermitian matrix with complex entries, then, in general, $\wh{N}(\bt)$ is not zero.
Now $n=1$, therefore, $\wh{N}(\bt) = \wh{N}_0(\bt)$ is the operator of multiplication by $\wh{\mu}(\bt)$,
where $\wh{\mu}(\bt)$ is the coefficient in the expansion for the first eigenvalue
$\wh{\lambda}(t,\bt) = \wh{\gamma}(\bt)t^2 +  \wh{\mu}(\bt)t^3+\dots$ of the operator $\wh{\A}(\k)=\wh{A}(t,\bt)$.
Calculation (see [BSu3, Subsection~10.3]) shows that
$$
\begin{aligned}
\wh{N}(\bt) =& \wh{\mu}(\bt) = - i \sum_{j,l,k=1}^d (a_{jlk} - a_{jlk}^*) \theta_j \theta_l \theta_k, \quad \bt \in {\mathbb S}^{d-1},
\cr
a_{jlk}:=&  |\Omega|^{-1} \intop_\Omega \Phi_j(\x)^* \langle g(\x) (\nabla \Phi_l(\x) + \e_l ), \e_k \rangle \,d\x, \quad j,l,k=1,\dots,d.
\end{aligned}
\eqno(15.3)
$$

The following example is borrowed from  [BSu3, Subsection~10.4].

\smallskip\noindent\textbf{Example 15.1.}
Let $d=2$ and $\Gamma = (2\pi \Z)^2$. Suppose that the matrix $g(\x)$ is given by
$$
g(\x) = \begin{pmatrix}
1 & i \beta'(x_1) \cr - i \beta'(x_1) & 1
\end{pmatrix},
$$
where $\beta(x_1)$ is a smooth $(2\pi)$-periodic real-valued function such that  $\int_0^{2\pi} \beta(x_1)\,dx_1=0$  and
$1 - (\beta'(x_1))^2 >0$. Then $\wh{N}(\bt) = - \alpha \pi^{-1} \theta_2^3$, where $\alpha = \int_0^{2\pi} \beta(x_1) (\beta'(x_1))^2\,dx_1$.
It is easy to give a concrete example where \hbox{$\alpha \ne 0$}: if $\beta(x_1) = c (\sin x_1 + \cos 2 x_1)$ with $0< c < 1/3$, then $\alpha = - (3\pi/2) c^3 \ne 0$.
In this example, $\wh{N}(\bt)= \wh{\mu}(\bt) \ne 0$ for all $\bt\in {\mathbb S}^1$ except for the points $(\pm 1,0)$.

\smallskip

Consider the Cauchy problem
$$
i \frac{\partial u_\eps(\x,\tau)}{\partial \tau} = \D^* g^\eps(\x) \D u_\eps(\x,\tau),\quad u_\eps(\x,0) = \phi(\x),
\eqno(15.4)
$$
where $\phi \in L_2(\R^d)$ is a given function. Let $u_0(\x,\tau)$ be the solution of the ``homogenized'' Cauchy problem
$$
i \frac{\partial u_0(\x,\tau)}{\partial \tau} = \D^* g^0  \D u_0 (\x,\tau),\quad u_0(\x,0) = \phi(\x).
\eqno(15.5)
$$
 Applying Theorem 14.1 and, in the ``real'' case,  applying Theorem 14.2, we arrive at  the following statement.

\smallskip\noindent\textbf{Proposition 15.2.}  \textit{Suppose that the assumptions of Subsection} 15.1
\textit{are satisfied. Let $u_\eps$ be the solution of problem} (15.4),  \textit{and let $u_0$ be the solution of problem} (15.5).

\noindent $1^\circ$. \textit{If $\phi \in H^s(\R^d)$ for some $0 \le s \le 3$, then for $\tau \in \R$ and $\eps>0$ we have}
$$
\|  u_\eps(\cdot,\tau) -  u_0(\cdot,\tau) \|_{L_2(\R^d)} \le \eps^{s/3} \wh{\mathfrak C}_1(s;\tau) \|\phi \|_{H^s(\R^d)},
$$
\textit{where the constant $\wh{\mathfrak C}_1(s;\tau)$ is given by} (13.11). \textit{If $\phi \in L_2(\R^d)$, then}
$$
\lim_{\eps \to 0} \|  u_\eps(\cdot,\tau) -  u_0(\cdot,\tau) \|_{L_2(\R^d)}=0,\quad \tau \in \R.
$$

\noindent $2^\circ$. \textit{Let $g(\x)$ be a symmetric matrix with real entries. If $\phi \in H^s(\R^d)$ with some $0 \le s \le 2$, then for  $\tau \in \R$ and $\eps>0$ we have}
$$
\|  u_\eps(\cdot,\tau) -  u_0(\cdot,\tau) \|_{L_2(\R^d)} \le \eps^{s/2} \wh{\mathfrak C}_2(s;\tau) \|\phi \|_{H^s(\R^d)},
$$
\textit{where the constant $\wh{\mathfrak C}_2(s;\tau)$ is given by} (13.14).

\smallskip
One can also apply the statements about the behavior of the solution for $|\tau| = \eps^{-\alpha}$
with $0< \alpha < 1$ (estimates of the form (14.3) in the general case and (14.4) in the ``real'' case).
It is possible to consider more general problem for the nonhomogeneous equation and apply Theorem~14.4 in the general case and Theorem~14.5 in the ``real'' case.

\smallskip\noindent\textbf{15.2. The periodic Schr\"odinger operator. Factorization.} (See [BSu1, Chapter~6, Subsection~1.1].)
In $L_2(\R^d)$, $d\ge 1$, we consider the Schr\"odinger operator
$$
{\mathcal H} = \D^* \check{g}(\x)\D + V(\x)
\eqno(15.6)
$$
with the $\Gamma$-periodic metric $\check{g}(\x)$ and potential $V(\x)$.
It is assumed that $\check{g}(\x)$ is a symmetric  $(d\times d)$-matrix-valued function with real entries, $V(\x)$ is a real-valued function, and
$$
\check{g}(\x) >0;\quad \check{g}, \check{g}^{-1} \in L_\infty,
$$
$$
V \in L_q(\Omega),\quad 2q >d \ \text{for}\ d\ge 2; \ \ q=1 \ \text{for}\ d=1.
\eqno(15.7)
$$
 The precise definition of the operator ${\mathcal H}$ is given in terms of the quadratic form
$$
h[u,u] = \intop_{\R^d} \left( \langle \check{g}(\x)\D u, \D u\rangle + V(\x)|u|^2 \right)\,d\x,\quad u \in H^1(\R^d),
\eqno(15.8)
$$
which, under our assumptions, is closed and lower semibounded. Adding an appropriate constant to $V(\x)$, we assume that the
\textit{point $\lambda_0=0$ is the bottom of the spectrum of} $\mathcal H$.

	 Under our assumptions, the equation $\D^* \check{g}(\x)\D \omega(\x) + V(\x) \omega(\x)=0$ has a positive $\Gamma$-periodic solution $\omega \in \wt{H}^1(\Omega)$.
	 Moreover, $\omega$ is a multiplier in $H^1(\R^d)$ and $\wt{H}^1(\Omega)$.
	 We fix the choice of $\omega$ by the normalization condition $\int_\Omega \omega^2(\x)\,d\x =|\Omega|$.
 After the substitution $u = \omega v$, the form  (15.8) turns into
 $$
h[u,u] = \intop_{\R^d} \omega^2(\x) \langle \check{g}(\x)\D v, \D v\rangle \,d\x,\quad u= \omega v,\quad v  \in H^1(\R^d).
$$
This means that the operator (15.6) admits the following factorization
$$
{\mathcal H} = \omega^{-1} \D^* \omega^2 \check{g} \D \omega^{-1}.
\eqno(15.9)
$$
 Thus, the operator $\mathcal H$ is represented in the form (6.4) with $n=1$, $m=d$, $b(\D) = \D$, $g = \omega^2 \check{g}$, and $f = \omega^{-1}$.

\smallskip\noindent\textbf{Remark 15.3.} The expression  (15.9) can be taken as the definition of the operator $\mathcal H$
for any $\Gamma$-periodic function  $\omega$ such that
$\omega(\x)>0$; $\omega,\omega^{-1} \in L_\infty$. The form  (15.6) can be recovered by the formula $V = - \omega^{-1} (\D^* \check{g} \D \omega)$.
The corresponding potential $V$ may be a singular distribution.

\smallskip
The operator  (15.9) and the operator (15.1) (with $g = \omega^2 \check{g}$) satisfy the identity ${\mathcal H}= \omega^{-1} \wh{\A} \omega^{-1}$.
Let $g^0$ be the effective matrix for the operator~(15.1). The function $Q= (ff^*)^{-1}$ takes the form $Q(\x)= \omega^2(\x)$.
By the normalization condition for $\omega$, we have $\overline{Q}=1$ and $f_0 = (\overline{Q})^{-1/2}= 1$.
Therefore, the operator~(10.3) takes the form
$$
{\mathcal H}^0 = \D^* g^0 \D.
$$
By Proposition 10.1($1^\circ$), the operator $\wh{N}_Q(\bt)$ is equal to zero: $\wh{N}_Q(\bt)=0$ for all $\bt \in {\mathbb S}^{d-1}$.

\smallskip\noindent\textbf{15.3. The nonstationary Schr\"odinger equation with a singular potential.}
Now we consider the operator
$$
{\mathcal H}_\eps = (\omega^\eps)^{-1} \D^* g^\eps \D (\omega^\eps)^{-1},\quad g^\eps = (\omega^\eps)^2 \check{g}^\eps.
\eqno(15.10)
$$
Under condition (15.7), the operator (15.10) can be written in the initial terms:
$$
{\mathcal H}_\eps =  \D^* \check{g}^\eps \D  + \eps^{-2} V^\eps.
\eqno(15.11)
$$
Note that the expression (15.11) contains a large factor  $\eps^{-2}$ at the rapidly oscillating potential  $V^\eps$.

We consider the Cauchy problem of the form (14.12):
$$
i \frac{\partial u_\eps(\x,\tau)}{\partial \tau} = {\mathcal H}_\eps  u_\eps(\x,\tau),\quad (\omega^\eps (\x))^{-1}u_\eps(\x,0) = \phi(\x),
\eqno(15.12)
$$
where $\phi \in L_2(\R^d)$. Let $u_0(\x,\tau)$ be the solution of the homogenized problem (see (14.13))
$$
i \frac{\partial u_0(\x,\tau)}{\partial \tau} = {\mathcal H}^0  u_0(\x,\tau),\quad u_0(\x,0) = \phi(\x).
\eqno(15.13)
$$

Applying Theorem 14.8, we arrive at the following resullt.

\smallskip\noindent\textbf{Proposition 15.4.}  \textit{Suppose that the assumptions of Subsections} 15.2 \textit{and} 15.3
\textit{are satisfied. Let $u_\eps$ be the solution of problem} (15.12),  \textit{and let $u_0$ be the solution of problem} (15.13).
\textit{If $\phi \in H^s(\R^d)$ for some $0 \le s \le 2$, then for $\tau \in \R$ and $\eps>0$ we have}
$$
\|  (\omega^\eps)^{-1} u_\eps(\cdot,\tau) -  u_0(\cdot,\tau) \|_{L_2(\R^d)} \le \eps^{s/2} {\mathfrak C}_2(s;\tau) \|\phi \|_{H^s(\R^d)},
$$
\textit{where the constant ${\mathfrak C}_2(s;\tau)$ is given by} (13.23). \textit{If $\phi \in L_2(\R^d)$, then}
$$
\lim_{\eps \to 0} \|  (\omega^\eps)^{-1} u_\eps(\cdot,\tau) -  u_0(\cdot,\tau) \|_{L_2(\R^d)}=0,\quad \tau \in \R.
$$

\smallskip
One can also apply the statement of Theorem~14.8 about the behavior of the solutions for $|\tau| = \eps^{-\alpha}$
with $0< \alpha <1$. It is also possible to consider more general Cauchy problem for the nonhomogeneous equation and to apply Theorem~14.11.

\smallskip\noindent\textbf{15.4. The nonstationary Schr\"odinger equation with a magnetic potential.}
In $L_2(\R^d)$, $d \ge 2$, we consider the periodic magnetic Schr\"odinger operator $\mathcal M$ with $\Gamma$-periodic metric $\check{g}(\x)$, magnetic potential ${\mathbf A}(\x)$,
and electric potential $V(\x)$:
$$
{\mathcal M} = (\D - {\mathbf A}(\x))^* \check{g}(\x) (\D - {\mathbf A}(\x)) + V(\x).
$$
Here $\check{g}(\x)$ is a symmetric  $(d\times d)$-matrix-valued function with real entries such that
$\check{g}(\x)>0$ and $\check{g},\check{g}^{-1}\in L_\infty$. If $d\ge 3$, we assume in addition that $\check{g} \in C^\alpha$ with some $0< \alpha <1$.
Suppose that ${\mathbf A}(\x)$ is a $\R^d$-valued function and $V(\x)$ is a real-valued function such that
$$
{\mathbf A} \in L_{2q}(\Omega),\quad V \in L_q(\Omega),\quad 2q >d.
$$
As usual, the precise definition of the operator is given in terms of the corresponding quadratic form.
Adding an appropriate constant to $V(\x)$, we assume that the bottom of the spectrum of $\mathcal M$ is the point $\lambda_0=0$.

According to  [Sh2], under the above assumptions and for sufficiently small (in the $L_{2q}(\Omega)$-norm) magnetic potential $\mathbf A$,
the operator $\mathcal M$ admits an appropriate factorization. Let us describe this factorization. Let ${\mathcal M}(\k)$ be
the family of operators in  $L_2(\Omega)$
that arise in the direct integral expansion for $\mathcal M$. The condition $\inf {\rm spec} \, {\mathcal M}=0$ means that for some $\k_0 \in \wt{\Omega}$
the point $\lambda_0=0$ is an eigenvalue of ${\mathcal M}(\k_0)$. If the magnetic potential is sufficiently small, then this point
$\k_0$ is unique and the eigenvalue $\lambda_0=0$ of the operator ${\mathcal M}(\k_0)$ is simple. Let $\eta(\x)$ be the corresponding eigenfunction normalized
by the condition $\int_\Omega |\eta(\x)|^2 \,d\x = |\Omega|$ (the phase factor is not important).
Then $\eta \in \wt{H}^1(\Omega)$, and $\eta, \eta^{-1} \in L_\infty$. As shown in [Sh2], $\eta$ is a multiplier in $H^1(\R^d)$ and  $\wt{H}^1(\Omega)$.
We denote
$$
\wt{\mathcal M} := [e^{-i\langle \k_0,\cdot \rangle}] {\mathcal M} [e^{i\langle \k_0,\cdot \rangle}].
$$
Clearly, the coefficients of the operator $\wt{\mathcal M}$ are periodic.
By Theorems 2.7 and 2.8 from [Sh2], if the norm $\|{\mathbf A}\|_{L_{2q}(\Omega)}$ is sufficiently small, then the operator $\wt{\mathcal M}$ admits the following factorization:
$$
\wt{\mathcal M} = (\eta(\x)^*)^{-1} \D^* g(\x) \D \eta(\x)^{-1}.
\eqno(15.14)
$$
Here the Hermitian $\Gamma$-periodic matrix-valued function $g(\x)$ is defined by
$$
g(\x) = |\eta(\x)|^2 \check{g}(\x) + i g_2(\x),
\eqno(15.15)
$$
and the antisymmetric matrix-valued function $g_2(\x)$ with real entries satisfies the equation
$$
({\rm div}\, g_2(\x))^t = -2 |\eta(\x)|^2 \check{g}(\x) ({\mathbf A}(\x) - \k_0) + 2\, {\rm Im}\, (\eta(\x)^* \check{g}(\x) \nabla \eta(\x)).
\eqno(15.16)
$$
As shown in [Sh2], we have
$$
g(\x)>0; \quad g, g^{-1} \in L_\infty.
$$
The operator~(15.14) is of the form~(6.4) with $n=1$, $m=d$, $b(\D)=\D$,  $g$ defined by (15.15), (15.16), and
$f = \eta^{-1}$. Let $g^0$ be the effective matrix for the operator $\wh{\mathcal A}= \D^* g \D$; in general, the effective matrix may have complex entries.
Now the function $Q = (ff^*)^{-1}$ takes the form $Q(\x)=|\eta(\x)|^2$. By the normalization condition on $\eta$, we have $\overline{Q}=1$, and then $f_0=1$.
The operator (10.3)  takes the form
$$
\wt{\mathcal M}^0 = \D^* g^0 \D.
$$
Let us describe the operator $\wh{N}_Q(\bt)$.
Let $\Phi_j$ be the $\Gamma$-periodic solution of problem (15.2). Since $n=1$ and $\overline{Q}=1$, then (see (5.9)) the operator $\wh{N}_Q(\bt)= \wh{N}_{0,Q}(\bt)$
acts as multiplication by $\mu(\bt)$, where $\mu(\bt)$ is the coefficient in the expansion
$\lambda(t,\bt) = \gamma(\bt) t^2 + \mu(\bt) t^3 + \dots$ for the first eigenvalue $\lambda(t,\bt)$ of the operator $\wt{\mathcal M}(\k)= \wt{\mathcal M}(t \bt)$.
A calculation shows that
$$
\begin{aligned}
\wh{N}_Q(\bt)= \mu(\bt) = -  i \sum_{j,l,k=1}^d (a_{jlk} - a_{jlk}^*) \theta_j \theta_l \theta_k +  2 \langle g^0 \bt, \bt \rangle \sum_{j=1}^d {\rm Im}\, (\overline{Q \Phi_j})\theta_j ,
\cr
\bt \in {\mathbb S}^{d-1},
\end{aligned}
$$
where the coefficients $a_{jlk}$ are defined by (15.3). In general, the operator $\wh{N}_Q(\bt)$ is not zero.

Now we consider the operators
$$
\wt{\mathcal M}_\eps =  ((\eta^\eps)^*)^{-1} \D^* g^\eps \D (\eta^\eps)^{-1}, \quad
{\mathcal M}_\eps = [e^{i \eps^{-1} \langle \k_0 ,\cdot \rangle}] \wt{\mathcal M}_\eps [e^{-i \eps^{-1} \langle \k_0 ,\cdot \rangle}].
$$
In the initial terms, we have
$$
{\mathcal M}_\eps =  (\D - \eps^{-1} {\mathbf A}^\eps)^* \check{g}^\eps (\D - \eps^{-1} {\mathbf A}^\eps) + \eps^{-2} V^\eps.
\eqno(15.17)
$$
Note that the expression (15.17) contains large factors $\eps^{-1}$ at the rapidly oscillating magnetic potential ${\mathbf A}^\eps$
and $\eps^{-2}$ at the electric potential $V^\eps$. Let $u_\eps(\x,\tau)$ be the solution of the Cauchy problem for the nonstationary magnetic Schr\"odinger equation:
$$
i \frac{\partial u_\eps(\x,\tau)}{\partial \tau} = {\mathcal M}_\eps  u_\eps(\x,\tau),
\quad (\eta^\eps(\x))^{-1} e^{-i \eps^{-1} \langle \k_0 ,\x \rangle} u_\eps(\x,0) = \phi(\x),
\eqno(15.18)
$$
 where  $\phi \in L_2(\R^d)$.
Then the function $v_\eps(\x,\tau)= e^{-i \eps^{-1} \langle \k_0 ,\x \rangle} u_\eps(\x,\tau)$ is the solution of the problem
$$
i \frac{\partial v_\eps(\x,\tau)}{\partial \tau} = \wt{\mathcal M}_\eps  v_\eps(\x,\tau),\quad (\eta^\eps(\x))^{-1} v_\eps(\x,0) = \phi(\x).
\eqno(15.19)
$$
We can apply Theorem 14.7. The effective problem is of the form
$$
i \frac{\partial v_0(\x,\tau)}{\partial \tau} = \D^* g^0 \D  v_0(\x,\tau),\quad  v_0(\x,0) = \phi(\x).
\eqno(15.20)
$$
Applying Theorem 14.7 to the problem (15.19),  we arrive at the following result.

\smallskip\noindent\textbf{Proposition 15.5.}  \textit{Suppose that the assumptions of Subsection} 15.4 \textit{are satisfied.
Let $u_\eps$ be the solution of problem} (15.18),  \textit{and let $v_0$ be the solution of problem} (15.20).
\textit{If $\phi \in H^s(\R^d)$ with some $0 \le s \le 3$, then for $\tau \in \R$ and $\eps>0$ we have}
$$
\|  (\eta^\eps)^{-1} e^{-i \eps^{-1} \langle \k_0 ,\x \rangle} u_\eps(\cdot,\tau) -  v_0(\cdot,\tau) \|_{L_2(\R^d)} \le \eps^{s/3} {\mathfrak C}_1(s;\tau) \|\phi \|_{H^s(\R^d)},
$$
\textit{where the constant ${\mathfrak C}_1(s;\tau)$ is given by} (13.21). \textit{If $\phi \in L_2(\R^d)$, then}
$$
\lim_{\eps \to 0} \|  (\eta^\eps)^{-1} e^{-i \eps^{-1} \langle \k_0 ,\x \rangle} u_\eps(\cdot,\tau) -  v_0(\cdot,\tau) \|_{L_2(\R^d)}=0,\quad \tau \in \R.
$$

\smallskip
One can also apply the statemment of Theorem 14.7 about the behavior of the solutions for $|\tau| = \eps^{-\alpha}$ with $0< \alpha <1$.
It is also possible to consider more general  Cauchy problem for the nonhomogeneous equation and to apply Theorem~14.10.

\smallskip\noindent\textbf{Remark 15.6.}  In [Sh1] it was shown that  in general (without the smallness condition on  $\mathbf A$)
the required factorization for the magnetic  Schr\"odinger operator is not valid. This leads to interesting effects in the corresponding
homogenization problem (see [Sh3]).

\section*{§16. Application of the general results: the nonstationary two-dimensional Pauli equation}

\smallskip\noindent\textbf{16.1. Definition and factorization of the two-dimensional Pauli operator.}
(See [BSu1, Chapter~6, Subsection~2.1].) Suppose that the magnetic potential is a vector-valued function
${\mathbf A}(\x)= \{A_1(\x),A_2(\x)\}$ in $\R^2$, where $A_j(\x)$ are $\Gamma$-periodic real-valued functions such that
$$
A_j \in L_\rho(\Omega),\quad \rho>2,\quad j=1,2.
\eqno(16.1)
$$
  Recall the standard notation for the Pauli matrices
  $$
  \sigma_1 = \begin{pmatrix} 0 & 1 \cr 1 & 0 \end{pmatrix},\quad
  \sigma_2 = \begin{pmatrix} 0 & -i \cr i & 0 \end{pmatrix},\quad
  \sigma_3 = \begin{pmatrix} 1 & 0 \cr 0 & -1 \end{pmatrix}.
  $$
In $L_2(\R^2;\C^2)$, we consider the operator
$$
{\mathcal D} = (D_1 - A_1)\sigma_1 + (D_2 - A_2)\sigma_2,\quad {\rm Dom}\, {\mathcal D} = H^1(\R^2;\C^2).
\eqno(16.2)
$$
By definition, the Pauli operator is the square of  $\mathcal D$:
$$
{\mathcal P}:= {\mathcal D}^2= \begin{pmatrix} P_- & 0 \cr 0 & P_+ \end{pmatrix}.
\eqno(16.3)
$$
Precisely,  $\mathcal P$ is the selfadjoint operator in  $L_2(\R^2;\C^2)$
corresponding to the closed quadratic form  $\| {\mathcal D}\u \|^2_{L_2(\R^2)}$, $\u \in H^1(\R^2;\C^2)$.
If ${\mathbf A}(\x)$ is Lipschitz, then the blocks $P_\pm$ of the operator (16.3) can be written as
$$
P_\pm = ( \D - {\mathbf A}(\x))^2 \pm B(\x),\quad B(\x):= \partial_1 A_2(\x) - \partial_2 A_1(\x).
$$

We use a well known factorization for the  two-dimensional Pauli operator. A gauge transformation allows us to assume that
the potential $\mathbf A$ is subject to the conditions
$$
{\rm div}\, {\mathbf A}(\x)=0, \quad \intop_\Omega {\mathbf A}(\x) \, d\x =0,
\eqno(16.4)
$$
and still satisfies (16.1). Under conditions (16.1) and (16.4),  there exists a (unique) real-valued $\Gamma$-periodic function  $\varphi$ such that

 $$
 \nabla \varphi(\x) = \{ A_2(\x),- A_1(\x)\}, \quad \intop_\Omega \varphi(\x) \, d\x =0.
 $$
Note that  $\varphi \in \wt{W}^1_\rho(\Omega) \subset C^\sigma$, $\sigma = 1 - 2 \rho^{-1}$.
 We put
 $$
 \omega_\pm(\x) := e^{\pm \varphi(\x)}.
 $$
The operators (16.2), (16.3) admit the following factorization:
$$
{\mathcal D} = f_\times(\x) b_\times(\D) f_\times(\x),
\eqno(16.5)
$$
$$
{\mathcal P} = f_\times(\x) b_\times(\D) g_\times(\x) b_\times(\D) f_\times(\x),
\eqno(16.6)
$$
where
 $$
 \begin{aligned}
 b_\times(\D) &= \begin{pmatrix} 0 & D_1 - i D_2 \cr D_1+ i D_2 & 0 \end{pmatrix},
 \cr
 f_\times(\x) &= \begin{pmatrix} \omega_+(\x) &  0  \cr 0  & \omega_-(\x) \end{pmatrix},
 \quad g_\times(\x) = f_\times(\x)^2=
 \begin{pmatrix} \omega^2_+(\x) &  0  \cr 0  & \omega_-^2(\x) \end{pmatrix}.
 \end{aligned}
 $$
 The blocks $P_\pm$ of the operator (16.3) can be written as
$$
P_+ = \omega_- (D_1+i D_2) \omega_+^2 (D_1 - i D_2) \omega_-, \quad
P_- = \omega_+ (D_1- i D_2) \omega_-^2 (D_1 + i D_2) \omega_+.
\eqno(16.7)
$$

\smallskip\noindent\textbf{Remark 16.1.} 1) We may take expressions (16.5), (16.6), (16.7)
as the definition of the operators $\mathcal D$, $\mathcal P$, and $P_\pm$, assuming that   $\omega_\pm(\x)$
are arbitrary $\Gamma$-periodic functions satisfying the conditions $\omega_\pm(\x)>0$;
$\omega_+, \omega_- \in L_\infty$, and $\omega_+(\x) \omega_-(\x)=1$.
2) Note that the operators $P_+$ and $P_-$ are unitarily equivalent.
Moreover, the operators $P_+(\k)$ and $P_-(\k)$ in $L_2(\Omega)$
are also unitarily equivalent for each $\k$.

\smallskip\noindent\textbf{16.2. The effective characteristics for the operators $P_\pm$. Homogenization.}
   The operators $P_\pm$ are of the form (6.4) with
   $d=2$, $m=n=1$, $b(\D)= D_1 \mp i D_2$, $g(\x)= \omega^2_\pm(\x)$, and $f(\x)= \omega_\mp(\x)$.
   The role of the operator $\wh{\A}$ for $P_\pm$ is played by the operatoor $\wh{\mathcal A}_\pm = (D_1\pm i D_2) \omega_\pm^2 (D_1 \mp  i D_2)$.
   Since $m=n$, the effective constant is given by
   $$
   g^0_\pm = \underline{\omega_\pm^2} = \left( |\Omega|^{-1} \int_\Omega \omega_\mp^2(\x) \,d\x \right)^{-1} =: \omega^2_{\pm,0}.
   \eqno(16.8)
   $$
The role of $Q(\x)$ for the operator $P_\pm$ is played by  $Q_\pm(\x)= \omega^2_\pm(\x)$.
Then, by  (16.8), $\overline{Q_\pm}= (g^0_\mp)^{-1}$. The role of $f_0$ is played by the constant $(\overline{Q_\pm})^{-1/2}= (g^0_\mp)^{1/2}=\omega_{\mp,0}$.
Next, the role of $\A^0$ for $P_\pm$ is played by the operator $P_\pm^0$, where
$$
\begin{aligned}
 P_+^0 &= \omega_{-,0} (D_1+i D_2) g^0_+  (D_1 - i D_2)  \omega_{-,0} = -\gamma \Delta,
 \cr
 P_-^0 &=  \omega_{+,0} (D_1- i D_2) g^0_-  (D_1 + i D_2) \omega_{+,0} = -\gamma \Delta.
 \end{aligned}
 $$
 Here
 $$
 \gamma := g^0_+ g^0_- = |\Omega|^2 \| \omega_+ \|_{L_2(\Omega)}^{-2} \| \omega_- \|_{L_2(\Omega)}^{-2}.
 \eqno(16.9)
 $$
 Let $\lambda_\pm(t,\bt)$ be the analytic  (in $t$) branch of the first eigenvalue of the operator $P_\pm(\k)$, and let
$\lambda_\pm(t,\bt)= \gamma_\pm(\bt) t^2 + \mu_\pm(\bt)t^3+\dots$ be the corresponding power series expansion.
Since the operators $P_+(\k)$ and $P_-(\k)$ are unitarily equivalent, then $\lambda_+(t,\bt)=\lambda_-(t, \bt)$, and  also
$\gamma_+(\bt)= \gamma_-(\bt)$, $\mu_+(\bt) = \mu_-(\bt)$.
As shown in  [BSu1, Chapter 6, \S 2], the numbers  $\gamma_\pm(\bt)$ do not depend on  $\bt$ and are given by
$$
\gamma_+(\bt)=\gamma_- (\bt)= \gamma,
$$
where $\gamma$ is defined by (16.9).
Now we describe the operator $\wh{N}_{Q,\pm}(\bt)$ that plays the role of $\wh{N}_Q(\bt)$ for $P_\pm$.
Let $v_\pm(\x)$ be the  $\Gamma$-periodic solution of the problem
$$
(D_1 \mp i D_2) v_\pm(\x) = g^0_\pm \omega_\mp^{2}(\x) -1,\quad \intop_\Omega v_\pm(\x)=0.
\eqno(16.10)
$$
Then
$$
\wh{N}_{Q,\pm}(\bt) = - 2 \gamma \left( \theta_1 {\rm Re} \,\overline{\omega^2_\pm v_\pm}  \pm \theta_2 {\rm Im} \,\overline{\omega^2_\pm v_\pm} \right),\quad \bt \in {\mathbb S}^1.
\eqno(16.11)
$$
According to (10.12),
$$
\mu_\pm(\bt) = - 2  g^0_\mp \gamma \left( \theta_1 {\rm Re} \,\overline{\omega^2_\pm v_\pm}  \pm \theta_2 {\rm Im} \,\overline{\omega^2_\pm v_\pm} \right),\quad \bt \in {\mathbb S}^1.
\eqno(16.12)
$$
Though we know that $\mu_+(\bt)=\mu_-(\bt)=:\mu(\bt)$, it is not evident to deduce this directly from (16.12).
Now we give an example where $\mu(\bt) \ne 0$.

\smallskip\noindent\textbf{Example 16.2.} Let $\Gamma = (2\pi \Z)^2$ and let $\omega^2_-(\x)= 1+ \alpha (\sin x_2 + 4 \sin 2 x_2)$, where $\alpha >0$ is sufficiently small.
Then, according to (16.8) and (16.10), we have
$g^0_+=1$ and $v_+(\x) = \alpha (\cos x_2 + 2 \cos 2x_2)$. Let us calculate $\overline{\omega_+^2 v_+}$:
$$
\begin{aligned}
\overline{\omega_+^2 v_+}
&= \frac{\alpha}{2\pi}
\intop_0^{2\pi} \frac{\cos x + 2 \cos 2x}{1+ \alpha(\sin x + 4 \sin 2x)}\,dx
= \frac{3\alpha}{8\pi}
\intop_0^{2\pi} \frac{\cos x }{1+ \alpha(\sin x + 4 \sin 2x)}\,dx
\cr
&= - \frac{6 \alpha^2}{\pi}
\intop_{-1}^{1} \frac{t \sqrt{1-t^2}}{(1+ \alpha t)^2 - 64 \alpha^2 t^2 (1 -t^2) }\,dt
\cr
&=
 \frac{24 \alpha^3}{\pi}
\intop_{0}^{1} \frac{t^2 \sqrt{1-t^2}}{\left((1+ \alpha t)^2 - 64 \alpha^2 t^2 (1 -t^2)\right) \left((1- \alpha t)^2 - 64 \alpha^2 t^2 (1 -t^2)\right) }\,dt.
\end{aligned}
$$
Obviously, for sufficiently small $\alpha>0$ (for instance, one can take $\alpha = \frac{1}{16}$)
the function in the last integral is positive, whence $\overline{\omega_+^2 v_+}>0$.
Then $\mu(\bt)= - 2 g^0_- \gamma \theta_1 \overline{\omega_+^2 v_+} \ne 0$ for $\theta_1 \ne 0$.

\smallskip
Now we consider the operators
$$
\begin{aligned}
P_{+,\eps}= \omega_-^\eps (D_1+iD_2) (\omega^\eps_+)^2 (D_1 - i D_2) \omega^\eps_-,
\cr
P_{-,\eps}= \omega_+^\eps (D_1- iD_2) (\omega^\eps_-)^2 (D_1 + i D_2) \omega^\eps_+.
\end{aligned}
\eqno(16.13)
$$
If $\mathbf A$ is Lipschitz, then the operators (16.13) can be written as
$$
P_{\pm,\eps} = (\D - \eps^{-1} {\mathbf A}^\eps)^2 \pm \eps^{-2} B^\eps.
$$
We consider the following Cauchy problems for the scalar functions $u_{\pm,\eps}(\x,\tau)$:
$$
i \frac{\partial u_{\pm,\eps}(\x,\tau)}{\partial \tau} = {P}_{\pm,\eps}  u_{\pm,\eps}(\x,\tau),\quad
\omega_\mp^\eps(\x)  u_{\pm,\eps}(\x,0) = \phi_\pm(\x),
\eqno(16.14)
$$
 where $\phi_\pm \in L_2(\R^2)$. We can apply Theorem 14.7. The corresponding homogenized problems are
$$
i \frac{\partial u_{\pm,0}(\x,\tau)}{\partial \tau} =  - \gamma \Delta  u_{\pm,0}(\x,\tau),\quad \omega_{\mp,0}  u_{\pm,0}(\x,0) = \phi_\pm(\x),
\eqno(16.15)
$$
where $\gamma$ is defined by (16.9), and $\omega_{\pm,0}$ are given by  (16.8).

\smallskip\noindent\textbf{Proposition 16.3.}  \textit{Suppose that the assumptions of Subsections} 16.1, 16.2 are satisfied.
\textit{Let $u_{\pm,\eps}$ be the solution of problem} (16.14),  \textit{and let $u_{\pm,0}$ be the solution of problem} (16.15).
\textit{If $\phi_\pm \in H^s(\R^2)$ with some $0 \le s \le 3$, then for $\tau \in \R$ and $\eps>0$ we have}
$$
\|  \omega_\mp^\eps  u_{\pm,\eps}(\cdot,\tau) -  \omega_{\mp,0} u_{\pm,0} (\cdot,\tau) \|_{L_2(\R^2)} \le \eps^{s/3} {\mathfrak C}_\pm(s;\tau) \|\phi_\pm \|_{H^s(\R^2)},
$$
\textit{where the constants ${\mathfrak C}_\pm(s;\tau)$ are of the form} (13.21). \textit{If $\phi_\pm \in L_2(\R^2)$, then}
$$
\lim_{\eps \to 0}
\|  \omega_\mp^\eps  u_{\pm,\eps}(\cdot,\tau) -  \omega_{\mp,0}  u_{\pm,0}(\cdot,\tau) \|_{L_2(\R^2)} =0, \quad \tau \in \R.
$$

\smallskip
One can also apply the statement of  Theorem~14.7 about the behavior of the solutions for $|\tau| = \eps^{-\alpha}$ with $0<\alpha <1$.
It is also possible to consider more general Cauchy problem for the nonhomogeneous equation and apply  Theorem 14.10.

\smallskip\noindent\textbf{16.3. The effective characteristics for the operator $\mathcal P$. Homogenization.}
 The operator $\mathcal P$ is of the form (6.4) with
   $d=2$, $m=n=2$, $b(\D)= b_\times(\D)$, $g(\x)= g_\times(\x)$, and $f(\x)= f_\times(\x)$.
   The role of the operator $\wh{\A}$ for $\mathcal P$ is played by the operator $\wh{\A}_\times = b_\times(\D) g_\times(\x) b_\times(\D)$.
   Since $m=n$, the effective matrix is given by
   $$
   g^0_\times = \underline{g_\times} = {\rm diag}\, \{g^0_+, g^0_-\},
   $$
   where $g^0_\pm$ are defined by (16.8). The role of  $Q$ is played by the matrix $Q_\times = f_\times^{-2}= g_\times^{-1}$.
   Then $\overline{Q_\times} = {\rm diag}\, \{(g^0_+)^{-1}, (g^0_-)^{-1}\}$.
   The role of $f_0$ is played by the matrix
   $$
   f_{\times,0} = {\rm diag}\, \{\omega_{+,0}, \omega_{-,0}\}.
   \eqno(16.16)
   $$
   Next, the operator ${\A}^0$ (see (10.3)) takes the form
   $$
   {\mathcal P}^0 =  f_{\times,0} b_{\times}(\D) g^0_\times b_{\times}(\D) f_{\times,0} = \begin{pmatrix} - \gamma \Delta & 0 \cr 0 & - \gamma \Delta \end{pmatrix}.
   $$
Let us describe the operator $\wh{N}_{Q,\times}(\bt)$ that  plays the role of $\wh{N}_Q(\bt)$ for $\mathcal P$.
A calculation shows that
$$
\wh{N}_{Q,\times}(\bt)= {\rm diag}\, \{ \wh{N}_{Q,+} (\bt),  \wh{N}_{Q,-} (\bt)\},
$$
where the operators $\wh{N}_{Q,\pm}(\bt)$ are defined by (16.11).
The first eigenvalue $\lambda(t,\bt)$ (analytic in $t$) of the operator ${\mathcal P}(\k)$
is of multiplicity two identically in $\k=t\bt$, because the blocks  $P_+(\k)$ and $P_-(\k)$
are unitarily equivalent. In the power series expansion $\lambda(t,\bt) = \gamma t^2 + \mu(\bt) t^3+\dots$,
the coefficient $\gamma$ is independent of $\bt$ and defined by (16.9), and the coefficient  $\mu(\bt)$ is defined by  (16.12).
Example 16.2 shows that, in the general case, the coefficient $\mu(\bt)$ is not zero.

Now we consider the operator
$$
{\mathcal P}_{\eps}= f_\times^\eps b_\times (\D) g_\times^\eps b_\times (\D) f_\times^\eps  = \begin{pmatrix} P_{-,\eps} & 0 \cr 0 & P_{+,\eps}\end{pmatrix},
$$
where the blocks are defined by (16.13).
Consider the Cauchy problem for a vector-valued function  $\u_{\eps}(\x,\tau)$:
$$
i \frac{\partial \u_{\eps}(\x,\tau)}{\partial \tau} = {\mathcal P}_{\eps}  \u_{\eps}(\x,\tau),\quad f_\times^\eps(\x)  \u_{\eps}(\x,0) = \bphi(\x),
\eqno(16.17)
$$
 where $\bphi \in L_2(\R^2;\C^2)$. Let $\bphi = {\rm col}\, \{\phi_-,\phi_+\}$. Clearly, we have  $\u= {\rm col}\, \{u_{-,\eps},u_{+,\eps}\}$,
 where $u_{\pm,\eps}$ are the solutions of problems (16.14).

We can apply Theorem 14.7 to problem (16.17). The corresponding homogenized problem is
$$
i \frac{\partial \u_{0}(\x,\tau)}{\partial \tau} =  - \gamma \Delta \u_{0}(\x,\tau),\quad f_{\times,0} \u_{0}(\x,0) = \bphi(\x),
\eqno(16.18)
$$
where $\gamma$ is defined by (16.9), and the matrix $f_{\times,0}$ is defined by (16.16).
Clearly, we have $\u_0 = {\rm col}\, \{u_{-,0},u_{+,0}\}$,  where $u_{\pm,0}$ are the solutions of problems (16.15).

\smallskip\noindent\textbf{Proposition 16.4.}  \textit{Supose that the assumptions of Subsections} 16.1--16.3
\textit{are satisfied. Let $\u_{\eps}$ be the solution of problem} (16.17),  \textit{and let $\u_{0}$ be the solution of  problem} (16.18).
\textit{If  $\bphi \in H^s(\R^2;\C^2)$ with some $0 \le s \le 3$, then for $\tau \in \R$ and $\eps>0$ we have}
$$
\|  f_\times^\eps  \u_{\eps}(\cdot,\tau) -  f_{\times,0} \u_0 (\cdot,\tau) \|_{L_2(\R^2)} \le \eps^{s/3} {\mathfrak C}_\times(s;\tau) \|\bphi \|_{H^s(\R^2)}.
$$
\textit{The constant ${\mathfrak C}_\times(s;\tau)$ is given by} (13.21). \textit{If $\bphi \in L_2(\R^2;\C^2)$, then}
$$
\lim_{\eps \to 0}
\|  f_\times^\eps  \u_{\eps}(\cdot,\tau) -  f_{\times,0} \u_0 (\cdot,\tau) \|_{L_2(\R^2)} =0, \quad \tau \in \R.
$$

\smallskip
One can also apply the statement of Theorem~14.7 about the behavior of the solution for $|\tau| = \eps^{-\alpha}$ with $0<\alpha < 1$.
Also it is possible to consider  more general Cauchy problem for the nonhomogeneous equation and apply  Theorem~14.10.

\end{document}